\documentclass[final,onefignum,onetabnum]{siamart171218}

\usepackage{lipsum}
\usepackage{amsfonts}
\usepackage{graphicx}
\usepackage{grffile}
\usepackage{epstopdf}
\usepackage{algorithmic}
\ifpdf
  \DeclareGraphicsExtensions{.eps,.pdf,.png,.jpg}
\else
  \DeclareGraphicsExtensions{.eps}
\fi

\newsiamremark{remark}{Remark}
\newsiamremark{hypothesis}{Hypothesis}
\crefname{hypothesis}{Hypothesis}{Hypotheses}
\newsiamthm{claim}{Claim}

\usepackage{algorithm}
\usepackage{nicefrac}
\usepackage{algorithmic}
\usepackage{paralist}
\usepackage{mathtools}
\usepackage{tcolorbox}
\usepackage{amsmath,amsfonts,amssymb,amsbsy,amstext,amscd,amsxtra,multicol}
\usepackage{pifont}

\usepackage{forloop}
\usepackage{graphicx}
\usepackage{wrapfig}
\usepackage{tikz}
\usepackage{makecell} 

\def \NN {\mathbb N}

\def \EE {\mathbb E}

\def\cE{{\mathcal{E}}}

\newcommand{\argmin}{\mathop{\arg\!\min}}
\newcommand{\circledOne}{\text{\ding{172}}}
\newcommand{\circledTwo}{\text{\ding{173}}}
\newcommand{\circledThree}{\text{\ding{174}}}
\newcommand{\circledFour}{\text{\ding{175}}}

\usepackage[colorinlistoftodos,textsize=scriptsize]{todonotes}

\def\pd#1{{\color{magenta}#1}} 

\def\tf{\widetilde{f}}

\def\tnmf{\widetilde{\nabla}^m f}
\newcommand{\la}{\langle}
\newcommand{\ra}{\rangle}
\def\O{\widetilde{O}}
\def\e{\varepsilon}

\def \R {\mathbb R}
\def\RR{\mathcal R}

\def\eqdef{\overset{\text{eqdef}}{=}}

\headers{An Accel. Method for Der.-Free Smooth Stoch. Convex Optimization}{E. Gorbunov, P. Dvurechensky, and A. Gasnikov}

\title{An Accelerated Method for Derivative-Free Smooth Stochastic Convex Optimization\thanks{Submitted to the editors 30 April, 2019.
\funding{
The work of Eduard Gorbunov \newstuff{in Section~2.3.\ was} supported by the Ministry of Science and Higher Education of the Russian Federation (Goszadaniye) No. 075-00337-20-03, project No. 0714-2020-0005.
}}}

\author{Eduard Gorbunov\thanks{Moscow Institute of Physics and Technology; National Research University Higher School of Economics  
  (\email{eduard.gorbunov@phystech.edu}, \url{https://eduardgorbunov.github.io/}).}
\and Pavel Dvurechensky\thanks{Weierstrass Institute for Applied Analysis and Stochastics; Institute for Information Transmission Problems RAS; National Research University Higher School of Economics 
  (\email{pavel.dvurechensky@wias-berlin.de}).}
\and Alexander Gasnikov\thanks{Moscow Institute of Physics and Technology;  Institute for Information Transmission Problems RAS; National Research University Higher School of Economics 
    (\email{gasnikov@yandex.ru})}}

\usepackage{amsopn}

\ifpdf
\hypersetup{
  pdftitle={An Accelerated Method for Derivative-Free Smooth Stochastic Convex Optimization},
  pdfauthor={E. Gorbunov, P. Dvurechensky, and A. Gasnikov}
}
\fi

\def\pd#1{{\color{black}#1}} 
\def\pdd#1{{\color{black}#1}} 
\def\newstuff#1{{\color{black}#1}} 

\begin{document}

\maketitle

\begin{abstract}
  We consider an unconstrained problem of \pd{minimizing a} smooth convex function which is only available through noisy observations of its values, the noise consisting of two parts. 
Similar to stochastic optimization problems, the first part is of  stochastic nature. 
\pd{The} second part \pd{is additive noise of unknown nature, but bounded in absolute value}. 
In the two-point feedback setting, i.e. when pairs of function values are available, we propose an accelerated derivative-free algorithm together with its complexity analysis.
The complexity bound of our derivative-free algorithm is only by a factor of $\sqrt{n}$ larger than the bound for accelerated gradient-based algorithms, where $n$ is the dimension of the decision variable. 
We also propose a non-accelerated derivative-free algorithm with a complexity bound similar to the stochastic-gradient-based algorithm, that is, our bound does not have any dimension-dependent factor \pdd{except logarithmic}.
\pdd{Notably}, if the difference between the starting point and the solution is a sparse vector, for both our algorithms, we obtain \pd{a} better complexity bound if the algorithm uses \pd{an} $1$-norm proximal setup, rather than the Euclidean proximal setup, which is a standard choice for unconstrained problems
\end{abstract}



\section{Introduction}
Derivative-free or zeroth-order optimization \pd{\cite{rosenbrock1960automatic,fabian1967stochastic,brent1973algorithms,spall2003introduction,conn2009introduction}} is one of the oldest areas in optimization, which constantly attracts attention \pd{from the} learning community, mostly in connection to online learning in the bandit setup \cite{bubeck2012regret} and reinforcement learning \cite{salimans2017evolution,choromanski2018structured,fazel2018global,choromanski2018optimizing}. 
We study stochastic derivative-free optimization problems in a two-point feedback situation, considered by \cite{agarwal2010optimal,duchi2015optimal,shamir2017optimal} in the learning community and by \cite{nesterov2017random,stich2011linear,ghadimi2013stochastic,ghadimi2016mini-batch,gasnikov2016gradient-free} in the optimization community. 
Two-point setup allows \pd{one} to prove complexity bounds, which typically coincide with the complexity bounds for gradient-based algorithms up to a small-degree polynomial of $n$, where $n$ is the dimension of the decision variable. 
On the contrary, problems with one-point feedback are harder and complexity bounds for such problems either have worse dependence on $n$, or worse dependence on the desired accuracy of the solution, see \cite{nemirovsky1983problem,protasov1996algorithms,flaxman2005online, agarwal2011stochastic,jamieson2012query,shamir2013complexity,liang2014lzeroth-order,bach16highly-smooth,bubeck2017kernel-based} and the references therein.

More precisely, we consider the following optimization problem
\begin{equation}
\label{eq:PrSt}
   \min_{x\in\R^n} \left\{ f(x) := \EE_{\xi}[F(x,\xi)] = \int_{\mathcal{X}}F(x,\xi)dP(x) \right\}, 
\end{equation}
where $\xi$ is a random vector with probability distribution $P(\xi)$, $\xi \in \mathcal{X}$, and the function $f(x)$ is closed  and  convex. Note that $F(x,\xi)$ can be non-convex \pd{in} $x$ with positive probability. Moreover, we assume that, \pd{almost sure w.r.t. distribution $P$}, the function $F(x,\xi)$ has gradient $g(x,\xi)$, which is $L(\xi)$-Lipschitz continuous with respect to the Euclidean norm. \pd{We assume that we know a constant $L_2\geqslant 0$ such that $\sqrt{\EE_{\xi} L(\xi)^2 } \leq L_2 < +\infty$.} Under \pd{these} assumptions, $\EE_{\xi}g(x,\xi) = \nabla f(x)$ and $f$ \pd{is $L_2$-smooth, i.e.} has $L_2$-Lipschitz continuous gradient with respect to the Euclidean norm. Also we assume that\pd{, for all $x$,}
\begin{equation}
\label{stoch_assumption_on_variance}
    \EE_{\xi}[\|g(x,\xi) - \nabla f(x)\|_2^2] \leqslant \sigma^2,
\end{equation}
where $\|\cdot\|_2$ is the Euclidean norm. We emphasize that, unlike \cite{duchi2015optimal}, we \textit{do not assume that $\EE_\xi\left[\|g(x,\xi)\|_2^2\right]$ is bounded} since it is not the case for many unconstrained optimization problems, e.g. for deterministic quadratic optimization problems. 

Finally, we assume that \pd{we are in the two-point feedback setup, which is also connected to the common  random  numbers assumption, see \cite{larson2019derivative-free} and references therein. Specifically,} an optimization procedure, given a pair of points $(x,y) \in \R^{2n}$, can obtain a pair of noisy stochastic realizations $(\tf(x,\xi),\tf(y,\xi))$ of the objective value $f$, which we refer to as \textit{oracle call}. Here
\begin{equation}
\label{eq:tf_def}
\tf(x,\xi) = F(x,\xi) + \eta(x,\xi), \quad |\eta(x,\xi)| \leqslant  \Delta, \;\forall x \in \R^n, \; \text{a.s.\ in } \xi,
\end{equation}
and \pd{there is a possibility to obtain an iid sample $\xi$ from $P$}. This makes our problem more complicated than problems studied in the literature. Not only \pd{do} we have  stochastic noise in problem \eqref{eq:PrSt}, but \pd{also} an additional noise $\eta(x,\xi)$, which can be adversarial. 

\pdd{O}ur model of the two-point feedback oracle is pretty general and covers deterministic exact oracle \pdd{or} even specific types of one-point feedback oracle. For example, if the function $F(x,\xi)$ is separable, i.e.\ $F(x,\xi) = f(x) + h(\xi)$, where $\EE_\xi\left[h(\xi)\right] = 0$, $|h(\xi)|\le \tfrac{\Delta}{2}$ for all $\xi$ and the oracle gives us $F(x,\xi)$ \pd{at a} given point $x$, then for all $\xi_1,\xi_2$ we can define $\widetilde{f}(x,\xi_1) = F(x,\xi_1)$ and $\widetilde{f}(y,\xi_2) = F(y,\xi_2) = F(y,\xi_1) + h(\xi_2) - h(\xi_1)$. Since $|h(\xi_2)-h(\xi_1)| \le |h(\xi_2)| + |h(\xi_1)| \le \Delta$ we can use representation \eqref{eq:tf_def} omitting dependence of $\eta(x,\xi_1)$ on $\xi_2$.
Moreover, such \pd{an} oracle can be \pd{encountered} in practice as rounding errors can be \pd{modeled as a process of adding a random bit modulo $2$ to the last or several last bits in machine number representation format (see \cite{gasnikov2016stochastic} for details).}


As it is known \cite{lan2012optimal,devolder2011stochastic,dvurechensky2016stochastic,gasnikov2016stochasticInter}, if \pdd{a} stochastic approximation $g(x,\xi)$ for the gradient of $f$ is available, an accelerated gradient method has oracle complexity bound (i.e. \pdd{the total} number of \pdd{stochastic} first-order oracle calls) $O\left(\max\left\{\sqrt{L_2R_2^2/\e},\sigma^2R_2^2/\e^2\right\}\right)$, where $\e$ is the target optimization error \pd{in terms of the objective residual}, the goal being to find such $\hat{x}$ that $\EE f(\hat{x}) - f^* \leq \e$. Here $f^*$ is the global optimal value of $f$\pd{, $R_2$ is such that $\|x_0-x^*\|_2 \leq R_2$ with $x^*$ being \textit{some} solution.}
The question, to which we give a positive answer in this paper, is as follows.

\textit{Is it possible to solve a stochastic optimization problem with the same $\e$-dependence in the iteration and sample complexity and only noisy observations of the objective value?} 

\pd{Many existing first- and zero-order methods are based on so-called proximal setup (see \cite{ben-tal2020lectures} and Subsection \ref{sec:preliminaries} for the precise definition). This includes a choice of some norm in $\R^n$ and a corresponding prox-function, which is strongly convex with respect to this norm. Standard gradient method for unconstrained problems such as \eqref{eq:PrSt} is obtained when one chooses the Euclidean $\|\cdot\|_2$-norm as the norm and squared Euclidean norm as the prox-function. We go beyond this conventional path and consider $\|\cdot\|_1$-norm in $\R^n$ and corresponding prox-function given in \cite{ben-tal2020lectures}. Yet this proximal setup is described in the textbook, we are not aware of any particular examples where it is used for unconstrained optimization problems. Notably, as we show in our analysis, this choice can lead to better complexity bounds. In what follows, we characterize these two cases by the choice of $\|\cdot\|_p$-norm with $p\in \{1,2\}$ and its conjugate $q \in \{2,\infty\}$, given by the identity $\tfrac{1}{p}+\tfrac{1}{q} = 1$.}

\subsection{Related Work}
\pd{Online optimization with two-point bandit feedback was considered in \cite{agarwal2010optimal}, where regret bounds were obtained. 
}
Non-smooth deterministic and stochastic problems in the two-point derivative-free \pdd{offline optimization} setting \pd{were} considered \pdd{in} \cite{nesterov2017random}.\footnote{We list the references in the order of the date of the first appearance, but not in the order of the date of official publication.} 
Non-smooth stochastic problems were considered \pdd{in} \cite{shamir2017optimal} and independently \pdd{in} \cite{bayandina2017gradient-free}, the latter paper considering also problems with additional noise of an unknown nature in the objective value.
\pdd{The authors of} \cite{duchi2015optimal} consider smooth stochastic optimization problems, yet under additional quite restrictive assumption $\EE_\xi\left[\|g(x,\xi)\|_{\pdd{q}}^2\right] \pd{<} +\infty$. Their bound was improved \pdd{in} \cite{gasnikov2016gradient-free,gasnikov2017stochastic} for the problems with non-Euclidean proximal setup and noise in the objective value.  Strongly convex problems with different smoothness assumptions were considered \pdd{in} \cite{gasnikov2017stochastic,bayandina2017gradient-free}. 
Smooth stochastic convex optimization problems, without the assumption that $\EE \|g(x,\xi)\|_2^2 < +\infty$, were studied \pdd{in} \cite{ghadimi2016mini-batch,ghadimi2013stochastic} for the Euclidean case.
Accelerated and non-accelerated derivative-free method \pdd{for smooth but} deterministic problems \pd{were} proposed in \cite{nesterov2017random} and extended in \cite{bogolubsky2016learning,dvurechensky2017randomized} for the case of additional bounded noise in the function value. 

Table~\ref{tbl:summ_and_comp} presents a detailed comparison of our results and most close results in the literature on two-point feedback derivative-free optimization and assumptions, under which they are obtained.  
\pd{
The first row corresponds to the non-smooth setting with the assumption that $\EE_\xi\left[\|g(x,\xi)\|_2^2\right] \leq M_2^2$, which mostly restricts the scope to constrained optimization problems on a convex set with the diameter $R_p$ measured by $\|\cdot\|_p$-norm. This setting is very well understood with the proposed methods being able to solve stochastic optimization problems with additional bounded noise in the objective value and to use non-Euclidean proximal setup. Importantly, non-Euclidean proximal setup corresponding to $p=1$, $q=\infty$ allows one to obtain a complexity bound with only logarithmic dependence on the dimension $n$. 
}

\pd{Rows 2-6 of Table~\ref{tbl:summ_and_comp} correspond to smooth problems with $L_2$-Lipschitz continuous gradient, which makes possible to apply Nesterov's acceleration and obtain better complexity bounds. In this case stochastic optimization problems are characterized by the variance $\sigma^2$ of the stochastic gradient, see \eqref{stoch_assumption_on_variance}.
For the smooth setting the full picture is not completely understood in the literature, and our goal is to obtain methods, which provide the full picture similarly to the non-smooth setting by combining stochastic optimization setup, additional bounded noise in the objective value, acceleration, and better complexity bounds achievable owing to the use of non-Euclidean proximal setup corresponding to $p=1$, $q=\infty$. Previous works for the smooth case consider only Euclidean case and either deterministic problems with additional bounded noise \cite{nesterov2017random,bogolubsky2016learning,dvurechensky2017randomized} or stochastic problems without additional bounded noise \cite{ghadimi2016mini-batch,ghadimi2013stochastic}.
}

{
\renewcommand{\arraystretch}{1.7}
\begin{table}
\small
\begin{center}
\begin{tabular}{|c|c|c|c|c|c|}
\hline
Method & \begin{tabular}{c}Assumptions \end{tabular} & Oracle complexity, $\O\left( \cdot \right)$ & $p=1$ & $\sigma^2$ & \pdd{$\Delta$}\\
\hline
\makecell{MD \\ {\tiny\cite{duchi2015optimal,gasnikov2016gradient-free, gasnikov2017stochastic,gasnikov2016stochastic}}\\
{\tiny\cite{shamir2017optimal,bayandina2017gradient-free}}}
&bound. gr. &$\tfrac{{\color{blue}n^{\tfrac{2}{q}}}{\color{red}M_2^2}R_p^2}{\varepsilon^2}$ & $\surd$ & $\surd$ & $\surd$ \\
\hline
\makecell{RSGF\\ {\tiny\cite{ghadimi2016mini-batch,ghadimi2013stochastic}}} & bound. var. 
& $\max\left\{\tfrac{{\color{red}n}L_2R_2^2}{\e},\tfrac{{\color{red}n}{\color{blue}\sigma^2}R_2^2}{\e^2}\right\}$ & $\times$ & $\surd$ & $\times$ \\
\hline
\makecell{RS \\ {\tiny\cite{nesterov2017random,bogolubsky2016learning}} } & 
& $\tfrac{{\color{red}n}L_2R_2^2}{\e}$& $\times$ & $\times$ & $\surd$  \\
\hline
\makecell{RDFDS\\  {[This paper]}} & bound. var. 
& $\max\left\{\tfrac{{\color{blue}n^{\tfrac{2}{q}}}L_2R_p^2}{\e},\tfrac{{\color{blue}n^{\tfrac{2}{q}}}{\color{blue}\sigma^2}R_p^2}{\e^2}\right\}$ & $\surd$ & $\surd$ & $\surd$  \\
\hline
\makecell{
AccRS\\ {\tiny\cite{nesterov2017random,dvurechensky2017randomized}}
}
& 
& ${\color{red}n}\sqrt{\tfrac{L_2R_2^2}{\varepsilon}}$ & $\times$ & $\times$ & $\surd$  \\
\hline
\makecell{ARDFDS \\ {[This paper]}} & bound. var. 
& $\max\left\{{\color{blue}n^{\tfrac{1}{2}+\tfrac{1}{q}}}\sqrt{\tfrac{L_2R_p^2}{\e}},\tfrac{{\color{blue}n^{\tfrac{2}{q}}}{\color{blue}\sigma^2}R_p^2}{\e^2}\right\}$ & $\surd$ & $\surd$ & $\surd$   \\
\hline
\end{tabular}
\caption{Comparison of oracle complexity (total number of \pd{zero-order} oracle calls) of different methods with two-point feedback 
for convex optimization problems. \pd{$R_p$ is such that $\|x_0-x^*\|_p \leq R_p$ with $x^*$ being \textit{some} solution.} In the column ``Assumptions'' we use ``bound. gr.'' \pdd{for} $\EE_\xi\left[\|g(x,\xi)\|_2^2\right] \le M_2^2$ and ``bound. var.'' \pdd{for} $\EE_{\pd{\xi}}\|g(x,\xi) - \nabla f(x)\|_2^2 \leqslant \sigma^2$. Column ``$p=1$'' corresponds to the support of non-Euclidean \pdd{proximal} setup, column ``$\sigma^2$'' to the support of stochastic optimization \pd{problems}, ``\pdd{$\Delta$}'' corresponds to the support of additional \pdd{bounded} noise of unknown nature. \pd{All the rows except the first one assume that $f$ is $L_2$-smooth. $\widetilde{O}(\cdot)$ notation means ${O}(\cdot)$ up to logarithmic factors in $n,\e$.}
\label{tbl:summ_and_comp}
}
\end{center}
\end{table}

}

We also mention the works \pd{\cite{nemirovsky1983problem,protasov1996algorithms,polyak1990optimal,dippon2003accelerated,flaxman2005online,saha2011improved,dekel2015bandit,gasnikov2017stochastic,agarwal2011stochastic,liang2014lzeroth-order,belloni2015escaping,bubeck2017kernel-based,shamir2013complexity,jamieson2012query,hazan2014bandit,bach16highly-smooth,jamieson2012query,bartlett2019simple,locatelli2018adaptivity,akhavan2020exploiting} where derivative-free optimization with one-point feedback is studied }in different settings, and works \cite{nesterov2005smooth,allen2014linear} on coupling non-accelerated methods to obtain acceleration, which inspired our work. \pd{After our paper appeared as a preprint, the papers \cite{berahas2019derivative-free,bollapragada2019adaptive} studied derivative-free quasi-Newton methods for problems with noisy function values, and the paper \cite{berahas2019theoretical} reported theoretical and empirical comparison of different gradient approximations for zero-order methods.} 
\pd{The authors of \cite{chen2020accelerated} combine accelerated derivative-free optimization with accelerated variance reduction technique for finite-sum convex problems.}
\pd{For a recent review of derivative-free optimization see \cite{larson2019derivative-free}.}
\pdd{We extend the proposed algorithms for a more general setting of inexact directional derivative oracle as well as for strongly convex problems in \cite{dvurechensky2018accelerated2}.}
\pdd{Mixed first-order/zero-order setting is considered in \cite{beznosikov2020derivative-Free} and zero-order methods for non-smooth saddle-point problems are developed in \cite{beznosikov2020gradient-free}.}


\subsection{Our Contributions}\label{sec:our_contrib}
As our main contribution, we propose an accelerated method for smooth stochastic derivative-free optimization \pdd{with two-point feedback}, which we call Accelerated Randomized Derivative-Free Directional Search (ARDFDS). Our method has the complexity bound 
\begin{equation}
\label{eq:ARDFDSComplInformal}
\O\left(\max\left\{n^{\tfrac12+\tfrac{1}{q}}\sqrt{\tfrac{L_2R_p^2}{\e}},\tfrac{n^{\tfrac{2}{q}}\sigma^2R_p^2}{\e^2}\right\}\right),
\end{equation}
where \pd{$\widetilde{O}$ hides logarithmic factor of the dimension}, \pdd{$R_p$ is such that $\|x_0-x^*\|_p \leq R_p$ with $x^*$ being an \textit{arbitrary} solution to \eqref{eq:PrSt} and $x_0$ being the starting point of the algorithm.}
\pd{We underline that our bounds hold for any solution. Thus, to obtain the best possible bound, one can consider the closest solution to the starting point.}
In the Euclidean case $p=q=2$, the first term in the above bound has better dependence \pd{on $\e$},  $L_2$ and $R_2$ than the \pd{first term in the} bound in \cite{ghadimi2016mini-batch,ghadimi2013stochastic}. Unlike these papers, our bound also covers the non-\pdd{E}uclidean case $p=1$, $q=\infty$ and, due to that, allows to obtain better complexity bounds. To illustrate this, let us \pd{consider an arbitrary solution $x^*$ to \eqref{eq:PrSt}}, start method from a point $x_0$ and define the sparsity $s$ of the vector $x_0 - x^*$, i.e. $\|x_0-x^*\|_1 \leq s \cdot \|x_0-x^*\|_2$ and $1 \leq s \leq \sqrt{n}$. Then the complexity of our method for $p=1$, $q=\infty$ is $\O\left(\max\left\{\sqrt{\tfrac{ns^2L_2\|x_0-x^*\|_2^2}{\e}},\tfrac{s^2\sigma^2\|x_0-x^*\|_2^2}{\e^2}\right\}\right)$, which is \textit{always} no worse than the complexity for $p=q=2$, \pd{which is} $\O\left(\max\left\{\sqrt{\tfrac{n^2L_2\|x_0-x^*\|_2^2}{\e}},\tfrac{n\sigma^2\|x_0-x^*\|_2^2}{\e^2}\right\}\right)$ and allows to gain up to $\sqrt{n}$ if $s$ is close to 1. Notably, this is done automatically, without any prior knowledge of $s$. \pd{An example of this situation can be a typical compressed sensing problem \cite{candes2006stable,donoho2006compressed} of recovering a sparse signal $x^*$ from noisy observations of a linear transform of $x^*$ via solving an optimization problem. In this case, if $x_0=0$ then $x_0-x^*$ is sparse by the problem assumption. Moreover, since our bounds hold for arbitrary solution $x^*$, to get better complexity estimate, one can use the bound obtained using the sparsest solution.}

Unlike \pd{previous works}, we consider additional, possibly adversarial noise $\eta(x,\xi)$ in the objective value and analyze how this noise affects the convergence rate estimates. \pdd{If the noise can be controlled and $\Delta$ can be made arbitrarily small, e.g. if the objective is calculated by an auxiliary procedure, we estimate how $\Delta$ should depend on the target accuracy $\e$ to ensure finding an $\e$-solution.} \pd{If the noise is uncontrolled, e.g. we only have an estimate for the noise level $\Delta$} and we cannot \pd{make $\Delta$ arbitrarily small}, we can run our algorithms \pd{and guarantee that they generate a point with expected objective residual bounded by a quantity dependent on $\Delta$}.
This is important when the objective is given as a solution to some auxiliary problem\pd{,} which can \pd{not} be solved exactly, e.g. in bi-level optimization or reinforcement learning. 
\pd{It should also be mentioned that our assumption $\EE_{\xi} [L(\xi)^2] \leq L_2$ for some $L_2$ is weaker than the assumption that there is $L_2$ s.t. $L(\xi) \leq L_2$ a.s. in $\xi$, which is used in \cite{ghadimi2016mini-batch,ghadimi2013stochastic}.} 

As our second contribution, we propose a non-accelerated Randomized Derivative-Free Directional Search (RDFDS) method with the complexity bound
\begin{equation}
\label{eq:RDFDSComplInformal}
\O\left(\max\left\{\tfrac{n^{\tfrac{2}{q}}L_2R_p^2}{\e},\tfrac{n^{\tfrac{2}{q}}\sigma^2R_p^2}{\e^2}\right\}\right),
\end{equation}
where, unlike \cite{ghadimi2016mini-batch,ghadimi2013stochastic}, the non-\pd{E}uclidean case $p=1$, $q=\infty$ \pdd{is also covered} with the gain in the complexity \pdd{of} up to the factor of $n$ \pdd{in comparison to the case $p=q=2$}. \pdd{Notably, in the non-\pd{E}uclidean case}, we obtain a nearly dimension independent \pd{($\widetilde{O}$ hides logarithmic factor of the dimension)} complexity bound despite we use only noisy function value observations.

\textbf{Why \pd{is it} important to improve the first term \pdd{under} the maximum?}
\begin{enumerate}
    \item\textbf{Acceleration when $n$ is \pdd{large}.} The first term \pdd{under the maximum} dominates the second term when $\sigma^2 \le \tfrac{\e^{\tfrac{3}{2}}n^{\tfrac{1}{2}-\tfrac{1}{q}}\sqrt{L_2}}{R_p}$ in the accelerated case and \pdd{when} $\sigma^2 \le \varepsilon L_2$ in the non-accelerated case, which could be met in practice if $\e$, $L_2$ and $n$ are \pdd{large} enough \pdd{compared to} $R_p$. \pdd{For example, if 
    $p=1,q=\infty$ and we would like} to \pdd{find} \pdd{an} $\varepsilon$-solution with $\varepsilon = 10^{-3}$ and $L_2 = 100$, $R_p = 10$, $n = 10 000$ (or \pdd{larger}), \pdd{and the variance satisfies mild assumption $\sigma^2 \le 10^{-1}$, then the complexity of ARDFDS is better than that of RDFDS.}
    \item\textbf{Better dimension dependence in the deterministic case.} We underline that even in the deterministic case \pdd{with} $\sigma = 0$ \pdd{and} without additive noise, both our non-accelerated and accelerated complexity bounds 
    for $p=1$ are new. Moreover, disregarding $\ln n$ factors, \pdd{for $s\in [1,\sqrt{n}]$,} the existing bounds \cite{nesterov2017random} are $n/s^2$ and $n/s$ times worse than our new bounds respectively in non-accelerated and accelerated cases. Importantly, in the non-accelerated case our bound is dimension-independent up to a $\ln n$ factor.
    \item\textbf{Parallel computation of mini-batches makes acceleration reasonable when $\sigma^2$ is not small.} 
    \pd{
    Even when the second term in  \eqref{eq:ARDFDSComplInformal} is dominating and, thus, the total computation time is proportional to the second term, using parallel computations we can force the total computation time to be proportional to the first term, underlining the importance of making it smaller via accelerating the method. The idea is to use parallel computations of mini-batches as follows.
    Instead of sampling one $\xi$ in each iteration of the algorithm 
    one can consider a mini-batch of size $r$, i.e. sample $r$ iid realizations of $\xi$ and average $r$ finite-difference approximations for the gradient to reduce the variance of this approximation from $\sigma^2$ to $\tfrac{\sigma^2}{r}$. 
    If one can have an access to at least $r$ processors, in each iteration all processors simultaneously in parallel can make a call to the zeroth-order oracle and calculate finite-difference approximation for the gradient. Then a processor chosen to be central calculates the average of these $r$ approximations, which gives a mini-batch approximation of the gradient. Since this work is done in parallel, it takes nearly the same amount of time as using a mini-batch of size $1$ in the standard approach. By choosing sufficiently large $r$, one can make the second term  in  \eqref{eq:ARDFDSComplInformal} (which is now proportional to $\tfrac{\sigma^2}{r}$) smaller than the first term.
    Hence, the total computation time will be proportional  to the first term under the maximum in \eqref{eq:ARDFDSComplInformal}. Such an acceleration can be achieved by a reasonable amount of processors. For example, if $\sigma^2=1$, which is not small, $n=10000$, $R_p = 10$, $\varepsilon = 10^{-3}$ and $L_2 = 100$, then it is sufficient to have $r=10^{2.5}\approx 316$ processors which is \pd{a small number compared to} modern supercomputers and clusters that often have $\sim 10^5-10^6$ processors.
    }
\end{enumerate}

\section{Algorithms for Stochastic Convex Optimization}
\subsection{Preliminaries}\label{sec:preliminaries}

\noindent \textbf{\pd{$p$-norm p}roximal setup.}
Let $p\in[1,2]$ and $\|x\|_p$ be the \pdd{$\|\cdot\|_p$}-norm in $\R^n$ defined as $\|x\|_p^p = \sum\limits_{i=1}^n|x_i|^p$. \pdd{Further, let} $\|\cdot\|_{q}$ be its dual, defined by $\|g\|_{q} = \max\limits_{x} \big\{ \la g, x \ra, \| x \|_p \leq 1 \big\}$, where $q \in [2,\infty]$ is the conjugate number to $p$, given by $\tfrac{1}{p} + \tfrac{1}{q} = 1$, and, for $q = \infty$, \pdd{by definition,} $\|x\|_\infty = \max\limits_{i=1,\ldots,n}|x_i|$. 
\newstuff{We also use $\|x\|_0$ to denote the number of non-zero components of $x\in\R^n$.} We choose a \textit{prox-function} $d(x)$, which is continuous
and $1$-strongly convex on $\R^n$ with respect to $\|\cdot\|_p$, i.e., for any $x, y \in \R^n$, $d(y)-d(x) -\la \nabla d(x) ,y-x \ra \geq \tfrac12\|y-x\|_p^2$. 
Without loss of generality, we assume that $\min\limits_{x\in \R^n} d(x) = 0$.
We define also the corresponding \textit{Bregman divergence} $V[z] (x) = d(x) - d(z) - \la \nabla d(z), x - z \ra$, \pd{for} $x, z \in \R^n$. Note that, by the \pd{$1$-}strong convexity of $d(\cdot)$,
\begin{equation}
\label{eq:VStrConv}
V[z] (x) \geq \tfrac{1}{2}\|x-z\|_p^2, \quad \pd{\forall} \; x, z \in \R^n.
\end{equation}
For $p=1$, we choose the prox-function (see \cite{ben-tal2020lectures}) $d(x) = \tfrac{\|x\|_\kappa^2{\pd{\exp(1)}}n^{(\kappa-1)(2-\kappa)/\kappa}\ln n}{2}$, where $\kappa=1 + \tfrac{1}{\ln n}$ and, for the case $p=2$, we choose the prox-function to be 
$d(x) = \tfrac{1}{2}\|x\|_2^2.$ 

\noindent \textbf{Main technical lemma.}
In our proofs of complexity bounds, we rely on the following lemma. The proof is rather technical and is provided in the appendix. 
\begin{lemma}
		\label{Lm:MainTechLM}
    Let $e \in RS_2(1)$, i.e. be a random vector uniformly distributed on the surface of the unit Euclidean sphere in $\R^n$, $p\in[1,2]$ and $q$ be given by $\tfrac{1}{p}+\tfrac{1}{q} = 1$. \pd{Define  $\rho_n = \min\{q-1,\,16\ln n - 8\}n^{2/q-1}$.} 
    Then, for $n \geqslant8$,
    \begin{align}\label{assumption_e_norm}
        \EE_e\|e\|_q^2 &\leq \rho_n, \; \text{\pd{and}} \\
    \label{assumption_e_product}
        \EE_e\left(\la s,e \ra^2 \|e\|_q^2 \right) &\leq \tfrac{6\rho_n}{n} \|s\|_2^2,     \quad \forall s \in \R^n.
    \end{align}
\end{lemma}

\noindent \textbf{Stochastic approximation of the gradient.} 
Based on the noisy observations \eqref{eq:tf_def} of the objective value, we form the following stochastic approximation of $\nabla f(x)$
\begin{equation}
\label{eq:FinDiffStocGrad}
\tnmf^t(x) = \frac{1}{m}\sum\limits_{i=1}^m\frac{\tf(x+te,\xi_i)-\tf(x,\xi_i)}{t}e,
\end{equation}
where $e \in RS_2(1)$, $\xi_i$, $i=1,...,m$ are independent realizations of $\xi$, $m$ is the \textit{mini-batch size}, $t$ is some small positive parameter, which we call \textit{smoothing parameter}.

\subsection{\newstuff{Accelerated Randomized Derivative-Free Directional Search}}
The method is listed as Algorithm~\ref{Alg:ARDFDS}. \pd{Following \cite{nesterov2017random,ghadimi2016mini-batch,ghadimi2013stochastic}  we assume that $L_2$ is known.} 
\pd{The possible choice of the smoothing parameter $t$ and mini-batch size $m$ are discussed below.}
\pd{Note that at every iteration the algorithm requires to solve an auxiliary minimization problem. As it is shown in \cite{ben-tal2020lectures}, for both cases $p=1$ and $p=2$ this minimization can be made explicitly in $O(n)$ arithmetic operations.}
\begin{algorithm}
	\caption{Accelerated Randomized Derivative-Free Directional Search (ARDFDS)}
	\label{Alg:ARDFDS}
	\begin{algorithmic}[1]
		\REQUIRE $x_0$~-- starting point; $N$~-- number of iterations; \pd{$L_2$~-- smoothness constant}; $m \geq 1$~-- mini-batch size; $t > 0$~-- smoothing parameter\pd{;} 
		\pd{$V(\cdot,\cdot)$~-- Bregman divergence}.
		\STATE $y_0 \leftarrow x_0, \, z_0 \leftarrow x_0$.
		\FOR{$k=0,\, \dots, \, N-1$.}
		\STATE Generate $e_{k+1} \in RS_2\left( 1 \right)$ independently from previous iterations and $\xi_i^{\pd{k+1}}$, $i=1,...,m$ -- independent realizations of $\xi$\pd{, which are also independent from previous iterations}. 
		\STATE $\tau_k \leftarrow \tfrac{2}{k+2}$, $x_{k+1} \leftarrow \tau_kz_k + (1-\tau_k)y_k $.
		\STATE Calculate $\tnmf^t(x_{k+1})$ \pd{using} \eqref{eq:FinDiffStocGrad} \pd{with $e=e_{k+1}$} and set $y_{k+1} \leftarrow x_{k+1}-\tfrac{1}{2L_2}\tnmf^t(x_{k+1})$. 
		\STATE $\alpha_{k} \leftarrow \tfrac{k+1}{96n^2\rho_n L_2}$, $z_{k+1} \leftarrow \argmin\limits_{z\in\R^n} \Big\{ \alpha_{k+1} n \left\langle \tnmf^t(x_{k+1}), \, z-z_k \right\rangle+V[z_{k}] \left( z \right)\Big\}.$
		\ENDFOR
		\RETURN $y_N$
	\end{algorithmic}
\end{algorithm}
%
\begin{theorem}
\label{Th:ARDFDSConv}
    Let ARDFDS be applied to solve problem \eqref{eq:PrSt}, 
    \pd{ $x^*$ be an arbitrary solution to \eqref{eq:PrSt} and $\Theta_p$ be such that $ V[z_0](x^*)\leq \Theta_p$.} 
    Then, for all $n\ge 8$,
    \begin{equation}\label{eq:ARDFDSConv}
        \begin{array}{rcll}
            \EE[f(y_N)] - f(x^*) 
        &\leqslant& \tfrac{384 n^2\rho_nL_2\Theta_p}{N^2} +\tfrac{384N}{nL_2}\tfrac{\sigma^2}{m} + \tfrac{12\sqrt{2n\Theta_p}}{N^2}\left( \tfrac{L_2t}{2} + \tfrac{2\Delta }{t} \right)&\\ &&\quad+ 
				 \tfrac{6N}{L_2} \left( L_2^2t^2+\tfrac{16\Delta^2}{t^2} \right)+ \tfrac{N^2}{24n\rho_nL_2} \left( L_2^2t^2 + \tfrac{16\Delta^2 }{t^2} \right).
        \end{array}
    \end{equation}
\end{theorem}

\pd{\noindent Before we prove Theorem \ref{Th:ARDFDSConv} in the next subsection, let us discuss its result. 
In the simple case $\Delta=0$, all the terms in the r.h.s. of \eqref{eq:ARDFDSConv} can be made smaller than $\e$ for any $\e \geqslant 0 $ by an appropriate choice of $N,m,t$. 
Thus, we consider a more interesting case and assume that the noise level satisfies $0< \Delta \leqslant L_2\Theta_p n^3\rho_n^2/2$. The second inequality is non-restrictive since by the $L_2$-smoothness $f(x_0)-f^*\leqslant L_2\Theta_p$ and it is natural to assume that the oracle error is smaller than the initial objective residual. In order to minimize the terms with $L_2^2t^2 + \tfrac{16\Delta^2}{t^2}$ in the r.h.s of \eqref{eq:ARDFDSConv}, we set $t$ as $2\sqrt{\tfrac{\Delta}{L_2}}$. Substituting this into the r.h.s. of  \eqref{eq:ARDFDSConv} and using that, by our assumption on $\Delta$,  $\Theta_p n^2\rho_nL_2 \geqslant \sqrt{2nL_2\Theta_p \Delta}$, we obtain the following inequality
\begin{equation}
\label{eq:ARDFDSConv_cor}
    \begin{array}{rl}
        \EE[f(y_N)] - f(x^*) 
        &\leqslant 
		\tfrac{408L_2\Theta_p n^2\rho_n}{N^2} + \tfrac{384N}{nL_2}\tfrac{\sigma^2}{m} + 48N\Delta + \tfrac{N^2\Delta}{3n\rho_n}
    \end{array}
\end{equation}
First, we consider the situation of controlled noise level $\Delta$ which can be made arbitrarily small. For example, the value of $f$ is defined as a solution of some auxiliary problem, which can be solved numerically with arbitrarily small accuracy $\Delta$. Then we have control over parameters $N,m,\Delta$ in the r.h.s of \eqref{eq:ARDFDSConv_cor} and can choose these parameters to make it smaller than $\varepsilon$. 
First, we choose $N$ to make the first term to be smaller than $\e$. After that we choose $m$ to make the second term smaller than $\e$. Finally, we choose $\Delta$ to make all the other terms smaller than $\e$. The resulting values of these parameters up to constants are given in Table \ref{tbl:ARDFDS_params}. As a summary, we have the following corollary of Theorem \ref{Th:ARDFDSConv}.}
\begin{corollary}\label{cor:ardfds_conv}
    \pd{
    Assume that the value of $\Delta$ can be controlled and satisfies $0 < \Delta \leqslant L_2\Theta_p n^3\rho_n^2/2$.
    Assume that for a given accuracy $\e \geqslant 0$ the values of the parameters $N(\e),m(\e),t(\e),\Delta(\e)$ satisfy relations stated in Table \ref{tbl:ARDFDS_params} and ARDFDS 
    is applied to solve problem \eqref{eq:PrSt}. 
    Then the output point $y_N$ satisfies $\EE\left[f(y_N)\right] - f(x^*) \leqslant \e$. Moreover, the overall number of oracle calls is $N(\e)m(\e)$ given in the same table.} 
\end{corollary}
{
{\small
\begin{table*}[t]
\centering
    \centering
    \begin{tabular}{|c|c|c|}
        \hline
         & $p=1$ & $p=2$ \\
        \hline
        $N(\e)$ & $\sqrt{\tfrac{ n\ln nL_2 \Theta_1}{\varepsilon}}$ & $\sqrt{\tfrac{ n^{2}L_2\Theta_2}{\varepsilon}}$ \\
        \hline
        $m(\e)$ & $\max\left\{1, \tfrac{\sigma^2}{\varepsilon^{3/2}}\cdot\sqrt{\tfrac{\Theta_1\ln n}{nL_2}}\right\}$ & $\max\left\{1,\tfrac{\sigma^2}{\varepsilon^{3/2}}\cdot\sqrt{\tfrac{\Theta_2}{L_2}}\right\}$ \\ 
		\hline	
        $\Delta(\e)$ & $\min\left\{\tfrac{\varepsilon^{3/2}}{\sqrt{ L_2\Theta_1n\ln n}},\, \tfrac{\varepsilon^2}{nL_2\Theta_1}\right\}$ & $\min\left\{\tfrac{\varepsilon^{3/2}}{n\sqrt{L_2\Theta_2}},\,  \tfrac{\varepsilon^2}{nL_2\Theta_2}\right\}$\\
        \hline	
		$t(\e)$ & $\min\left\{\tfrac{\varepsilon^{3/4}}{\sqrt[4]{L_2^3\Theta_1 n\ln n}},\, \tfrac{\varepsilon}{L_2\sqrt{n\Theta_1}}\right\}$ & $\min\left\{\tfrac{\varepsilon^{3/4}}{\sqrt[4]{n^2L_2^3\Theta_2}},\,  \tfrac{\varepsilon}{L_2\sqrt{n\Theta_2}}\right\}$ \\
        \hline
        \pd{$N(\e)m(\e)$} & $\max\left\{\sqrt{\tfrac{ n\ln nL_2\Theta_1}{\varepsilon}}, \tfrac{\sigma^2\Theta_1 \ln n}{\varepsilon^2}\right\}$ & $\max\left\{\sqrt{\tfrac{ n^{2}L_2\Theta_2}{\varepsilon}}, \tfrac{\sigma^2\Theta_2 n}{\varepsilon^2}\right\}$\\
        \hline
    \end{tabular}
    \caption{\pd{Summary of the values for $N,m,\Delta,t$ and the total number of function value evaluations $Nm$ guaranteeing for the cases $p=1$ and $p=2$ that Algorithm~\ref{Alg:ARDFDS} outputs $y_N$ satisfying $\EE\left[f(y_N)\right] - f(x_*) \le \varepsilon$. Numerical constants are omitted for simplicity.} \label{tbl:ARDFDS_params}
     }
\end{table*}}
}

\pd{Note that in the case of uncontrolled noise level $\Delta$, the values of this parameter stated in Table \ref{tbl:ARDFDS_params} can be seen as the maximum value of the noise level which can be tolerated by the method still allowing it to achieve $\EE\left[f(y_N)\right] - f(x^*) \leqslant \e$.} 

{\color{black}
Next, we consider the case of uncontrolled noise level $\Delta$ and estimate the smallest expected objective residual which can be guaranteed in theory. First, we focus on the following three terms in the r.h.s. of \eqref{eq:ARDFDSConv_cor}, for simplicity disregarding the numerical constants,
\begin{equation}
    \label{eq:smallest_Delta_rhs}
    \tfrac{L_2\Theta_p n^2\rho_n}{N^2}  + N\Delta +  \tfrac{N^2\Delta}{n\rho_n},
\end{equation}
and consider two cases a) $N\leqslant n\rho_n$ and b) $N\geqslant n\rho_n$. In the case a), we have that the third term in \eqref{eq:smallest_Delta_rhs} is dominated by the second one. Minimizing then in $N$ the upper bound $\tfrac{L_2\Theta_p n^2\rho_n}{N^2}  + N\Delta$ for \eqref{eq:smallest_Delta_rhs}, we obtain the optimal number $N_{a)}$ of steps and minimal possible value $\e_{a)}$ of this upper bound. Moreover inequality $N_{a)}\leqslant n\rho_n$  turns out to be equivalent to $\Delta \geqslant \frac{L_2\Theta_p }{n\rho_n^2}$. In the case b) the second term in \eqref{eq:smallest_Delta_rhs} is dominated by the third one. Minimizing then in $N$ the upper bound $\tfrac{L_2\Theta_p n^2\rho_n}{N^2}+ \tfrac{N^2\Delta}{n\rho_n}$ for \eqref{eq:smallest_Delta_rhs}, we obtain the optimal number $N_{b)}$ of steps and minimal possible value $\e_{b)}$ of this upper bound. Moreover inequality $N_{b)}\geqslant n\rho_n$  turns out to be equivalent to $\Delta \leqslant \frac{L_2\Theta_p }{n\rho_n^2}$.
Now we can choose $m_{a)} = \tfrac{N_{a)}}{nL_2}\tfrac{\sigma^2}{\e_{a)}}$ and $m_{b)} = \tfrac{N_{b)}}{nL_2}\tfrac{\sigma^2}{\e_{b)}}$ to make the second term in the r.h.s. of \eqref{eq:ARDFDSConv_cor} to be of the same order as the smallest achievable error $\e_{a)}$ or $\e_{b)}$ in the case a) or b) respectively.
Finally, we check that $\Delta \geqslant \frac{L_2\Theta_p }{n\rho_n^2}$ is equivalent to the case a) and inequalities $\e_{a)} \geqslant \e_{b)}$, $N_{a)} \leqslant N_{b)}$,  $m_{a)} \leqslant m_{b)}$. This means that the smallest possible error is $\max\{\e_{a)},\e_{b)}\}$ and it is achieved in the number of iterations $\min\{N_{a)},N_{b)}\}$ with batch size $\min\{m_{a)},m_{b)}\}$. The corresponding values of the parameters are given in Table \ref{tbl:ARDFDS_params_uncontrolled} and we summarize the result as follows.
\begin{corollary}\label{cor:ardfds_conv_uncontr}
    Assume that $\Delta$ is known and satisfies $0 < \Delta \leqslant L_2\Theta_p n^3\rho_n^2/2$, the parameters $N(\Delta),m(\Delta),t(\Delta)$ satisfy relations stated in Table \ref{tbl:ARDFDS_params_uncontrolled} and ARDFDS 
    is applied to solve problem \eqref{eq:PrSt}. 
    Then the output point $y_N$ satisfies $\EE\left[f(y_N)\right] - f(x^*) \leqslant \e(\Delta)$, where $\e(\Delta)$ satisfies the corresponding relation in the same table. Moreover, the overall number of oracle calls is $N(\Delta)m(\Delta)$ given in the same table.
\end{corollary}

{
{\small
\begin{table*}[t]
\centering
    \centering
    \begin{tabular}{|c|c|c|}
        \hline
         & $p=1$ & $p=2$ \\
        \hline
        $t(\Delta)$ & $\sqrt{\Delta/L_2}$ & $\sqrt{\Delta/L_2}$ \\
        \hline
        $N(\Delta)$ & $\min\left\{ \left(\tfrac{L_2\Theta_1 n}{\Delta}\right)^{\nicefrac{1}{3}},\left(\tfrac{L_2\Theta_1 n}{\Delta}\right)^{\nicefrac{1}{4}}\right\}$ 
        & $\min\left\{ \left(\tfrac{L_2\Theta_2 n^2}{\Delta}\right)^{\nicefrac{1}{3}},\left(\tfrac{L_2\Theta_2 n^3}{\Delta}\right)^{\nicefrac{1}{4}}\right\}$ \\
        \hline
        $m(\Delta)$ & $\min\left\{\tfrac{\sigma^2}{nL_2\Delta}, \tfrac{\sigma^2}{(n^{4}L_2^{5}\Theta_1\Delta^{3})^{\nicefrac{1}{4}}} \right\}$ & $\min\left\{\tfrac{\sigma^2}{nL_2\Delta}, \tfrac{\sigma^2}{(n^{3}L_2^{5}\Theta_2\Delta^{3})^{\nicefrac{1}{4}}} \right\}$ \\ 
        \hline	
		$\e(\Delta)$ & $\max\left\{ \left(L_2\Theta_1 n\Delta^2\right)^{\nicefrac{1}{3}},\sqrt{nL_2\Theta_1\Delta}\right\}$ & $\max\left\{ \left(L_2\Theta_2 n^2\Delta^2\right)^{\nicefrac{1}{3}},\sqrt{nL_2\Theta_2\Delta}\right\}$ \\
        \hline
        $Nm$ & $\min\left\{\sigma^2\left(\tfrac{\Theta_1 }{n^2L_2^2\Delta^4}\right)^{\nicefrac{1}{3}}, \tfrac{\sigma^2}{nL_2\Delta}  \right\}$ & $\min\left\{\sigma^2\left(\tfrac{\Theta_2 }{nL_2^2\Delta^4}\right)^{\nicefrac{1}{3}}, \tfrac{\sigma^2 }{L_2\Delta}  \right\}$\\
        \hline
    \end{tabular}
    \caption{Summary of the values for $N,m,t$ and the total number of function value evaluations $Nm$ guaranteeing for the cases $p=1$ and $p=2$ that Algorithm~\ref{Alg:ARDFDS} outputs $y_N$ with minimal possible expected objective residual $\e$ if the oracle noise level $\Delta$ is uncontrolled. Numerical constants and logarithmic factors in $n$ are omitted for simplicity. \label{tbl:ARDFDS_params_uncontrolled}}
\end{table*}}
}
}

\pdd{Using an additional ``light-tail'' assumption that $ \EE_{\xi}[\exp(\|g(x,\xi) - \nabla f(x)\|_2^2/\sigma^2)] \leqslant \exp(1)$ and techniques of \cite{gorbunov2019optimal} our algorithm and analysis can be extended to obtain results in terms of probability of large deviations. For example, in the case of controlled noise level $\Delta$ this means that our algorithm outputs a point $y_N$ which satisfies $\mathbb{P} \{f(y_N) - f(x^*) \leqslant \e\} \geqslant 1-\delta$, where $\delta \in (0,1)$ is the confidence level, for the price of extra $\ln \frac{1}{\delta}$ factor in $N$ and $m$.}

\newstuff{\pdd{In the several next subsections} we provide the full proof of Theorem~\ref{Th:ARDFDSConv} consisting of the four following parts. We start with the technical result providing us with inequalities relating the approximation \eqref{eq:FinDiffStocGrad} with the stochastic gradient \pdd{$g(x,\xi)$} and full gradient \pdd{$\nabla f(x)$}. 
\pdd{The next two parts are in the spirit of Linear Coupling method of \cite{allen2014linear}. Namely, we analyze the progress of the Gradient Descent step (line 5 of ARDFDS) and estimate the progress of the Mirror Descent step (line 6 of ARDFDS).}
\pdd{In the final fourth part,}  
we combine all previous parts and \pdd{finish the proof} of the theorem. We emphasize that in the last part we use a \pdd{careful} analysis of the recurrent inequalities for $\EE[\|x^* - z_k\|_p]$ (see Lemma~\ref{stoh:technical_lemma}, proved in Appendix \ref{S:Tech_result}) in order to bound the terms related to the noise in \pdd{the objective} values.}

\subsubsection{Inequalities for Gradient Approximation}\label{sec:ineq}
The proof of the main theorem relies on the following technical result, which connects finite-difference approximation \eqref{eq:FinDiffStocGrad} of the stochastic gradient with the stochastic gradient itself and also with $\nabla f$. \newstuff{This lemma plays a central role in our analysis \pdd{providing correct} 
dependence of the complexity bounds \pdd{on the dimension.}}
\begin{lemma}
\label{Lm:FinDiffToGrad}
    For all $x, s \in\R^n$, we have 
    \begin{eqnarray}
        \EE_{e}\| \tnmf^t(x)\|_q^2 &\leqslant& \tfrac{12\rho_n}{n}\|g^m(x,\vec{\xi}_m)\|_2^2+\tfrac{\rho_nt^2}{m}\sum\limits_{i=1}^mL(\xi_i)^2+\tfrac{16\rho_n \Delta^2}{t^2},\label{eq:tnfmq2}\\
        \EE_{e}\|\tnmf^t(x)\|_2^2 &\geqslant& \tfrac{1}{2n}\|g^m(x,\vec{\xi}_m)\|_2^2 - \tfrac{t^2}{2m}\sum\limits_{i=1}^mL(\xi_i)^2-\tfrac{8 \Delta^2}{t^2},\label{eq:tnfm22}\\
        \EE_{e} \langle \tnmf^t(x), s\rangle &\geqslant& \tfrac{1}{n}\la g^m(x,\vec{\xi}_m),s\ra - \tfrac{t\|s\|_p}{2m\sqrt{n}}\sum\limits_{i=1}^m L(\xi_i) - \tfrac{2\Delta \|s\|_p}{t\sqrt{n}},\label{eq:tnfms}\\
        \hspace{3em}\EE_{e} \|\langle \nabla f(x),e \ra e - \tnmf^t(x)\|_2^2 &\leqslant& \tfrac{2}{n}\|\nabla f(x) -  g^m(x,\vec{\xi}_m)\|_2^2 + \tfrac{t^2}{m}\sum\limits_{i=1}^mL(\xi_i)^2+\tfrac{16 \Delta^2}{t^2},\label{eq:nf-tnfm22}
    \end{eqnarray}
		where $g^m(x,\vec{\xi}_m) := \tfrac{1}{m}\sum\limits_{i=1}^mg(x,\xi_i)$, $\Delta$ is defined in \eqref{eq:tf_def}, $L(\xi)$ is the Lipschitz constant of $g(\cdot,\xi)$, which is the gradient of $F(\cdot,\xi)$.
\end{lemma}
\begin{proof}
First of all, we rewrite $\tnmf^t(x)$ as follows
\begin{equation}
\notag
\tnmf^t(x) = \left(\left\la g^m(x,\vec{\xi}_m),e \right\ra + \tfrac{1}{m}\sum\limits_{i=1}^m\theta(x,\xi_i,t,e)\right)e,
\end{equation}
where $\theta(x,\xi_i,t,e) = \tfrac{F(x+te,\xi_i) - F(x,\xi_i)}{t} - \la g(x,\xi_i),\, e\ra + \tfrac{\pdd{\eta}(x+te,\xi_i) - \pdd{\eta}(x,\xi_i)}{t}$,  $i=1,\ldots,m$. By \pdd{the $L(\xi)$-}smoothness of $F(\cdot,\xi)$ and \eqref{eq:tf_def}, we have
\begin{equation}
\label{eq:thetaEst}
|\theta(x,\xi_i,t,e)| \leq \tfrac{L(\xi)t}{2} + \tfrac{2\Delta}{t}.
\end{equation}

\textbf{Proof of \eqref{eq:tnfmq2}}.
\begin{equation}
	\begin{array}{rl}
		\EE_{e}\| \tnmf^t(x)\|_q^2
		&= \EE_{e}\Big\|\left(\left\la g^m(x,\vec{\xi}_m),e \right\ra + \tfrac{1}{m}\sum\limits_{i=1}^m\theta(x,\xi_i,t,e)\right)e\Big\|_q^2 \\
		&\overset{\circledOne}{\leqslant} 2\EE_{e}\|\la g^m(x,\vec{\xi}_m),e \ra e\|_q^2 + 2\EE_{e}\left\|\tfrac{1}{m}\sum\limits_{i=1}^m\theta(x,\xi_i,t,e)e\right\|_q^2 \\
		&\overset{\circledTwo}{\leqslant} \tfrac{12\rho_n}{n}\|g^m(x,\vec{\xi}_m)\|_2^2+\tfrac{2\rho_n}{m}\sum\limits_{i=1}^m\left(\tfrac{L(\xi_i)t}{2} + \tfrac{2\Delta}{t}\right)^2\\
		&\leqslant \tfrac{12\rho_n}{n}\|g^m(x,\vec{\xi}_m)\|_2^2+\tfrac{\rho_nt^2}{m}\sum\limits_{i=1}^mL(\xi_i)^2+\tfrac{16\rho_n \Delta^2}{t^2},
	\end{array}
\end{equation}
where $\circledOne$ holds since $\|x+y\|_q^2 \leqslant 2\|x\|_q^2 + 2\|y\|_q^2, \forall x,y\in\R^n$; $\circledTwo$ follows from inequalities \eqref{assumption_e_norm}, \eqref{assumption_e_product}, \eqref{eq:thetaEst} and the fact that, for any $a_1,a_2,\ldots,a_m > 0$, it holds that $\left(\sum\limits_{i=1}^ma_i\right)^2 \leqslant m\sum\limits_{i=1}^ma_i^2$.

\textbf{Proof of \eqref{eq:tnfm22}}.
\begin{equation}
	\begin{array}{rl}
		\EE_{e}\| \tnmf^t(x)\|_2^2
		&= \EE_{e}\Big\|\left(\left\la g^m(x,\vec{\xi}_m),e \right\ra + \tfrac{1}{m}\sum\limits_{i=1}^m\theta(x,\xi_i,t,e)\right)e\Big\|_2^2 \\
		&\overset{\circledOne}{\geqslant} \tfrac{1}{2}\EE_{e}\|\la g^m(x,\vec{\xi}_m),e \ra e\|_2^2 - \tfrac{1}{m}\sum\limits_{i=1}^m\left(\tfrac{L(\xi_i)t}{2} + \tfrac{2\Delta}{t}\right)^2\\ 
		&\overset{\circledTwo}{\geqslant} \tfrac{1}{2n}\|g^m(x,\vec{\xi}_m)\|_2^2 - \tfrac{t^2}{2m}\sum\limits_{i=1}^mL(\xi_i)^2-\tfrac{8 \Delta^2}{t^2},
	\end{array}
\end{equation}
where $\circledOne$ follows from \eqref{eq:thetaEst} and inequality $\|x+y\|_2^2 \geqslant \tfrac{1}{2}\|x\|_2^2 - \|y\|_2^2, \forall x,y\in\R^n$; $\circledTwo$ follows from $e \in RS_2(1)$ and Lemma~B.10 in \cite{bogolubsky2016learning}, stating that, for any $s \in \R^n$, $\EE\la s,e\ra^2 = \tfrac{1}{n} \|s\|_2^2$.

\textbf{Proof of \eqref{eq:tnfms}}.
\begin{equation}
	\begin{array}{rl}
		\EE_{e} \langle \tnmf^t(x), s\rangle  &= \EE_{e} \langle \la g^m(x,\vec{\xi}_m),e\ra e, s\rangle + \EE_{e}  \tfrac{1}{m}\sum\limits_{i=1}^m\theta(x,\xi_i,t,e) \la e, s \ra \\
		&\overset{\circledOne}{\geqslant} \tfrac{1}{n}\la g^m(x,\vec{\xi}_m),s\ra - \tfrac{1}{m}\sum\limits_{i=1}^m\left( \tfrac{L(\xi_i)t}{2} + \tfrac{2\Delta}{t} \right) \EE_{e}  |\la e,s \ra | \\
		&\overset{\circledTwo}{\geqslant} \tfrac{1}{n}\la g^m(x,\vec{\xi}_m),s\ra - \tfrac{t\|s\|_p}{2m\sqrt{n}}\sum\limits_{i=1}^m L(\xi_i) - \tfrac{2\Delta \|s\|_p}{t\sqrt{n}}
	\end{array}
\end{equation}
where $\circledOne$ follows from $\EE_{e}[n\la g, e \ra e] = g,\, \forall g\in\R^n$ and \eqref{eq:thetaEst}; 
$\circledTwo$ follows from Lemma B.10 in \cite{bogolubsky2016learning}, since $\EE|\la s,e\ra| \leq \sqrt{\EE\la s,e\ra^2}$, and the fact that $\|x\|_2 \leqslant \|x\|_p$ for $p\leqslant2$.

\textbf{Proof of \eqref{eq:nf-tnfm22}}.
\begin{equation}
	\begin{array}{r}
		\EE_{e} \|\langle \nabla f(x),e \ra e - \tnmf^t(x)\|_2^2\\
		= \EE_{e} \left\|\la \nabla f(x),e \ra e - \la g^m(x,\vec{\xi}_m),e \ra e - \tfrac{1}{m}\sum\limits_{i=1}^m\theta(x,\xi_i,t,e)  e\right\|_2^2 \\
		\overset{\circledOne}{\leqslant}  2\EE_{e} \left\|\la \nabla f(x) -  g^m(x,\vec{\xi}_m),e \ra e \right\|_2^2 + 
		2\EE_{e} \left\| \tfrac{1}{m}\sum\limits_{i=1}^m\theta(x,\xi_i,t,e)  e\right\|_2^2  \\
		\overset{\circledTwo}{\leqslant} \tfrac{2}{n}\|\nabla f(x) -  g^m(x,\vec{\xi}_m)\|_2^2 + \tfrac{t^2}{m}\sum\limits_{i=1}^mL(\xi_i)^2+\tfrac{16 \Delta^2}{t^2},
	\end{array}
\end{equation}
where $\circledOne$ holds since $\|x+y\|_2^2 \leqslant 2\|x\|_2^2 + 2\|y\|_2^2, \forall x,y\in\R^n$; 
$\circledTwo$ follows from $e \in S_2(1)$ and Lemma~B.10 in \cite{bogolubsky2016learning}, and \eqref{eq:thetaEst}.
\end{proof}

\subsubsection{\newstuff{Progress of the Gradient Descent Step}}\label{sec:progress}
The following lemma estimates the progress in step \newstuff{5} of ARDFDS, which is a gradient step.
\begin{lemma}
\label{Lm:GradStep}
    Assume that $y = x - \tfrac{1}{2L_2}\tnmf^t(x)$.
		Then, 
		\begin{equation}		
		\label{eq:GradStepProgr}
		\| g^m(x,\vec{\xi}_m)\|_2^2 \leq 8nL_2 (f(x) - \EE_ef(y)) + 8 \|\nabla f(x) -  g^m(x,\vec{\xi}_m)\|_2^2+ \tfrac{5nt^2}{m}\sum\limits_{i=1}^mL(\xi_i)^2+\tfrac{80n \Delta^2}{t^2},
		\end{equation}	
		where $g^m(x,\vec{\xi}_m)$ is defined in Lemma \ref{Lm:FinDiffToGrad}, $\Delta$ is defined in \eqref{eq:tf_def}, $L(\xi)$ is the Lipschitz constant of $g(\cdot,\xi)$, which is the gradient of $F(\cdot,\xi)$.
\end{lemma}
\begin{proof}
Since $\tnmf^t(x)$ is collinear to $e$, we have that, for some $\gamma \in \R$, $y-x = \gamma e$. Then, since $\|e\|_2=1$,
\begin{equation*}
    \begin{array}{c}
	    \langle \nabla f(x), y-x \rangle = \langle \nabla f(x), e \rangle \gamma = \langle \nabla f(x), e \rangle \langle e, y-x \rangle = \langle \langle \nabla f(x),e \rangle e, \, y-x\rangle.
	\end{array}
\end{equation*}
From this and $L_2$-smoothness of $f$, we obtain
\begin{equation*}
	\begin{array}{rl}
	    f(y) &\leqslant f(x) + \langle \langle \nabla f(x),e \rangle e, \, y-x\rangle + \tfrac{L_2}{2}||y-x||_2^2 \\
	    &= f(x) + \langle  \tnmf^t(x), \, y-x\rangle + L_2||y-x||_2^2 + \langle \langle \nabla f(x),e \rangle e - \tnmf^t(x), \, y-x\rangle\\
	    &\quad - \tfrac{L_2}{2}||y-x||_2^2 \\
        &\overset{\circledOne}{\leqslant} f(x) + \langle \tnmf^t(x), \, y-x\rangle + L_2||y-x||_2^2  + \tfrac{1}{2L_2}\| \langle \nabla f(x),e \rangle e - \tnmf^t(x)\|_2^2, 
	\end{array}
\end{equation*}
where $\circledOne$ follows \pd{from} the Fenchel inequality $\la s,z \ra - \tfrac{\zeta}{2}\|z\|_2^2 \leq \tfrac{1}{2\zeta}\|s\|_2^2$.
Using $y = x-\tfrac{1}{2L_2}\tnmf^t(x)$, we get
\begin{equation*}
    \begin{array}{c}
        \tfrac{1}{4L_2} \|\tnmf^t(x)\|_2^2 \leqslant f(x)-f(y) + \tfrac{1}{2L_2}\|\langle \nabla f(x),e \rangle e - \tnmf^t(x)\|_2^2.
    \end{array}
\end{equation*}
Taking the expectation in $e$ we obtain
\begin{equation*}
    \begin{array}{c}
    \tfrac{1}{4L_2} \left(\tfrac{1}{2n}\|g^m(x,\vec{\xi}_m)\|_2^2 - \tfrac{t^2}{2m}\sum\limits_{i=1}^mL(\xi_i)^2-\tfrac{8 \Delta^2}{t^2}\right)   \overset{\eqref{eq:tnfm22}}{\leqslant} \tfrac{1}{4L_2} \EE_e\|\tnmf^t(x)\|_2^2 \\
		\leqslant f(x)-\EE_ef(y) + \tfrac{1}{2L_2}\EE_e\|\langle \nabla f(x),e \rangle e - \tnmf^t(x)\|_2^2 \\
		 \overset{\eqref{eq:nf-tnfm22}}{\leqslant} f(x)-\EE_ef(y) + \tfrac{1}{2L_2} \left( \tfrac{2}{n}\|\nabla f(x) -  g^m(x,\vec{\xi}_m)\|_2^2 + \tfrac{t^2}{m}\sum\limits_{i=1}^mL(\xi_i)^2+\tfrac{16 \Delta^2}{t^2}\right).
    \end{array}
\end{equation*}
Rearranging the terms, we obtain the statement of the lemma.
\end{proof}

\subsubsection{Progress of the Mirror Descent Step}\label{sec:mirr_des_stp}
The following lemma estimates the progress in step \newstuff{6} of ARDFDS, which is a Mirror Descent step.
\begin{lemma}
\label{Lm:MDStep}
    For $z_+ = \argmin\limits_{u\in\R^n} \left\{ {\alpha n \left\langle \tnmf^t(x), \, u-z \right\rangle +V[z] \left( u \right)}\right\}$ we have
		\begin{equation}		
		\label{eq:MDStepProgr}
		\begin{array}{rl}
		\alpha \la g^m(x,\vec{\xi}_m),z-u\ra &\leqslant 6\alpha^2n \rho_n \|g^m(x,\vec{\xi}_m)\|_2^2 + V[z](u) - \EE_{e}[V[z_{+}](u) \\
		&+ \tfrac{\alpha^2n^2\rho_n}{2} \left( \tfrac{t^2}{m}\sum\limits_{i=1}^mL(\xi_i)^2+\tfrac{16\Delta^2}{t^2} \right)\\
		&+ \alpha \sqrt{n} \|z-u\|_p \left( \tfrac{t}{2m}\sum\limits_{i=1}^m L(\xi_i) + \tfrac{2\Delta }{t} \right),
		\end{array}
		\end{equation}	
		where $g^m(x,\vec{\xi}_m)$ is defined in Lemma \ref{Lm:FinDiffToGrad}, $\Delta$ is defined in \eqref{eq:tf_def}, $L(\xi)$ is the Lipschitz constant of $g(\cdot,\xi)$, which is the gradient of $F(\cdot,\xi)$.
\end{lemma}
\begin{proof}
For all $u\in \R^n$, we have
\begin{equation}	
\label{lemma110:basic_estimations}
	\begin{array}{rl}
		\alpha n\langle \tnmf^t(x), \, z-u\rangle
		&=  \alpha n \langle \tnmf^t(x), \,  z-z_{+}\rangle + \alpha n \langle \tnmf^t(x), \, z_{+}-u\rangle\\
		& \hspace{-5em}\overset{\circledOne}{\leqslant}  \alpha n \langle \tnmf^t(x), \, z-z_{+}\rangle + \langle -\nabla V[z](z_{+}), \,  z_{+}-u\rangle\\
		&\hspace{-5em}\overset{\circledTwo}{=}  \alpha n \langle \tnmf^t(x), \, z-z_{+}\rangle
		+ V[z](u) - V[z_{+}](u) - V[z](z_{+})\\
		&\hspace{-5em}\overset{\circledThree}{\leqslant} \left(\alpha n \langle \tnmf^t(x), \,  z-z_{+}\rangle - \tfrac{1}{2}\|z-z_{+}\|_p^2\right) + V[z](u) - V[z_+](u)\\
		&\hspace{-5em}\overset{\circledFour}{\leqslant} \tfrac{\alpha^2n^2}{2}\| \tnmf^t(x)\|_q^2 + V[z](u) - V[z_{+}](u),
	\end{array}
\end{equation}
where $\circledOne$ follows from the definition of $z_+$, whence $\langle \nabla V[z](z_{+}) +  \alpha n \tnmf^t(x), \, u - z_{+}\rangle \geqslant 0$ for all $u\in \R^n$; $\circledTwo$ follows from the ``magic identity'' Fact 5.3.2 in \cite{ben-tal2020lectures} for the Bregman divergence; $\circledThree$ follows from \eqref{eq:VStrConv}; and $\circledFour$ follows from the Fenchel inequality $\zeta\la s,z \ra - \tfrac{1}{2}\|z\|_p^2 \leq \tfrac{\zeta^2}{2}\|s\|_q^2$. Taking expectation in $e$, applying \eqref{eq:tnfms} with $s = z-u$ and \eqref{eq:tnfmq2}, we get
\begin{equation}
\label{lemma110:averaged_basic_estimation}
	\begin{array}{c}
	\alpha n \left( \tfrac{1}{n}\la g^m(x,\vec{\xi}_m),z-u\ra - \tfrac{t\|z-u\|_p}{2m\sqrt{n}}\sum\limits_{i=1}^m L(\xi_i) - \tfrac{2\Delta \|z-u\|_p}{t\sqrt{n}} \right)\\ \leqslant
		\alpha n\EE_{e}\langle \tnmf^t(x), \, z-u\rangle
		\leqslant \tfrac{\alpha^2n^2}{2}\EE_{e}\| \tnmf^t(x)\|_q^2 + V[z](u) - \EE_{e}[V[z_{+}](u)] \\
		\leqslant \tfrac{\alpha^2n^2}{2} \left( \tfrac{12\rho_n}{n}\|g^m(x,\vec{\xi}_m)\|_2^2+\tfrac{\rho_nt^2}{m}\sum\limits_{i=1}^mL(\xi_i)^2+\tfrac{16\rho_n \Delta^2}{t^2} \right) + V[z](u) - \EE_{e}[V[z_{+}](u)].
	\end{array}
\end{equation}
Rearranging the terms, we obtain the statement of the lemma.
\end{proof}

\subsubsection{Proof of Theorem \ref{Th:ARDFDSConv}}\label{sec:proof1}
First, we prove the following lemma, which estimates the per-iteration progress of the whole algorithm.
\begin{lemma}
\label{Lm:OneStep}
		Let $\{x_k,y_k,z_k,\alpha_k,\tau_k\}$, $k \geqslant 0$ be generated by ARDFDS. Then, for all $u \in \R^n$,
		\begin{equation}
		\label{eq:OneStep}
        \begin{array}{rl}
            48 n^2 \rho_n L_2 \alpha_{k+1}^2\EE_{e,\xi}[f(y_{k+1})\mid \cE_k,\Xi_k] - (48 n^2 \rho_n L_2 \alpha_{k+1}^2-\alpha_{k+1})f(y_k)&\\ - V[z_k](u) + \EE_{e,\xi}[V[z_{k+1}](u)\mid \cE_k,\Xi_k] - \RR_{k+1} &\leqslant \alpha_{k+1}f(u),
        \end{array}
    \end{equation}
		\begin{equation}		
		\label{eq:ErrTermDef}
		\begin{array}{c}
		    \RR_{k+1} = \tfrac{48 \alpha_{k+1}^2 n \rho_n \sigma^2}{m}
		+\tfrac{61\alpha_{k+1}^2n^2\rho_n}{2} \left( L_2^2t^2+\tfrac{16\Delta^2}{t^2} \right)
		+ \alpha_{k+1} \sqrt{n} \|z_k-u\|_p \left( \tfrac{L_2t}{2} + \tfrac{2\Delta }{t} \right),
		\end{array}
		\end{equation}
		where $\Delta$ is defined in \eqref{eq:tf_def},  
		$\cE_k$ and $\Xi_k$ denote the history of realizations of $e_1,\ldots,e_k$ and $\xi_{1}^{1},\ldots,\xi_{m}^{1},\ldots,\xi_{1}^{k},\ldots,\xi_{m}^{k}$ respectively, up to the step $k$.
\end{lemma}
\begin{proof}
Combining \eqref{eq:GradStepProgr} and \eqref{eq:MDStepProgr}, we obtain
\begin{equation}		
		\label{eq:OneStepPr1}
		\begin{array}{rl}
		\alpha \la g^m(x_{k+1},\vec{\xi}_{m}^{\hspace{0.2em}k+1}),z-u\ra &\leqslant 48\alpha^2n^2 \rho_n L_2 (f(x_{k+1}) - \EE_ef(y_{k+1})) \\
		&\hspace{-13em}+ V[z_k](u) - \EE_{e}[V[z_{k+1}](u)]+ 48 \alpha^2n \rho_n\|\nabla f(x_{k+1}) -  g^m(x_{k+1},\vec{\xi}_{m}^{\hspace{0.2em}k+1})\|_2^2 \\
		&\hspace{-13em}+ \tfrac{61\alpha^2n^2\rho_n}{2} \left( \tfrac{t^2}{m}\sum\limits_{i=1}^mL(\xi_i^{k+1})^2+\tfrac{16\Delta^2}{t^2} \right)+ \alpha \sqrt{n} \|z_k-u\|_p \left( \tfrac{t}{2m}\sum\limits_{i=1}^m L(\xi_i^{k+1}) + \tfrac{2\Delta }{t} \right),
		\end{array}
\end{equation}
where $g^m(x,\vec{\xi}_m)$ is defined in Lemma \ref{Lm:FinDiffToGrad} and the expectation in $e$ is conditional on $\cE_k$. 
By the definition of $g^m(x,\vec{\xi}_m)$ and \eqref{stoch_assumption_on_variance}, \pdd{for all $x \in \R^n$,} $\EE_{\xi}g^m(x,\vec{\xi}_m) = \nabla f(x)$ and 
$\EE_{\xi}\|\nabla f(x) -  g^m(x,\vec{\xi}_{m})\|_2^2 \leq \tfrac{\sigma^2}{m}$.
Using these two facts and taking the expectation in $\vec{\xi}_{m}^{\hspace{0.2em}k+1}$ conditional on $\Xi_k$, we obtain 
\begin{equation}
		\label{eq:OneStepPr2}
    \begin{array}{rl}
        \alpha_{k+1}\la \nabla f(x_{k+1}) , z_k-u \ra  
	    &\leqslant 48 \alpha_{k+1}^2 n^2 \rho_n L_2 \left(f(x_{k+1}) - \EE_{e,\xi}[f(y_{k+1})\mid \cE_k,\Xi_k]\right) \\
		&+ V[z_k](u) - \EE_{e,\xi}[V[z_{k+1}](u)\mid \cE_k,\Xi_k] + \RR_{k+1}.
    \end{array}
\end{equation}
Further,
\begin{equation*}
    \begin{array}{rl}
        \alpha_{k+1}&(f(x_{k+1}) - f(u)) \leqslant \alpha_{k+1}\langle \nabla f(x_{k+1}) ,\, x_{k+1} - u \rangle\\
        &= \alpha_{k+1}\langle\nabla f(x_{k+1}), \, x_{k+1}-z_k\rangle + \alpha_{k+1}\langle\nabla f(x_{k+1}), \, z_{k}-u\rangle\\
	    &\overset{\circledOne}{=} \tfrac{(1-\tau_k)\alpha_{k+1}}{\tau_k}\langle\nabla f(x_{k+1}), \, y_{k}-x_{k+1}\rangle +\alpha_{k+1}\langle\nabla f(x_{k+1}), \, z_{k}-u\rangle\\
		&\overset{\circledTwo}{\leqslant} \tfrac{(1-\tau_k)\alpha_{k+1}}{\tau_k} (f(y_k)-f(x_{k+1})) + \alpha_{k+1}\langle\nabla f(x_{k+1}), \, z_{k}-u\rangle \\
		&\overset{\eqref{eq:OneStepPr2}}{\leqslant} \tfrac{(1-\tau_k)\alpha_{k+1}}{\tau_k} (f(y_k)-f(x_{k+1}))\\ 
		&+  48 \alpha_{k+1}^2 n^2 \rho_n L_2 \left(f(x_{k+1}) - \EE_{e,\xi}[f(y_{k+1})\mid \cE_k,\Xi_k]\right)  \\
		&+ V[z_k](u) - \EE_{e,\xi}[V[z_{k+1}](u)\mid \cE_k,\Xi_k] + \RR_{k+1} \\
		&\overset{\circledThree}{=} (48 \alpha_{k+1}^2 n^2 \rho_n L_2-\alpha_{k+1})f(y_k)\\
		& - 48 \alpha_{k+1}^2 n^2 \rho_n L_2\EE_{e,\xi}[f(y_{k+1})\mid \cE_k,\Xi_k] +\alpha_{k+1}f(x_{k+1})\\
		&+ V[z_k](u) - \EE_{e,\xi}[V[z_{k+1}](u) \mid \cE_k,\Xi_k] + \RR_{k+1}.
    \end{array}
\end{equation*}
\pdd{Here} $\circledOne$ is since $x_{k+1} := \tau_kz_k+(1-\tau_k)y_k\; \Leftrightarrow\; \tau_k(x_{k+1}-z_k) = (1-\tau_k)(y_k - x_{k+1})$, $\circledTwo$ follows from the convexity of $f$ and inequality $1-\tau_k\geqslant0$, and $\circledThree$ is since $\tau_k = \tfrac{1}{48\alpha_{k+1}n^2\rho_n L_2}$.
Rearranging the terms, we obtain the statement of the lemma.
\end{proof}

\begin{proof}[Proof of Theorem \ref{Th:ARDFDSConv}]
Note that \\
\noindent $48n^2\rho_nL_2\alpha_{k+1}^2 - \alpha_{k+1} + \tfrac{1}{192n^2\rho_nL_2} = 48n^2\rho_nL_2\alpha_{k}^2$ \pdd{since}
\begin{equation*}
    \begin{array}{rl}
        48n^2\rho_nL_2\alpha_{k+1}^2 & - \alpha_{k+1} + \tfrac{1}{192n^2\rho_nL_2} = \tfrac{(k+2)^2}{192n^2\rho_nL_2} - \tfrac{k+2}{96n^2\rho_nL_2} + \tfrac{1}{192n^2\rho_nL_2}\\
        &= \tfrac{k^2+4k+4-2k-4+1}{192n^2\rho_nL_2}= \tfrac{(k+1)^2}{192n^2\rho_nL_2}= 48n^2\rho_nL_2\alpha_{k}^2.
    \end{array}
\end{equation*}
Taking, \pdd{for any $1 \leqslant l\leqslant N$,}  \pdd{the} full expectation $\EE[\cdot] = \EE_{e_1,\ldots,e_N,\xi_{1}^{1},\ldots,\xi_{m}^{1},\ldots,\xi_{1}^{N},\ldots,\xi_{m}^{N}}[\cdot]$ \pdd{in} both sides of \eqref{eq:OneStep} for $k=0,\ldots,l-1$ and telescoping the obtained inequalities\footnote{Note that $\alpha_1 = \tfrac{2}{96n^2\rho_nL_2} = \tfrac{1}{48n^2\rho_nL_2}$ and therefore $48n^2\rho_nL_2\alpha_1^2 - \alpha_1 = 0$.}, we have, 
\begin{equation}
\label{theo_after_telescoping_mini_gr_free}
    \begin{array}{rl}
        48n^2\rho_nL_2\alpha_{l}^2\EE[f(y_l)] + \sum\limits_{k=1}^{l-1}\tfrac{1}{192n^2\rho_nL_2}\EE[f(y_k)] - V[z_0](u)&\\ + \EE[V[z_l](u)] 
				- \zeta_1\sum\limits_{k=0}^{l-1}\alpha_{k+1}\EE[\|u-z_k\|_p] -\zeta_2\sum\limits_{k=0}^{l-1}\alpha_{k+1}^2   &\leqslant \sum\limits_{k=0}^{l-1}\alpha_{k+1}f(u),
    \end{array}
\end{equation}
where we denoted 
\begin{equation}
\label{eq:zeta1zeta2Def}
\zeta_1:= \sqrt{n}\left( \tfrac{L_2t}{2} + \tfrac{2\Delta }{t} \right), \quad \zeta_2 := 48 n \rho_n \tfrac{\sigma^2}{m}
		+\tfrac{61n^2\rho_n}{2} \left( L_2^2t^2+\tfrac{16\Delta^2}{t^2} \right).
\end{equation}
\pdd{Since $u$ in \eqref{theo_after_telescoping_mini_gr_free} is arbitrary, we set $u=x^*$}, where $x^*$ is a solution to \eqref{eq:PrSt}, \pd{ use the inequality $\Theta_p \geqslant V[z_0](x^*)$,
 and define} $R_k:=\EE[\|x^*-z_k\|_p]$. Also, from \eqref{eq:VStrConv}, we have that $\zeta_1\alpha_1R_0 \leq \tfrac{\sqrt{2\Theta_p}\zeta_1}{48n^2\rho_nL_2}$. To simplify the notation, we define $B_l := \zeta_2\sum\limits_{k=0}^{l-1}\alpha_{k+1}^2 + \Theta_p + \tfrac{\sqrt{2\Theta_p}\zeta_1}{48n^2\rho_nL_2}$ . Since $\sum\limits_{k=0}^{l-1}\alpha_{k+1} = \tfrac{l(l+3)}{192n^2\rho_nL_2}$ and, for all $i=1,\ldots,N$, $f(y_i) \leqslant f(x^*)$, we get from \eqref{theo_after_telescoping_mini_gr_free} that
\begin{equation}
\label{theo_basic_inequality_mini_gr_free}
    \begin{array}{rl}
        \tfrac{(l+1)^2}{192n^2\rho_nL_2}\EE[f(y_l)] &\leqslant f(x^*)\left(\tfrac{(l+3)l}{192n^2\rho_nL_2} - \tfrac{l-1}{192n^2\rho_nL_2}\right)
        + B_l\\
        &- \EE[V[z_l](x^*)] + \zeta_1\sum\limits_{k=1}^{l-1}\alpha_{k+1}R_k,\\
        0 \leqslant \tfrac{(l+1)^2}{192n^2\rho_nL_2}\left(\EE[f(y_l)] - f(x^*)\right) &\leqslant B_l - \EE[V[z_l](x^*)] + \zeta_1\sum\limits_{k=1}^{l-1}\alpha_{k+1}R_k,
    \end{array}
\end{equation}
which gives  $\EE[V[z_l](x^*)] \leqslant B_l + \zeta_1\sum\limits_{k=1}^{l-1}\alpha_{k+1}R_k$ and
\begin{equation}\label{theo_distance_estimation_mini_gr_free}
        \tfrac{1}{2}\left(\EE[\|z_l-x^*\|_p]\right)^2 \leqslant \tfrac{1}{2}\EE[\|z_l-x^*\|_p^2] \leqslant \EE[V[z_l](x^*)] 
        \leqslant B_l + \zeta_1\sum\limits_{k=1}^{l-1}\alpha_{k+1}R_k,
\end{equation}
whence, $R_l \leqslant \sqrt{2}\cdot\sqrt{B_l + \zeta_1\sum\limits_{k=1}^{l-1}\alpha_{k+1}R_k}$.
\pd{This recurrent sequence of $R_l$'s is analyzed separately in Appendix \ref{S:Tech_result}.  
Applying Lemma~\ref{stoh:technical_lemma}} 
 with $a_0 = \zeta_2\alpha_1^2 + \Theta_p + \tfrac{\sqrt{2\Theta_p}\zeta_1}{48n^2\rho_nL_2}, a_{k} = \zeta_2\alpha_{k+1}^2, b = \zeta_1$ for $k=1,\ldots,N-1$, we obtain
\begin{equation}\label{theo_inequality_for_induction_mini_gr_free}
    \begin{array}{c}
        B_l + \zeta_1\sum\limits_{k=1}^{l-1}\alpha_{k+1}R_k \leqslant \left(\sqrt{B_l} + \sqrt{2}\zeta_1\cdot\tfrac{l^2}{96n^2\rho_nL_2}\right)^2, \; l=1,\ldots,N
    \end{array}
\end{equation} 
Since $V[z](x^*) \geqslant 0$, by inequality \eqref{theo_basic_inequality_mini_gr_free}, for $l=N$ and the definition of $B_l$, 
we have
\begin{equation}\label{theo_pre_final_mini_gr_free}
    \begin{array}{rl}
        \tfrac{(N+1)^2}{192n^2\rho_nL_2}\left(\EE[f(y_N)] - f(x^*)\right) &\leqslant \left(\sqrt{B_N} + \sqrt{2}\zeta_1\cdot\tfrac{N^2}{96n^2\rho_nL_2}\right)^2\\
        &\hspace{-12em}\overset{\circledOne}{\leqslant} 2B_N + 4\zeta_1^2\cdot\tfrac{N^4}{(96n^2\rho_nL_2)^2}= 2\zeta_2\sum\limits_{k=0}^{l-1}\alpha_{k+1}^2 + 2\Theta_p + \tfrac{\sqrt{2\Theta_p}\zeta_1}{24n^2\rho_nL_2} + \tfrac{4\zeta_1^2N^4}{(96n^2\rho_nL_2)^2}\\
        &\overset{\circledTwo}{\leqslant} 2\Theta_p + \tfrac{\sqrt{2\Theta_p}\zeta_1}{24n^2\rho_nL_2} + \tfrac{2\zeta_2(N+1)^3}{(96n^2\rho_nL_2)^2}  + \tfrac{4\zeta_1^2N^4}{(96n^2\rho_nL_2)^2}
    \end{array}
\end{equation}
where $\circledOne$ is due to the fact that, $\forall a,b\in\R,\quad (a+b)^2 \leqslant 2a^2+2b^2$ and $\circledTwo$ is because $\sum\limits_{k=0}^{N-1}\alpha_{k+1}^2 = \tfrac{1}{(96n^2\rho_nL_2)^2}\sum\limits_{k=2}^{N+1}k^2 \leqslant \tfrac{1}{(96n^2\rho_nL_2)^2} \cdot \tfrac{(N+1)(N+2)(2N+3)}{6} \leqslant \tfrac{1}{(96n^2\rho_nL_2)^2} \cdot \tfrac{(N+1)2(N+1)3(N+1)}{6} = \tfrac{(N+1)^3}{(96n^2\rho_nL_2)^2}$. Dividing \eqref{theo_pre_final_mini_gr_free} by $\tfrac{(N+1)^2}{192n^2\rho_nL_2}$ and substituting $\zeta_1,\zeta_2$ from \eqref{eq:zeta1zeta2Def}, we obtain
\begin{equation*}
    \begin{array}{rl}
        \EE[f(y_N)] - f(x^*) &\leqslant \tfrac{384\Theta_p n^2\rho_nL_2}{(N+1)^2} + \tfrac{12\sqrt{2\Theta_p}}{(N+1)^2}\zeta_1 + \tfrac{384(N+1)\zeta_2}{(96n^2\rho_nL_2)^2} + \tfrac{N^4\zeta_1^2}{12n^2\rho_nL_2(N+1)^2}\\
        &\leqslant \tfrac{384\Theta_p n^2\rho_nL_2}{N^2} + \tfrac{12\sqrt{2n\Theta_p}}{N^2}\left( \tfrac{L_2t}{2} + \tfrac{2\Delta }{t} \right)+ \tfrac{384N}{nL_2}\tfrac{\sigma^2}{m} \\
				&+ \tfrac{6N}{L_2} \left( L_2^2t^2+\tfrac{16\Delta^2}{t^2} \right) + \tfrac{N^2}{24n\rho_nL_2} \left( L_2^2t^2 + \tfrac{16\Delta^2 }{t^2} \right).
    \end{array}
\end{equation*}

\end{proof}

\subsection{\newstuff{Randomized Derivative-Free Directional Search}}
\pdd{Our non-accelerated method} is listed as Algorithm \ref{Alg:RDFDS}.
\pd{Following \cite{nesterov2017random,ghadimi2016mini-batch,ghadimi2013stochastic}  we assume that $L_2$ is known.} 
\pd{The possible choice of the smoothing parameter $t$ and mini-batch size $m$ are discussed below.}
\pd{Note that at every iteration the algorithm requires to solve an auxiliary minimization problem. As it is shown in \cite{ben-tal2020lectures}, for both cases $p=1$ and $p=2$ this minimization can be made explicitly in $O(n)$ arithmetic operations.}
\begin{algorithm}
	\caption{Randomized Derivative-Free Directional Search (RDFDS)}
	\label{Alg:RDFDS}
	\begin{algorithmic}[1]
		\REQUIRE $x_0$~-- starting point; $N$~-- number of iterations; \pd{$L_2$~-- smoothness constant}; $m \geqslant 1$~-- mini-batch size; $t > 0$~-- smoothing parameter\pd{;} $\alpha = \tfrac{1}{48n\rho_n L_2}$~-- stepsize\pd{; $V(\cdot,\cdot)$~-- Bregman divergence}.
		\FOR{$k=0,\, \dots, \, N-1$.}
		\STATE Generate $e_{k+1} \in RS_2\left( 1 \right)$ independently from previous iterations and $\xi_i^{\pdd{k+1}}$, $i=1,...,m$ -- independent  realizations of $\xi$\pd{, which are also independent from previous iterations}. 
		\STATE \pd{Calculate $\tnmf^t(x_{k})$ \pd{using} \eqref{eq:FinDiffStocGrad} \pd{with $e=e_{k+1}$}.}
		\STATE $x_{k+1} \leftarrow \argmin\limits_{x\in\R^n} \Big\{ \alpha n \left\langle \tnmf^t(x_{k}), \, x-x_k \right\rangle+V[x_{k}] \left( x \right)\Big\}.$
		\ENDFOR
		\RETURN$\bar{x}_N \leftarrow \tfrac{1}{N}\sum\limits_{k=0}^{N-1}x_k$.
	\end{algorithmic}
\end{algorithm}

\begin{theorem}\label{theorem_convergence_mini_gr_free_non_acc}
    Let RDFDS 
    be applied to solve problem \eqref{eq:PrSt}, \pd{ $x^*$ be an arbitrary solution to \eqref{eq:PrSt}, and $\Theta_p$ be such that $ V[z_0](x^*)\leq \Theta_p$.} 
    Then
    \begin{equation}\label{theo_main_result_mini_gr_free}
        \begin{array}{rcl}
        \EE[f(\bar{x}_N)] - f(x_*) &\leqslant& \tfrac{384n\rho_nL_2\Theta_p}{N} + \tfrac{2\sigma^2}{L_2m} + \left(\tfrac{n}{6L_2} + \tfrac{N}{3L_2\rho_n} \right)\left(\tfrac{L_2^2t^2}{2} + \tfrac{8\Delta^2}{t^2}\right) \\
        &&\quad+ \tfrac{8\sqrt{2n\Theta_p}}{N}\left(\tfrac{L_2t}{2} + \tfrac{2\Delta}{t}\right),\qquad \forall n\ge 8.
    \end{array}
    \end{equation}
\end{theorem}
\begin{proof}[Proof of Theorem \ref{theorem_convergence_mini_gr_free_non_acc}]
    The proof of this result is \pdd{rather} similar to the proof of Theorem~\ref{Th:ARDFDSConv}. First of all,
    \begin{equation}\label{theo_ds_der_free:basic_estimations}
	\begin{array}{r}
		\alpha n\langle \tnmf^t(x_{k}), \, x_{k}-x_*\rangle
		\\=  \alpha n \langle \tnmf^t(x_{k}), \,  x_{k}-x_{k+1}\rangle + \alpha n \langle \tnmf^t(x_{k}), \, x_{k+1}-x_*\rangle\\
		\overset{\circledOne}{\leqslant}  \alpha n \langle \tnmf^t(x_{k}), \, x_{k}-x_{k+1}\rangle + \langle -\nabla V[x_k](x_{k+1}), \,  x_{k+1}-x_*\rangle\\
		\overset{\circledTwo}{=}  \alpha n \langle \tnmf^t(x_{k}), \, x_{k}-x_{k+1}\rangle
		+ V[x_k](x_*) - V[x_{k+1}](x_*) - V[x_k](x_{k+1})\\
		\overset{\circledThree}{\leqslant} \left(\alpha n \langle \tnmf^t(x_{k}), \,  x_{k}-x_{k+1}\rangle - \tfrac{1}{2}\|x_k-x_{k+1}\|_p^2\right) + V[x_k](x_*) - V[x_{k+1}](x_*)\\
		\leqslant \tfrac{\alpha^2n^2}{2}\| \tnmf^t(x_{k}\|_q^2 + V[x_k](x_*) - V[x_{k+1}](x_*),
	\end{array}
\end{equation}
where $\circledOne$ follows from $\langle \nabla V[x_k](x_{k+1}) +  \alpha n \tnmf^t(x_{k}), \, x - x_{k+1}\rangle \geqslant 0$ for all $x\in \R^n$, $\circledTwo$ follows from ``magic identity'' Fact 5.3.2 in \cite{ben-tal2020lectures} for Bregman divergence, and $\circledThree$ is \pdd{since} $V[x](y) \geqslant \tfrac{1}{2}\|x-y\|_p^2$. Taking conditional expectation $\EE_{e}[\;\cdot\mid \cE_k]$ \pdd{in} both sides of \eqref{theo_ds_der_free:basic_estimations} we get
\begin{equation}\label{theo_ds_der_free:averaged_basic_estimation}
	\begin{array}{rl}
		\alpha n\EE_{e}[\langle \tnmf^t(x_{k}), \, x_{k}-x_*\rangle\mid \cE_k] 
		&\leqslant \tfrac{\alpha^2n^2}{2}\EE_{e}[\| \tnmf^t(x_{k})\|_q^2\mid \cE_k]\\
		&+ V[x_k](x_*) - \EE_{e}[V[x_{k+1}](x_*)\mid \cE_k]
	\end{array}
\end{equation}
From \eqref{theo_ds_der_free:averaged_basic_estimation}, \eqref{eq:tnfmq2} and \eqref{eq:tnfms} for $s = x_k-x_*$, we obtain
\begin{equation*}
    \begin{array}{rl}
        \langle g^m(x_{k},\vec{\xi}_{m}^{\hspace{0.2em}k+1}), \, x_{k}-x_*\rangle &\leqslant 24\alpha^2n\rho_nL_2(f(x_k) - f(x_*))\\
        &\hspace{-10em}+ 12\alpha^2n\rho_n\|\nabla f(x_k) - g^m(x_{k},\vec{\xi}_{m}^{\hspace{0.2em}k+1})\|_2^2+ \alpha^2n^2\rho_n\cdot\tfrac{t^2}{2m}\sum\limits_{i=1}^mL_2(\xi_{i}^{k+1})^2 + \tfrac{8\alpha^2n^2\rho_n\Delta^2}{t^2}\\ 
        &\hspace{-10em}+ \alpha\sqrt{n}\|x_k-x_*\|_p\cdot\tfrac{t}{2m}\sum\limits_{i=1}^m L_2(\xi_{i}^{k+1}) + \tfrac{2\alpha\Delta\sqrt{n}\|x_k-x_*\|_p}{t}\\
        &\hspace{-10em}+ V[x_k](x_*) - \EE_{e}[V[x_{k+1}](x_*)\mid \cE_k].
    \end{array}
\end{equation*}
Taking conditional expectation $\EE_{\xi}[\;\cdot\mid \Xi_{k}]$
\pdd{in} the both sides of \pdd{the} previous inequality and using \pdd{the} convexity of $f$ and \eqref{stoch_assumption_on_variance}, we have
\begin{equation}\label{theo_ds_der_free:pre_final}
    \begin{array}{rl}
        \underbrace{(\alpha - 24\alpha^2n\rho_nL_2)}_{\alpha/4}\left(f(x_k) - f(x_*)\right) &\leqslant  12\alpha^2n\rho_n\tfrac{\sigma^2}{m}
        + \alpha^2n^2\rho_n\left(\tfrac{L_2^2t^2}{2}
        +\tfrac{8\Delta^2}{t^2}\right)\\
        & \hspace{-13em}+ \alpha\sqrt{n}\|x_k-x_*\|_p\left(\tfrac{L_2t}{2} + \tfrac{2\Delta}{t}\right)+ V[x_k](x_*) - \EE_{e,\xi}[V[x_{k+1}](x_*)\mid \cE_k,\Xi_k],
    \end{array} 
\end{equation}
\pdd{since} $\alpha = \tfrac{1}{48n\rho_nL_2}$. Denote
\begin{equation}\label{eq:zeta_1_zeta_2_def_non_acc}
    \zeta_1 = \tfrac{L_2t}{2} + \tfrac{2\Delta}{t},\quad \zeta_2 = \tfrac{L_2^2t^2}{2}+\tfrac{8\Delta^2}{t^2}.
\end{equation}
Note that
\begin{equation}\label{eq:zeta_2_zeta_1_non_acc}
    \zeta_1^2 = \left(\tfrac{L_2t}{2} + \tfrac{2\Delta}{t}\right)^2 \leqslant 2\cdot\tfrac{L_2^2t^2}{4} + 2\cdot\tfrac{4\Delta^2}{t^2} = \zeta_2.
\end{equation}
Taking \pdd{for any $1 \leqslant l\leqslant N$,}  \pdd{the} full expectation $\EE[\cdot] = \EE_{e_1,\ldots,e_N,\xi_{1}^{1},\ldots,\xi_{m}^{1},\ldots,\xi_{1}^{N},\ldots,\xi_{m}^{N}}[\cdot]$ \pdd{in} both sides of
inequalities \eqref{theo_ds_der_free:pre_final} for $k=0,\ldots,l-1$ and summing them, we get
\begin{equation}\label{theo_ds_der_free:after_telescoping}
    \begin{array}{rl}
        0 \leqslant \tfrac{N\alpha}{4}\left(\EE[f(\bar{x}_l)] - f(x_*)\right) &\leqslant l\cdot12\alpha^2n\rho_n\tfrac{\sigma^2}{m} + l\alpha^2n^2\rho_n\zeta_2\\
        &+\alpha\sqrt{n}\zeta_1\sum\limits_{k=0}^{l-1}\EE[\|x_k-x_*\|_p]
        \pd{+ V[x_0](x_*)}
        - \EE[V[x_{l}](x^*)],
    \end{array}
\end{equation}
where $\bar{x}_l =
\tfrac{1}{l}\sum\limits_{k=0}^{l-1}x_k$.
From the previous inequality, \pd{since $V[z_0](x^*) \leqslant \Theta_p$,} we get
\begin{equation}\label{theo_ds_der_free:distance_estimation}
    \begin{array}{c}
        \tfrac{1}{2}\left(\EE[\|x_l-x_*\|_p]\right)^2 \leqslant \tfrac{1}{2}\EE[\|x_l-x_*\|_p^2] \leqslant \EE[V[x_l](x_*)]\\
        \leqslant \Theta_p + l\cdot12\alpha^2n\rho_n\tfrac{\sigma^2}{m} + l\alpha^2n^2\rho_n\zeta_2
        +\alpha\sqrt{n}\delta\zeta_1\sum\limits_{k=0}^{l-1}\EE[\|x_k-x_*\|_p],
    \end{array}
\end{equation}
whence, $\forall l\leqslant N$, we obtain
\begin{equation}
    \EE[\|x_{\pd{l}}-x_*\|_p] \leqslant \sqrt{2}\sqrt{\Theta_p + l\cdot12\alpha^2n\rho_n\tfrac{\sigma^2}{m} + l\alpha^2n^2\rho_n\zeta_2
    +\alpha\sqrt{n}\zeta_1\sum\limits_{k=0}^{l-1}\EE[\|x_k-x_*\|_p]}.
\end{equation}
Denote $R_k = \EE[\|x^*-x_k\|_p]$ for $k=0,\ldots,N$. 
\pd{The recurrent sequence of $R_k$'s is analyzed separately in Appendix \ref{S:Tech_result}.  
Applying Lemma~\ref{stoh:technical_lemma_non_acc}} 
\pdd{with} $a_0 = \Theta_p + \alpha\sqrt{n}\zeta_1\EE[\|x_0-x_*\|_p] \leqslant \Theta_p + \alpha\sqrt{2n\Theta_p}\zeta_1, a_{k} = 12\alpha^2n\rho_n\tfrac{\sigma^2}{m}+\alpha^2n^2\rho_n\zeta_2, b = \sqrt{n}\zeta_1$ for $k=1,\ldots,N-1$ we have for $l=N$
\begin{equation*}
    \begin{array}{c}
        \tfrac{N\alpha}{4}\left(\EE[f(\bar{x}_N)] - f(x_*)\right)\\ 
        \leqslant \left(\sqrt{\Theta_p+N\cdot12\alpha^2n\rho_n\tfrac{\sigma^2}{m} + N\alpha^2n^2\rho_n\zeta_2 + \alpha\sqrt{2n\Theta_p}\zeta_1} + \sqrt{2n}\zeta_1\alpha N\right)^2\\
        \overset{\circledOne}{\leqslant} 2\Theta_p + 24N\alpha^2n\rho_n\tfrac{\sigma^2}{m} + 2N\alpha^2n^2\rho_n\zeta_2 + 2\alpha\sqrt{2n\Theta_p}\zeta_1 + 4n\zeta_1^2\alpha^2N^2,
    \end{array}
\end{equation*}
whence
\begin{equation*}
    \begin{array}{rl}
        \EE[f(\bar{x}_N)] - f(x_*) &\overset{\eqref{eq:zeta_2_zeta_1_non_acc}}{\leqslant} \tfrac{384n\rho_nL_2\Theta_p}{N} + \tfrac{2\sigma^2}{L_2m} + \tfrac{n\zeta_2}{6L_2} + \tfrac{8\sqrt{2n\Theta_p}\zeta_1}{N} + \tfrac{\zeta_2N}{3L_2\rho_n}\\
        &\overset{\eqref{eq:zeta_1_zeta_2_def_non_acc}}{\leqslant} \tfrac{384n\rho_nL_2\Theta_p}{N} + \tfrac{2\sigma^2}{L_2m} + \left(\tfrac{n}{6L_2} + \tfrac{N}{3L_2\rho_n} \right)\left(\tfrac{L_2^2t^2}{2} + \tfrac{8\Delta^2}{t^2}\right)\\
        &+ \tfrac{8\sqrt{2n\Theta_p}}{N}\left(\tfrac{L_2t}{2} + \tfrac{2\Delta}{t}\right),
    \end{array}
\end{equation*}
where we used also that $\alpha = \tfrac{1}{48n\rho_nL_2}$.
\end{proof}
\pd{
Similarly to the discussion above concerning the ARDFDS and its convergence theorem, we can formulate corollaries for the RDFDS in the case of controlled and uncontrolled noise level $\Delta$.
In the simple case $\Delta=0$, all the terms in the r.h.s. of \eqref{theo_main_result_mini_gr_free} can be made smaller than $\e$ for any $\e \geqslant 0 $ by an appropriate choice of $N,m,t$. 
Thus, we consider a more interesting case and assume that the noise level satisfies $0<\Delta \leqslant L_2\Theta_p n \rho_n^2/2$, the second inequality being non-restrictive. In order to minimize the term with $\tfrac{L_2^2t^2}{2} + \tfrac{8\Delta^2}{t^2}$ in the r.h.s of \eqref{theo_main_result_mini_gr_free}, we set $t=2\sqrt{\tfrac{\Delta}{L_2}}$. Substituting this into the r.h.s. of \eqref{theo_main_result_mini_gr_free} and using that, by our assumption on $\Delta$,  $\sqrt{2nL_2\Theta_p\Delta}\leqslant n\rho_nL_2\Theta_p$, we obtain an upper bound for $\EE[f(\bar{x}_N)] - f(x_*)$. Following the same steps as in the proof of Corollaries \ref{cor:ardfds_conv} and \ref{cor:ardfds_conv_uncontr}, we obtain the following results for RDFDS.
}

\begin{corollary}\label{cor:rdfds_conv}
    \pd{
    Assume that the value of $\Delta$ can be controlled and satisfies $0<\Delta \leqslant L_2\Theta_p n \rho_n^2/2$.
    Assume that for a given accuracy $\e \geqslant 0$ the values of the parameters $N(\e),m(\e),t(\e),\Delta(\e)$ satisfy the relations stated in Table \ref{tbl:RDFDS_params} and RDFDS is applied to solve problem \eqref{eq:PrSt}.
    Then the output point $\bar{x}_N$ satisfies $\EE\left[f(\bar{x}_N)\right] - f(x^*) \le \e$. Moreover, the overall number of oracle calls is $N(\e)m(\e)$ given in the same table.}
\end{corollary}
{
{\small
\begin{table*}[t]
\centering
    \centering
    \begin{tabular}{|c|c|c|}
        \hline
         & $p=1$ & $p=2$ \\
        \hline
        $N(\e)$ & $\tfrac{\ln nL_2\Theta_1}{\varepsilon}$ & $\tfrac{ nL_2\Theta_2}{\varepsilon}$\\
        \hline
        $m(\e)$ & $\max\left\{1,\tfrac{\sigma^2}{L_2\varepsilon}\right\}$ & $\max\left\{1,\tfrac{\sigma^2}{L_2\varepsilon}\right\}$ \\ 
		\hline	
        $\Delta(\e)$ & $\min\left\{\tfrac{\varepsilon}{n},\, \tfrac{\varepsilon^2}{nL_2\Theta_1}\right\}$ & $\min\left\{\tfrac{\varepsilon}{n},\, \tfrac{\varepsilon^2}{nL_2\Theta_2}\right\}$\\
        \hline	
        $t(\e)$ & $\min\left\{\sqrt{\tfrac{\varepsilon}{nL_2}},\, \tfrac{\varepsilon}{\sqrt{nL_2^2\Theta_1}}\right\}$ & $\min\left\{\sqrt{\tfrac{\varepsilon}{nL_2}},\, \tfrac{\varepsilon}{\sqrt{nL_2^2\Theta_1}}\right\}$\\
        \hline
        \pd{$N(\e)m(\e)$}
        & $\max\left\{\tfrac{L_2\Theta_1\ln n}{\varepsilon}, \tfrac{\sigma^2\Theta_1 \ln n}{\varepsilon^2}\right\}$ & $\max\left\{\tfrac{nL_2\Theta_2}{\varepsilon}, \tfrac{n\sigma^2\Theta_2}{\varepsilon^2}\right\}$\\
        \hline
    \end{tabular}
    \caption{
    \pd{Summary of the values for $N,m,\Delta,t$ and the total number of function value evaluations $Nm$ guaranteeing for the cases $p=1$ and $p=2$ that Algorithm~\ref{Alg:RDFDS} outputs $\bar{x}_N$ satisfying $\EE\left[f(\bar{x}_N)\right] - f(x_*) \le \varepsilon$. Numerical constants are omitted for simplicity.}\label{tbl:RDFDS_params}
    }
\end{table*}}
}
\pd{Note that in the case of uncontrolled noise level $\Delta$, the values of this parameter stated in Table \ref{tbl:RDFDS_params} can be seen as the maximum value of the noise level which can be tolerated by the method still allowing it to achieve $\EE\left[f(\bar{x}_N)\right] - f(x^*) \leqslant \e$. For a more general case of uncontrolled noise level $\Delta$, we obtain the following Corollary.}
\pd{
\begin{corollary}\label{cor:rdfds_conv_uncontr}
    Assume that $\Delta$ is known and satisfies $0<\Delta \leqslant L_2\Theta_p n \rho_n^2/2$, the parameters $N(\Delta),m(\Delta),t(\Delta)$ satisfy relations stated in Table \ref{tbl:RDFDS_params_uncontrolled} and RDFDS 
    is applied to solve problem \eqref{eq:PrSt}.
    Then the output point \newstuff{$\bar{x}_N$} satisfies $\EE\left[f(\bar{x}_N)\right] - f(x^*) \leqslant \e(\Delta)$, where $\e(\Delta)$ satisfies the corresponding relation in the same table. Moreover, the overall number of oracle calls is $N(\Delta)m(\Delta)$ given in the same table.
\end{corollary}
\pdd{Similarly to ARDFDS, RDFDS and its analysis can be extended to obtain convergence in terms of probability of large deviations under additional ``light-tail'' assumption.}
{
{\small
\begin{table*}[t]
\centering
    \centering
    \begin{tabular}{|c|c|c|}
        \hline
         & $p=1$ & $p=2$ \\
        \hline
        $t(\Delta)$ & $\sqrt{\Delta/L_2}$ & $\sqrt{\Delta/L_2}$ \\
        \hline
        $N(\Delta)$ & $\min\left\{ \tfrac{L_2\Theta_1 }{n\Delta},\sqrt{\tfrac{L_2\Theta_1 }{n\Delta}}\right\}$
        & $\min\left\{ \tfrac{L_2\Theta_2 }{\Delta},\sqrt{\tfrac{nL_2\Theta_2 }{\Delta}}\right\}$ \\
        \hline
        $m(\Delta)$ & $\min\left\{\tfrac{\sigma^2}{nL_2\Delta}, \tfrac{\sigma^2}{\sqrt{nL_2^{3}\Theta_1\Delta}} \right\}$ & $\min\left\{\tfrac{\sigma^2}{nL_2\Delta}, \tfrac{\sigma^2}{\sqrt{nL_2^{3}\Theta_2\Delta}} \right\}$ \\
        \hline    
        $\e(\Delta)$ & $\max\left\{ n\Delta,\sqrt{nL_2\Theta_1\Delta}\right\}$
        & $\max\left\{ n\Delta,\sqrt{nL_2\Theta_2\Delta}\right\}$ \\
        \hline
        $N(\Delta)m(\Delta)$ & $\min\left\{\tfrac{\sigma^2\Theta_1 }{n^2\Delta^2}, \tfrac{\sigma^2}{nL_2\Delta} \right\}$ & $\min\left\{\tfrac{\sigma^2\Theta_2 }{n\Delta^2}, \tfrac{\sigma^2}{L_2\Delta} \right\}$\\
        \hline
    \end{tabular}
    \caption{Summary of the values for $N,m,t$ and the total number of function value evaluations $Nm$ guaranteeing for the cases $p=1$ and $p=2$ that Algorithm~\ref{Alg:RDFDS} outputs \newstuff{$\bar{x}_N$} with minimal possible expected objective residual $\e$ if the oracle noise $\Delta$ is uncontrolled. Numerical constants and logarithmic factors in $n$ are omitted for simplicity. \label{tbl:RDFDS_params_uncontrolled}}
\end{table*}}
}
}

\subsection{Role of the algorithms parameters}
\label{sec:role_of_parameters}
{\textbf{Role of $\Delta$ and $t$.}} We \pd{would like} to mention that \textit{there is no need to know the noise level $\Delta$ to run our algorithms}. As it can be seen from \eqref{eq:ARDFDSConv}, the ARDFDS method is robust in the sense of \cite{nemirovski2009robust} to the choice of the smoothing parameter $t$. Namely, if we under/overestimate \pdd{$\Delta$}
by a constant factor, the corresponding terms in the convergence rate will increase only by a constant factor. \pdd{Similar remark holds for the assumption that $L_2$ is known.}

Our Theorems~\ref{Th:ARDFDSConv} and \ref{theorem_convergence_mini_gr_free_non_acc} are applicable in two situations, the noise being a) controlled and b) uncontrolled.
\begin{enumerate}
    \item[a)] Our assumptions on the noise level in Tables~\ref{tbl:ARDFDS_params} and \ref{tbl:RDFDS_params} can be met in practice. For example, in \cite{bogolubsky2016learning}, the objective function is defined by some auxiliary problem and its value can be calculated with accuracy $\Delta$ at the cost proportional to $\ln\tfrac{1}{\Delta}$, which would result in only a $\ln\tfrac{1}{\e}$ factor in the total complexity of our methods in this paper \pdd{combined with the method in \cite{bogolubsky2016learning} for approximating the function value.} 
    \item[b)] 
    \pdd{The minimum guaranteed accuracy $\e(\Delta)$ in Tables~\ref{tbl:ARDFDS_params_uncontrolled} and \ref{tbl:RDFDS_params_uncontrolled} can not be arbitrarily small, which} is reasonable: one can not solve the problem with better accuracy than the accuracy of the available information. 
    \pd{Interestingly, the minimal possible accuracy for the accelerated method could be larger than for the non-accelerated method, which means that accelerated methods are less robust to noise (cf. full gradient methods \cite{devolder2014first,kamzolov2020universal}). To illustrate this, let us, for simplicity neglect the numerical constants and consider a case with $\sigma=0$, $L_2=1$, $\Theta_p=1$, $t=2\sqrt{\Delta}$, and large $N \gg n\rho_n$. Then the main terms in the r.h.s. of \eqref{eq:ARDFDSConv} are $\tfrac{n^2\rho_n}{N^2} +   \tfrac{N^2\Delta}{n\rho_n}$. Minimizing in $N$, we have the minimal accuracy of the order $\sqrt{n\Delta}$.
    Similarly, the main terms in the r.h.s. of 
     \eqref{theo_main_result_mini_gr_free} are 
     $     \tfrac{n\rho_n}{N} +  \tfrac{N\Delta}{\rho_n} $. Minimizing in $N$, we have the minimal accuracy of the order $\Delta^{\nicefrac{1}{2}}n^{\nicefrac{1}{2}}\rho_n^{\newstuff{\nicefrac{1}{2}}} < \sqrt{n\Delta}$, which is smaller than for the accelerated method.}
\end{enumerate}

\textbf{Role of $\sigma^2$.} Although, all the related works, which we are aware of, assume $\sigma^2$ to be known, adaptivity to the variance $\sigma^2$ is a very important direction of future work. \pdd{Note that similarly to the robustness to $\Delta$, our method is robust to $\sigma^2$.}

\section{Experiments}\label{sec:experiments_main}
\newstuff{We performed several numerical experiments to illustrate our theoretical results. In particular, we compared our methods with the Euclidean and \pd{$1$}-norm proximal setups and the RSGF method from \cite{ghadimi2013stochastic} applied to two problems: minimizing Nesterov's function and logistic regression problem. For all the results reported below we tuned parameters $\alpha_k$ and $\alpha$ for ARDFDS and RDFDS respectively and the stepsize parameter for RSGF. We use \_E and \_NE in the plots to refer to the methods with \pd{$2$-norm and $1$}-norm proximal setups respectively and RSGF to refer to the method from \cite{ghadimi2013stochastic}. The code was written in Python using standard libraries, see the details at \url{https://github.com/eduardgorbunov/ardfds}.}

\subsection{Experiments with Nesterov's function}
\newstuff{We tested our methods on the problem of minimizing Nesterov's function \cite{nesterov2004introduction} defined as:
$$f(x) = \tfrac{L_2}{4}\left(\tfrac{1}{2}\left[(x^1)^2 + \sum\limits_{i=1}^{n-1}(x^i - x^{i+1})^2 + (x^n)^2\right] - x^1\right),$$
where $x^i$ is $i$-th component of vector $x\in\R^n$. 
$f$ is convex, $L_2$-smooth w.r.t. the Euclidean norm and attains its minimal value $f^* = \tfrac{L_2}{8}\left(-1+\tfrac{1}{n+1}\right)$ at the point $x^* = (x^{*,1},\ldots,x^{*,n})^\top$ such that $x^{*,i} = 1 - \tfrac{i}{n+1}$. Moreover, the lower complexity bound for first-order methods in smooth convex optimization is attained \cite{nesterov2004introduction} on this function.
}

\newstuff{We add stochastic noise to this function and consider $F(x,\xi) = f(x) + \xi\langle a, x\rangle$, where $\xi$ is Gaussian random variable with mean $\mu = 0$ and variance $\sigma^2$, $a\in\R^n$ is some vector in the unit Euclidean sphere, i.e. $\|a\|_2^2 = 1$. This implies that $f(x) = \EE_\xi\left[F(x,\xi)\right]$ and $F(x,\xi)$ is $L_2$-smooth in $x$ w.r.t. the Euclidean norm since $g(x,\xi) - g(y,\xi) = \nabla f(x) - \nabla f(y)$. Moreover, $\EE_\xi g(x,\xi) = \nabla f(x)$ and $\EE_\xi\left[\|g(x,\xi) - \nabla f(x)\|_2^2\right] = \|a\|_2^2\EE_\xi\left[\xi^2\right] = \sigma^2$ for all $x\in\R^n$. Also we introduce an additive noise $\eta(x) = \Delta\sin\left(\|x-x^*\|_2^{-2}\right)$. It is clear that $|\eta(x)| \le \Delta$ for all $x\in\R^n$. Overall, we are in the setting described in Introduction with $\widetilde{f}(x,\xi) = F(x,\xi) + \eta(x) = f(x) + \xi\langle a,x\rangle + \Delta\sin\left(\|x-x_*\|_2^{-2}\right)$.}

\newstuff{
We compare our methods with the Euclidean and \pd{$1$}-norm proximal setups as well as the RSGF method from \cite{ghadimi2013stochastic} applied to this problem for different sparsity levels of $x_0 - x^*$ and different values of $n$, $\sigma$ and $\Delta$. For all tests we use $L_{\pd{2}} = 10$, adjust starting point $x_0$ such that $f(x_0)-f(x^*) \sim 10^{2}$ and choose $t = \max\{10^{-8},2\sqrt{\nicefrac{\Delta}{L_{\pd{2}}}}\}$. The second term under the maximum in the definition of $t$ corresponds to the optimal choice of $t$ for given $\Delta$ and $L_{\pd{2}}$, i.e., it minimizes the right-hand sides of \eqref{eq:ARDFDSConv} and \eqref{theo_main_result_mini_gr_free}, and the first term under the maximum is needed to prevent unstable computations when $t$ is too small.}

\subsubsection{Experiments with different sparsity levels}
\newstuff{In this set of experiments we considered different choices of the starting point $x_0$ with different sparsity levels of $x_0-x^*$, i.e., for $n = 100, 500, 1000$ we picked such starting points $x_0$ that vector $x_0-x^*$ has $1, \nicefrac{n}{10}, \nicefrac{n}{2}$ non-zero components. In particular, we shift first $1, \nicefrac{n}{10}, \nicefrac{n}{2}$ components of $x^*$ by some constant to obtain $x_0$. In order to isolate the effect of the sparsity from effects coming from the stochastic nature of $\widetilde{f}(x,\xi)$ and noise $\eta(x)$ we choose $\sigma = \Delta = 0$. Our results are reported in Figure~\ref{fig:nesterov_sparsity_experiment}.
\begin{figure}
\centering 
\includegraphics[scale=0.14]{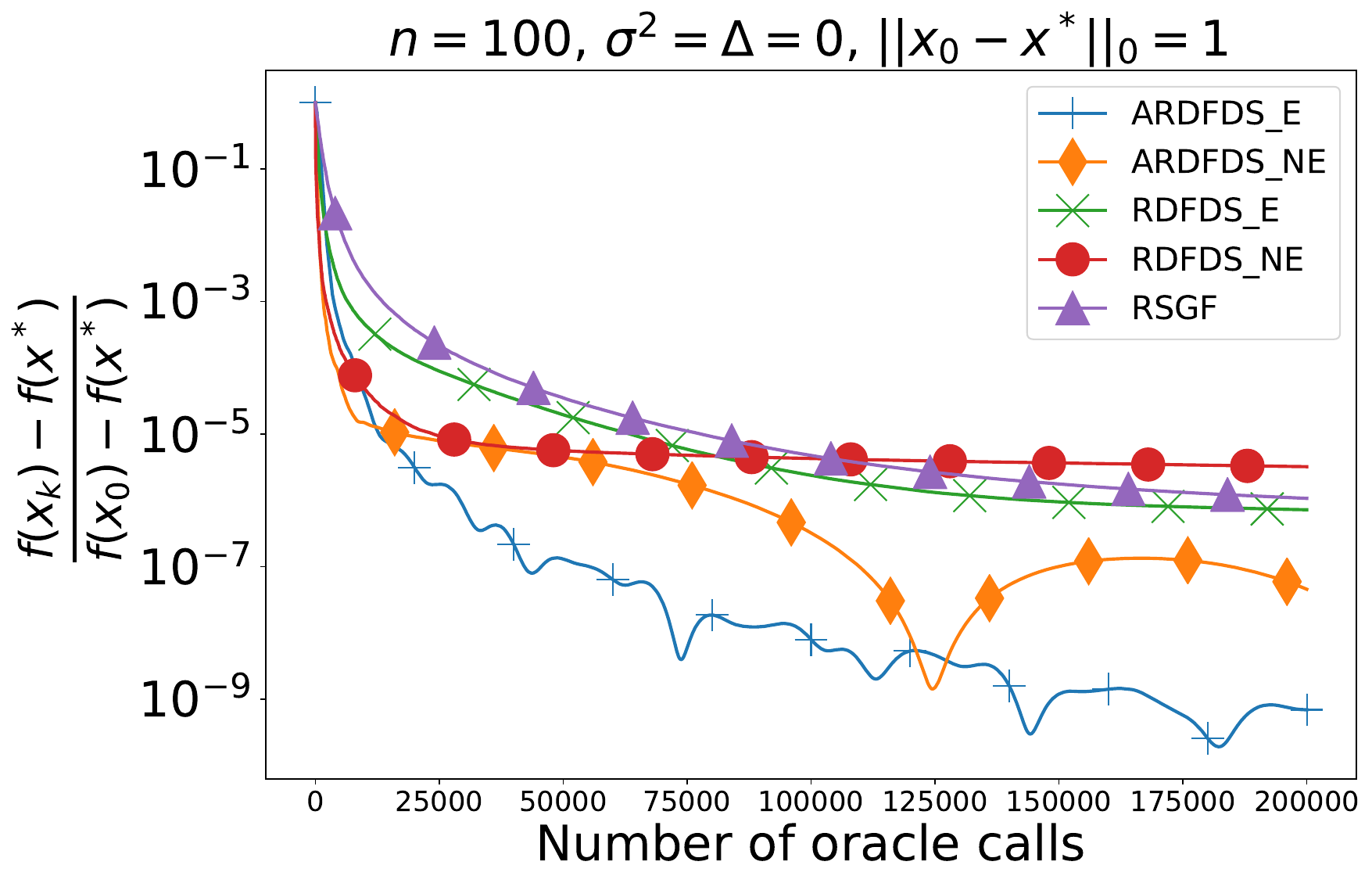}
\includegraphics[scale=0.14]{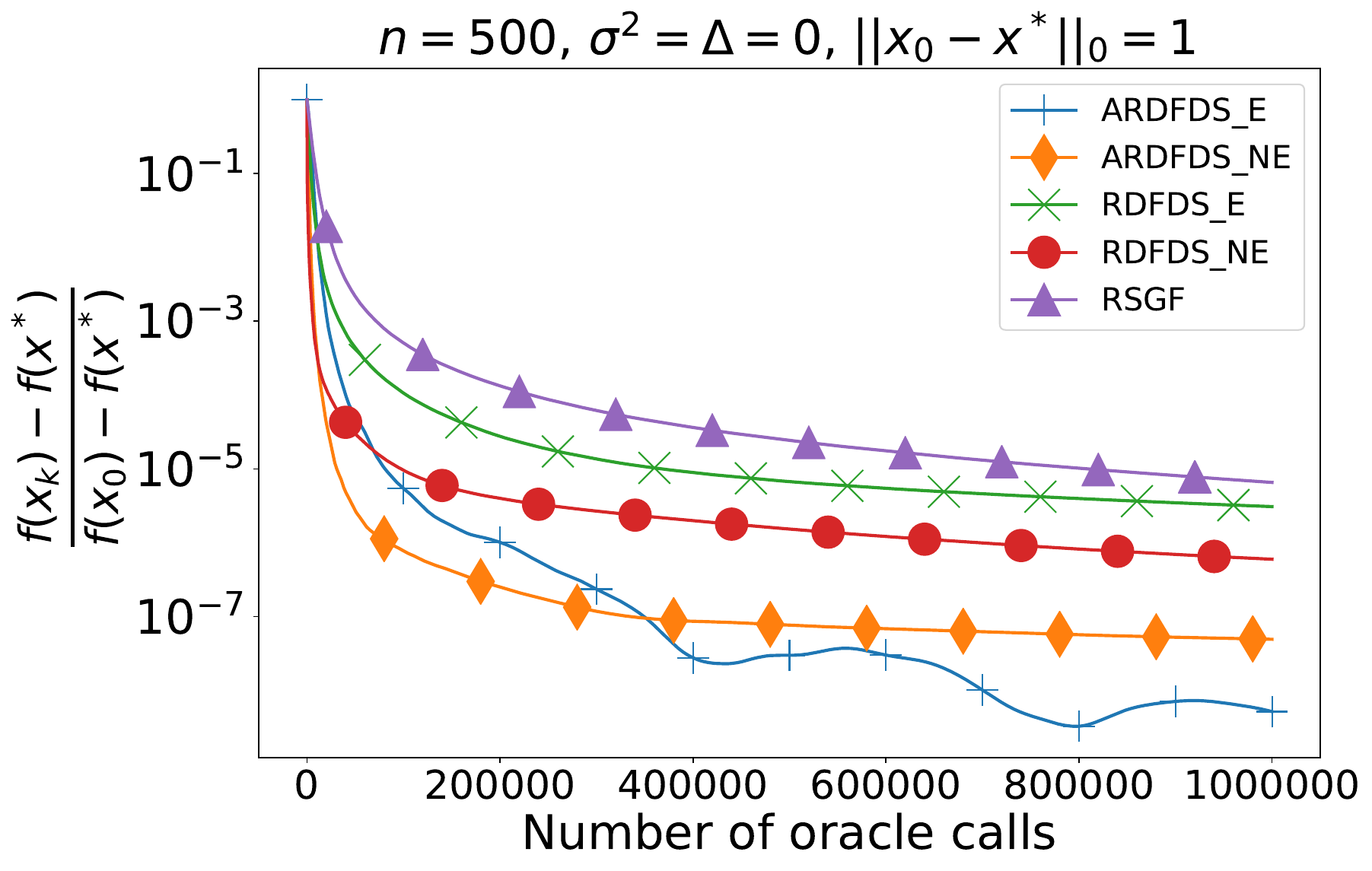}
\includegraphics[scale=0.14]{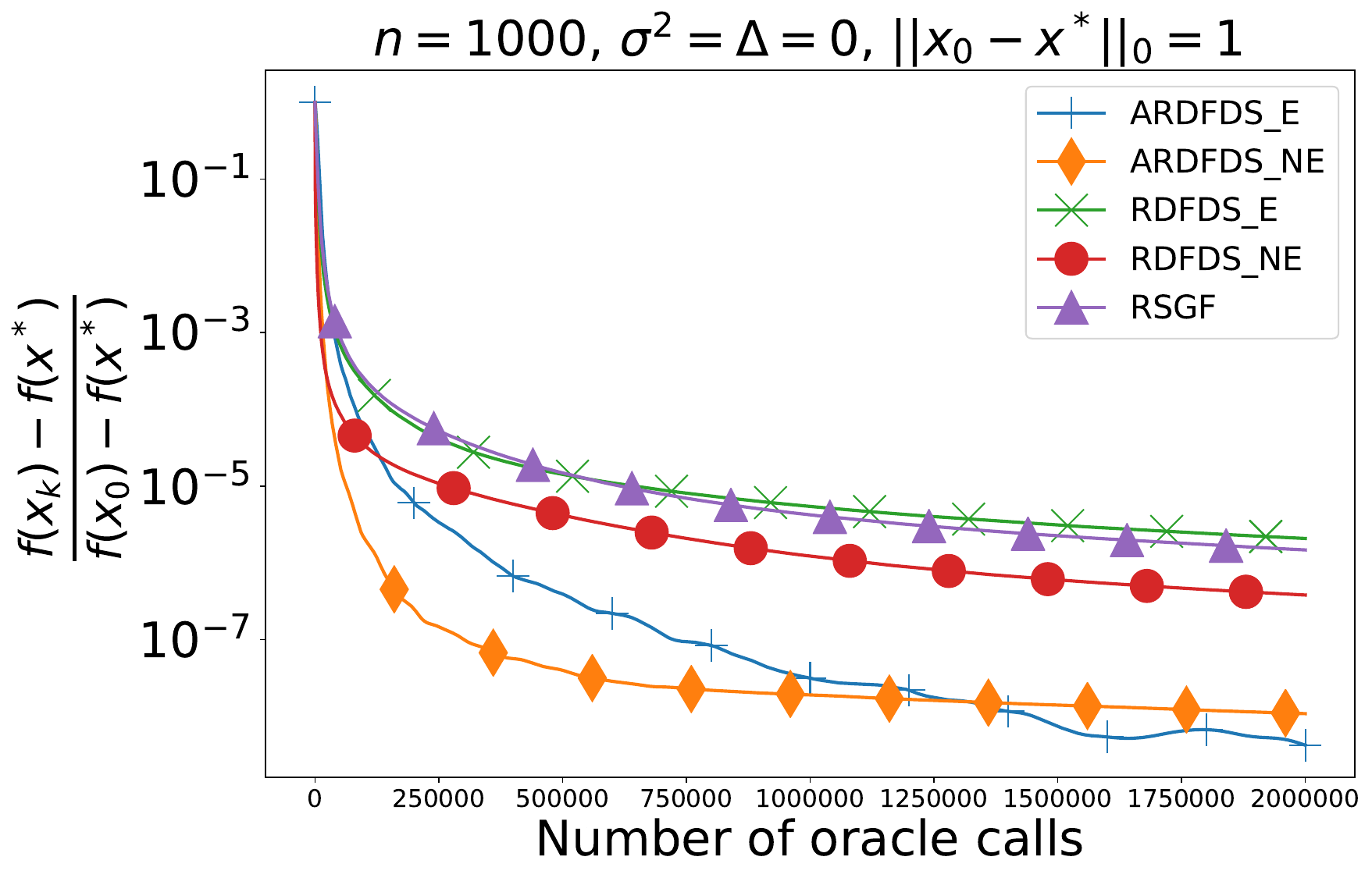}
\includegraphics[scale=0.14]{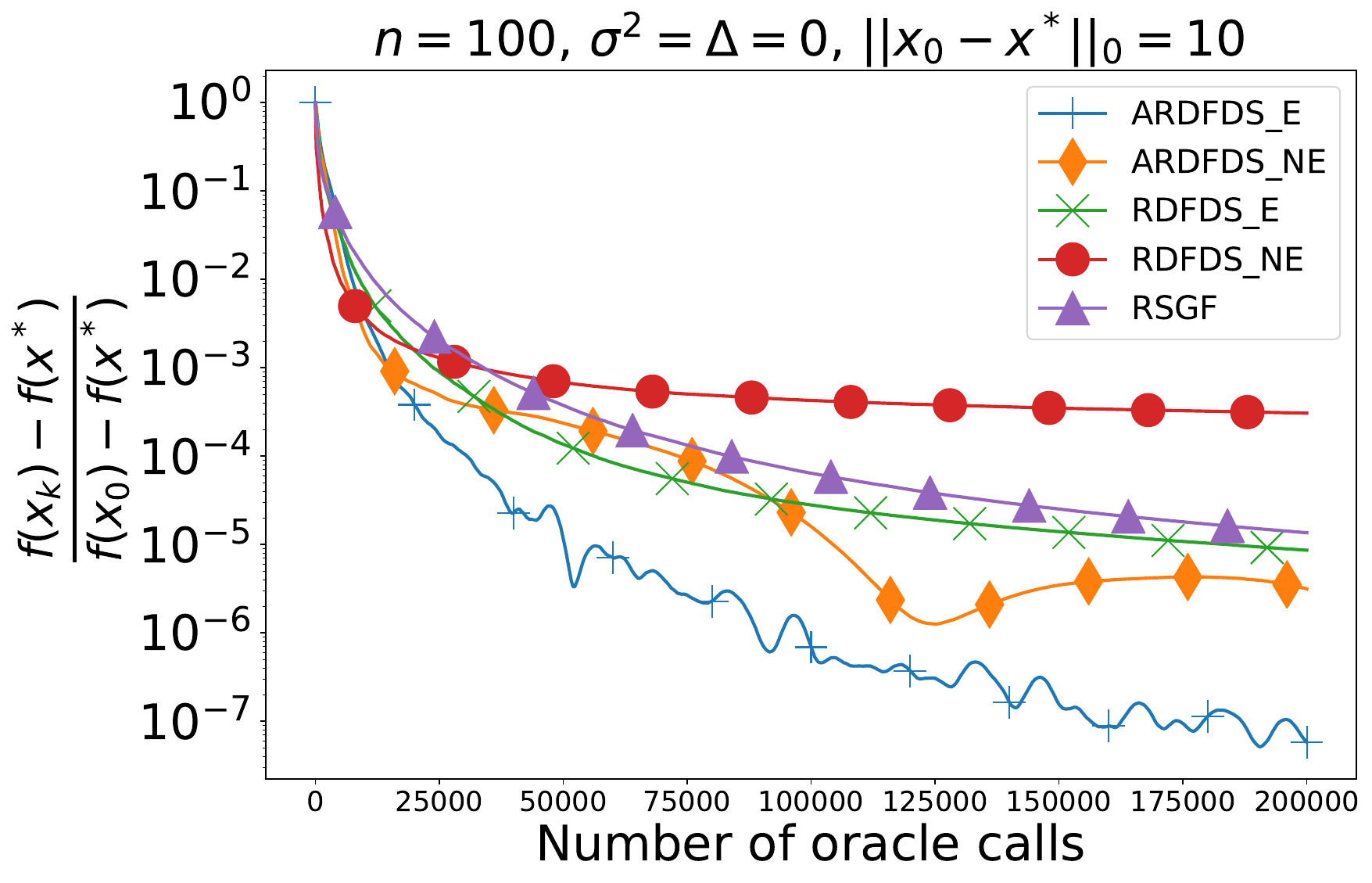}
\includegraphics[scale=0.14]{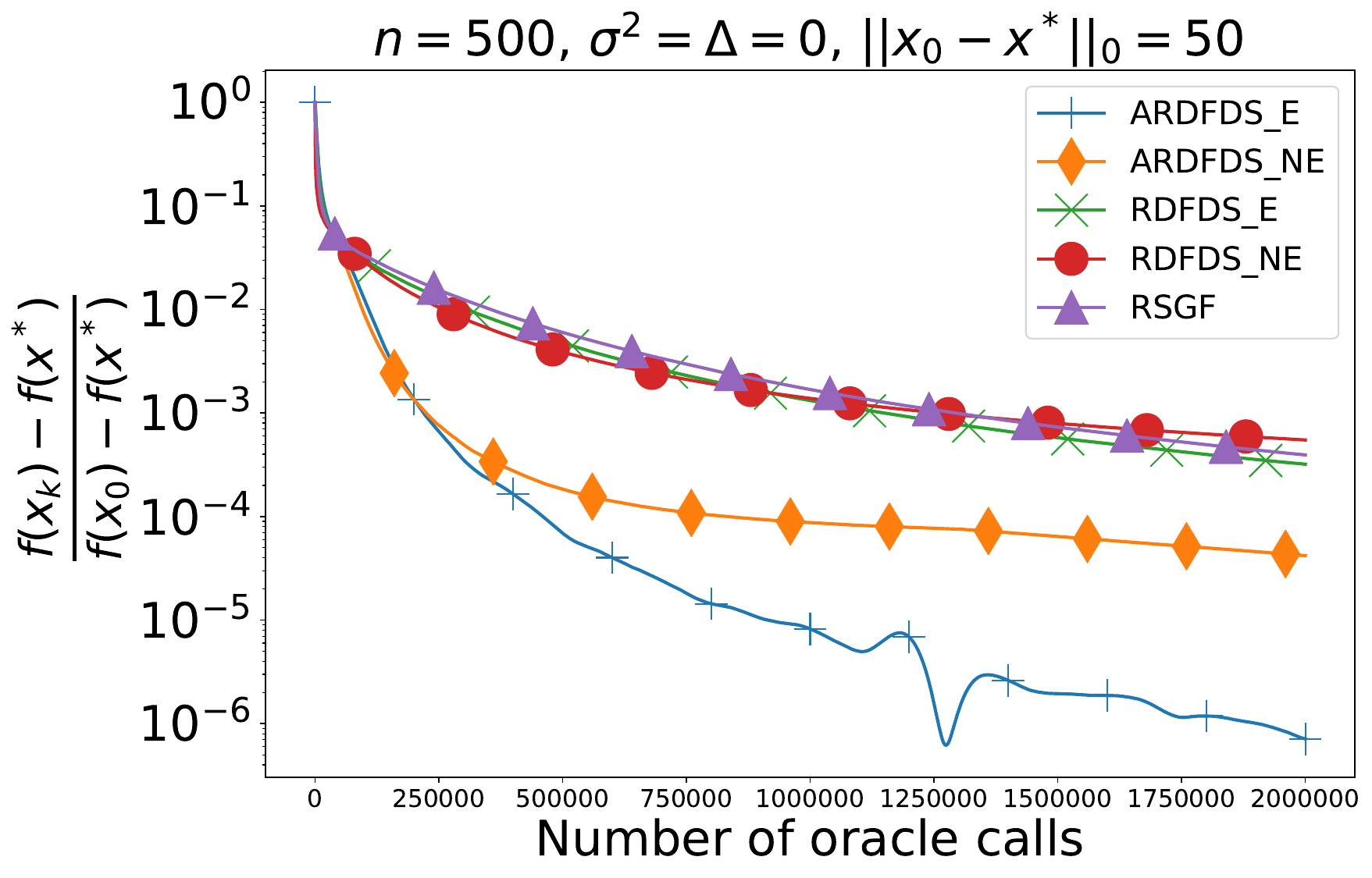}
\includegraphics[scale=0.14]{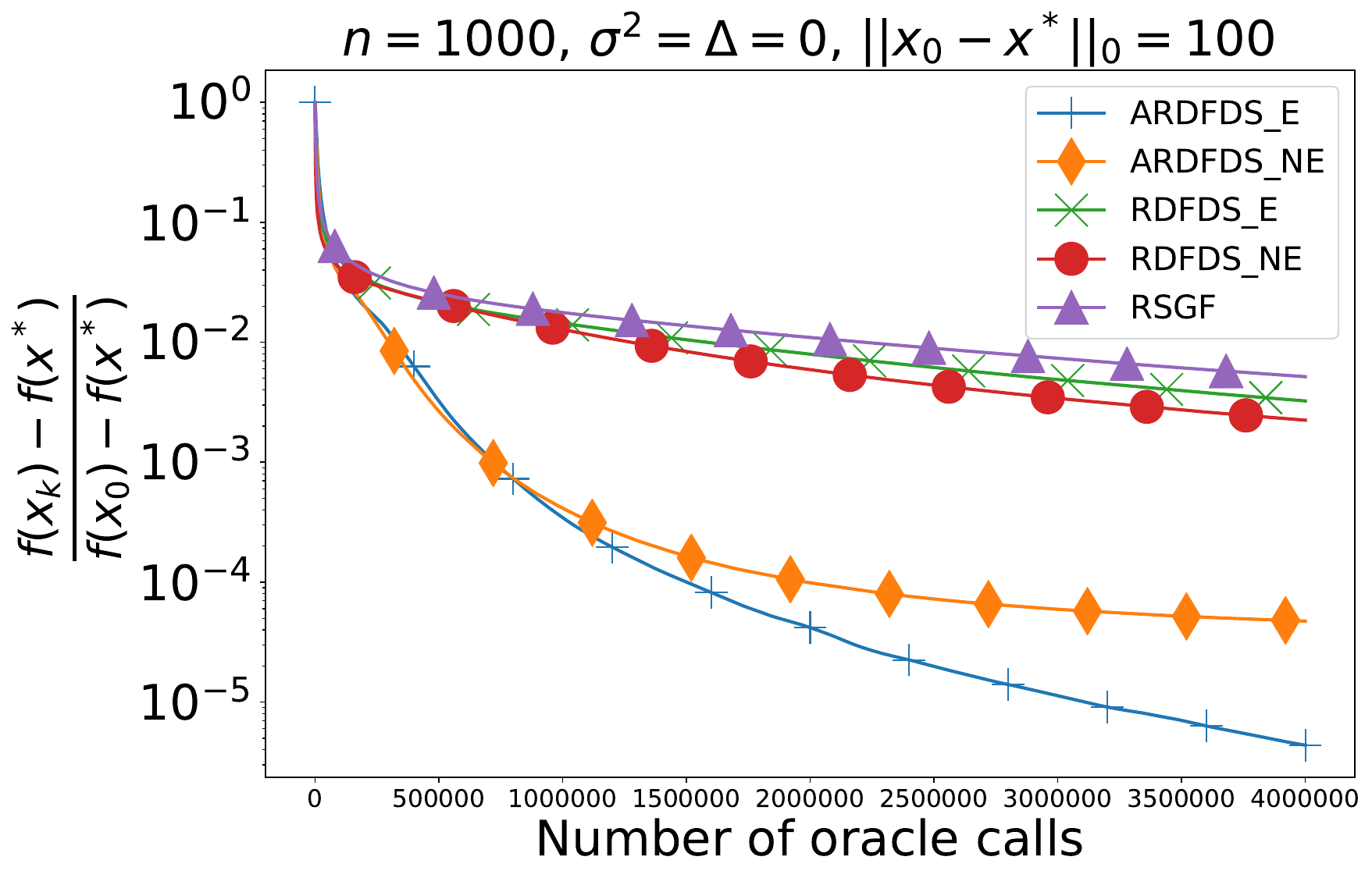}
\includegraphics[scale=0.14]{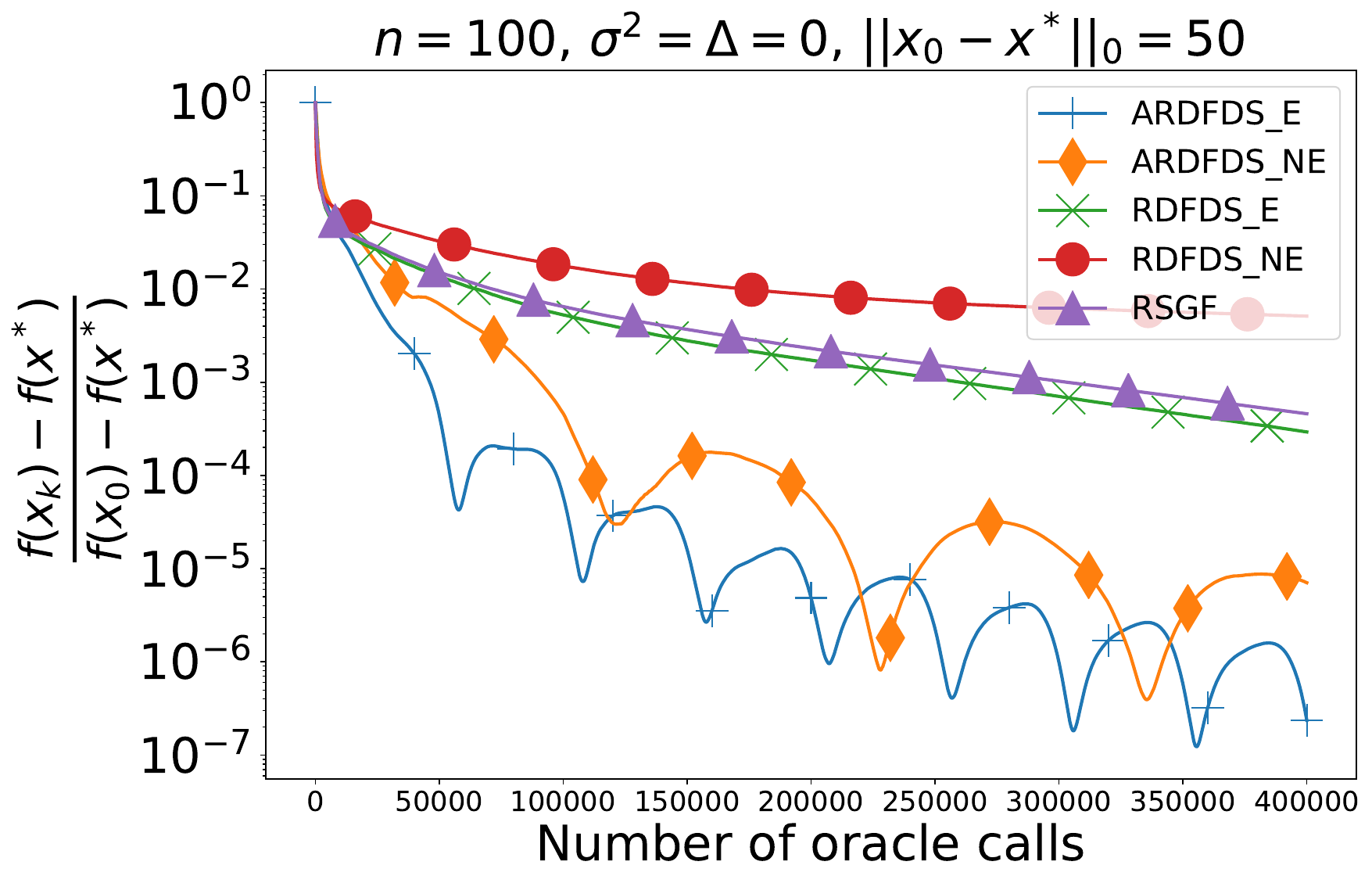}
\includegraphics[scale=0.14]{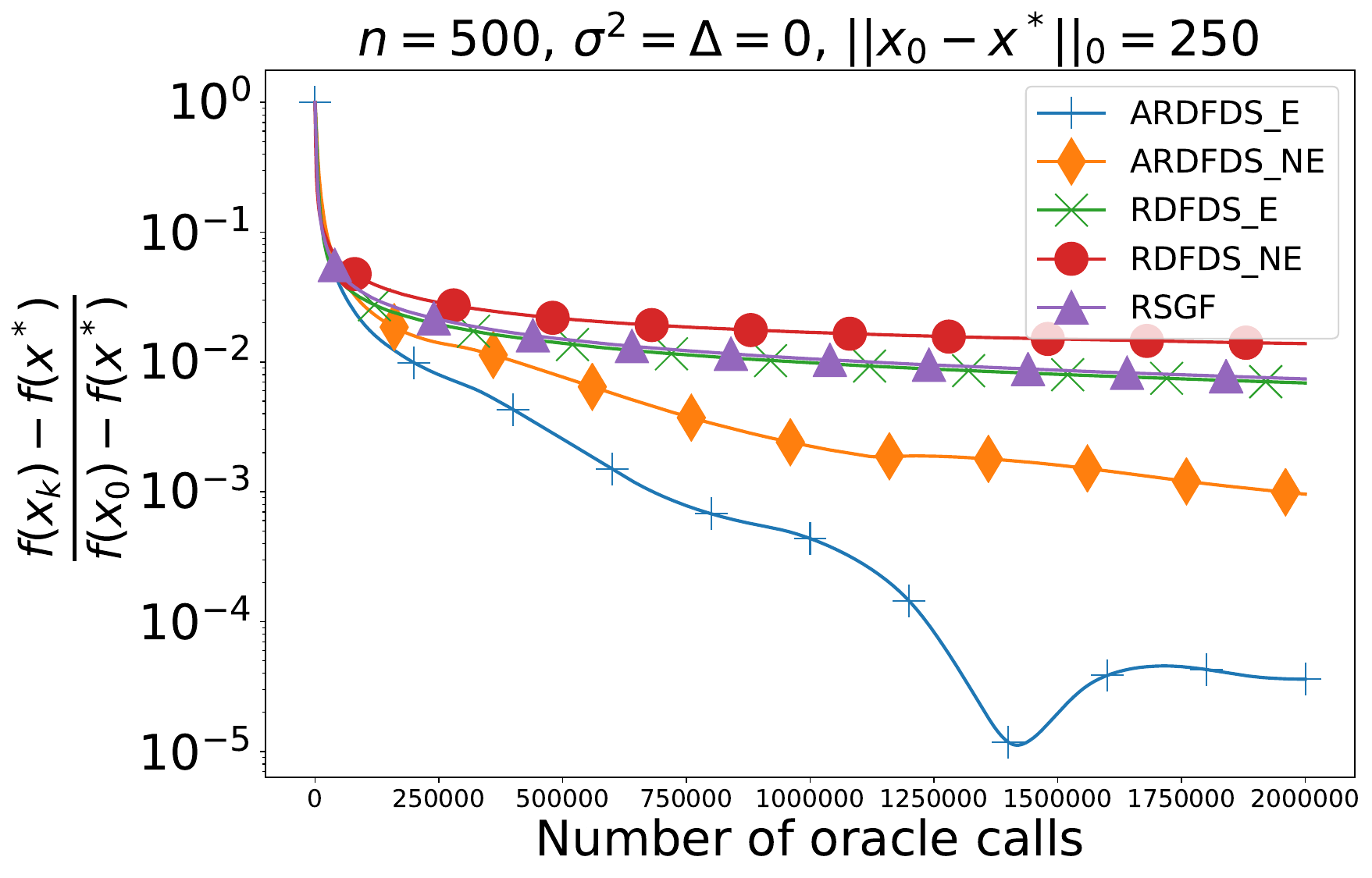}
\includegraphics[scale=0.14]{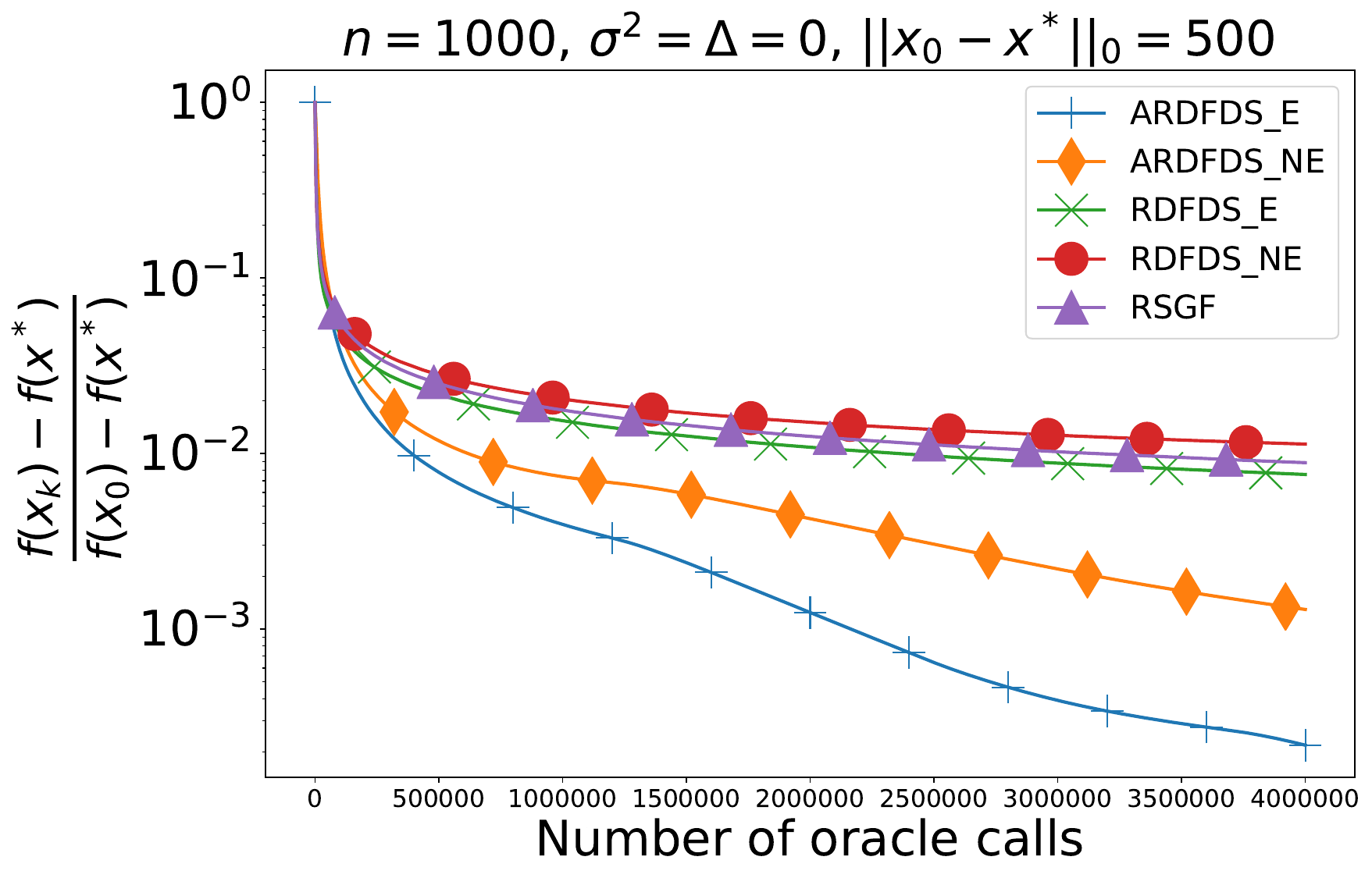}
\caption{\newstuff{Numerical results for minimizing Nesterov's function for different sparsity levels of $x_0-x^*$ and dimensions $n$ of the problem.}}
\label{fig:nesterov_sparsity_experiment}
\end{figure}
As the theory predicts\pd{, our methods with $p = 1$ work increasingly better than our methods with $p = 2$ as} $n$ is growing when $\|x_0 - x^*\|_0$ is small.}

\subsubsection{Experiments with different variance} 
\label{sec:exp_nesterov_sigma}
\newstuff{In this subsection we report the numerical results with different values of $\sigma^2$. For each choice of the dimension $n$ we used two values of $\sigma^2$: $\sigma_{\text{small}}^2 = \tfrac{\e^{\nicefrac{3}{2}}\sqrt{nL_2}}{\|x_0-x^*\|_1}$ and $\sigma_{\text{big}}^2 = 10000\sigma_{\text{small}}^2$ with $\e = 10^{-3}$. As one can see from Table~\ref{tbl:ARDFDS_params}, when $\sigma^2 = \sigma_{\text{small}}^2$ the first term under the maximum in the complexity bound is dominating (up to logarithmic factors). This implies that ARDFDS with $p=1$ is guaranteed to find an $\e$-solution even with the mini-batch size $m=1$ (up to logarithmical factors). We choose $\e = 10^{-3}$, $\Delta = 0$ and $x_0$ such that it differs from $x^*$ only in the first component and run the experiments for $n = 100, 500, 1000$ and $\sigma^2 = \sigma_{\text{small}}^2, \sigma_{\text{big}}^2$, see Figures~\ref{fig:nesterov_small_sigma}~and~\ref{fig:nesterov_big_sigma}.
\begin{figure}
\centering 
\includegraphics[scale=0.14]{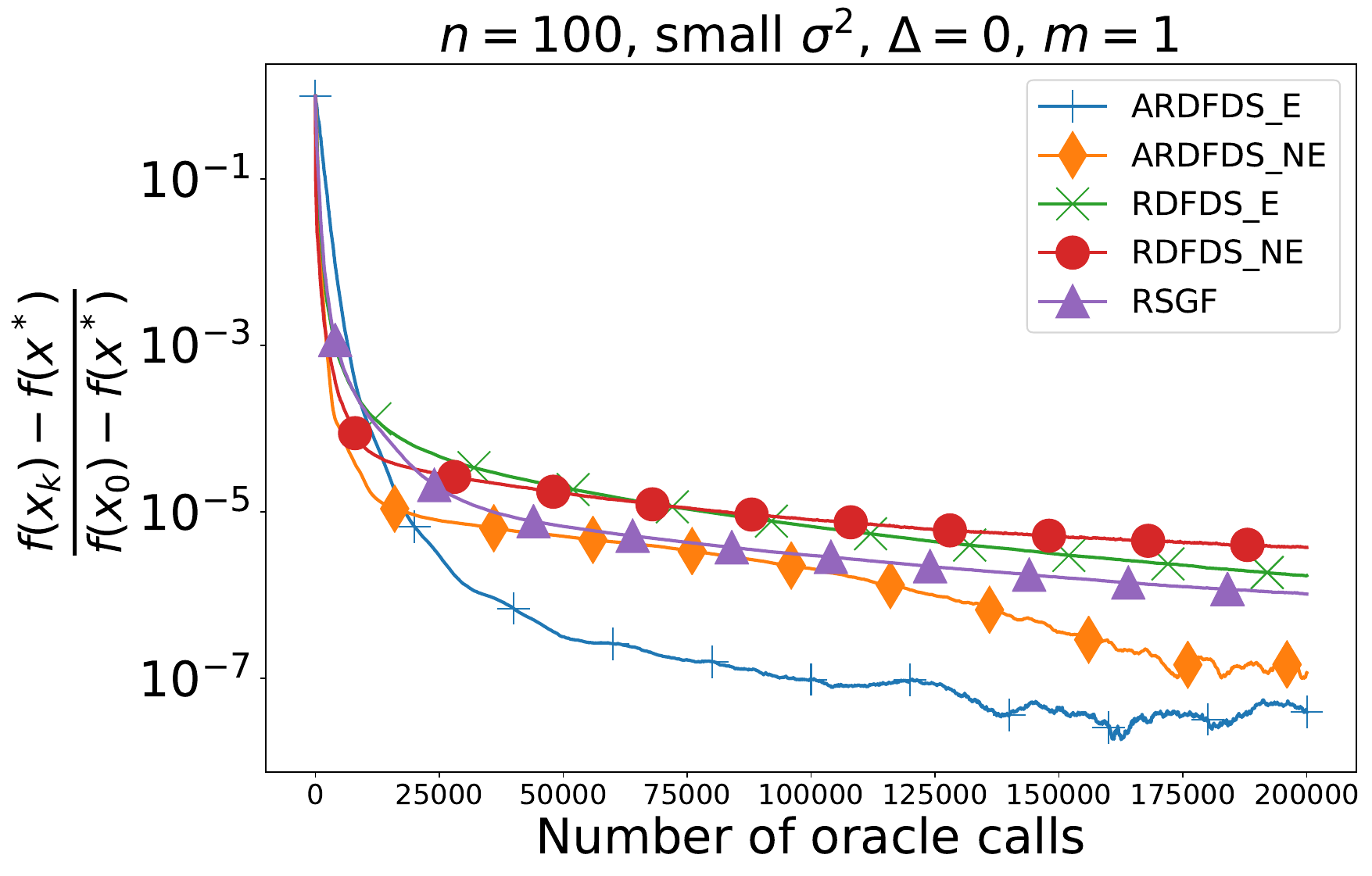}
\includegraphics[scale=0.14]{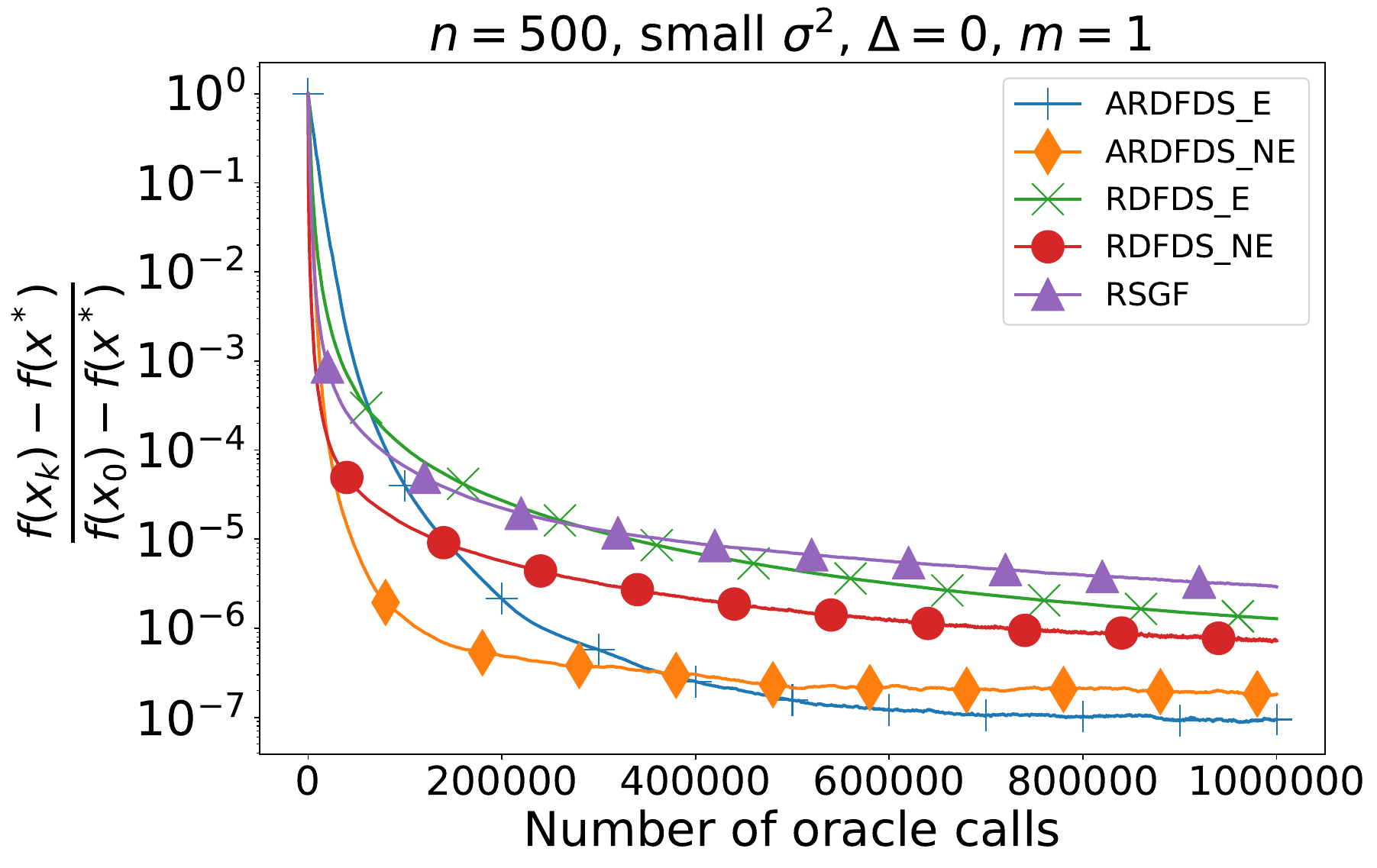}
\includegraphics[scale=0.14]{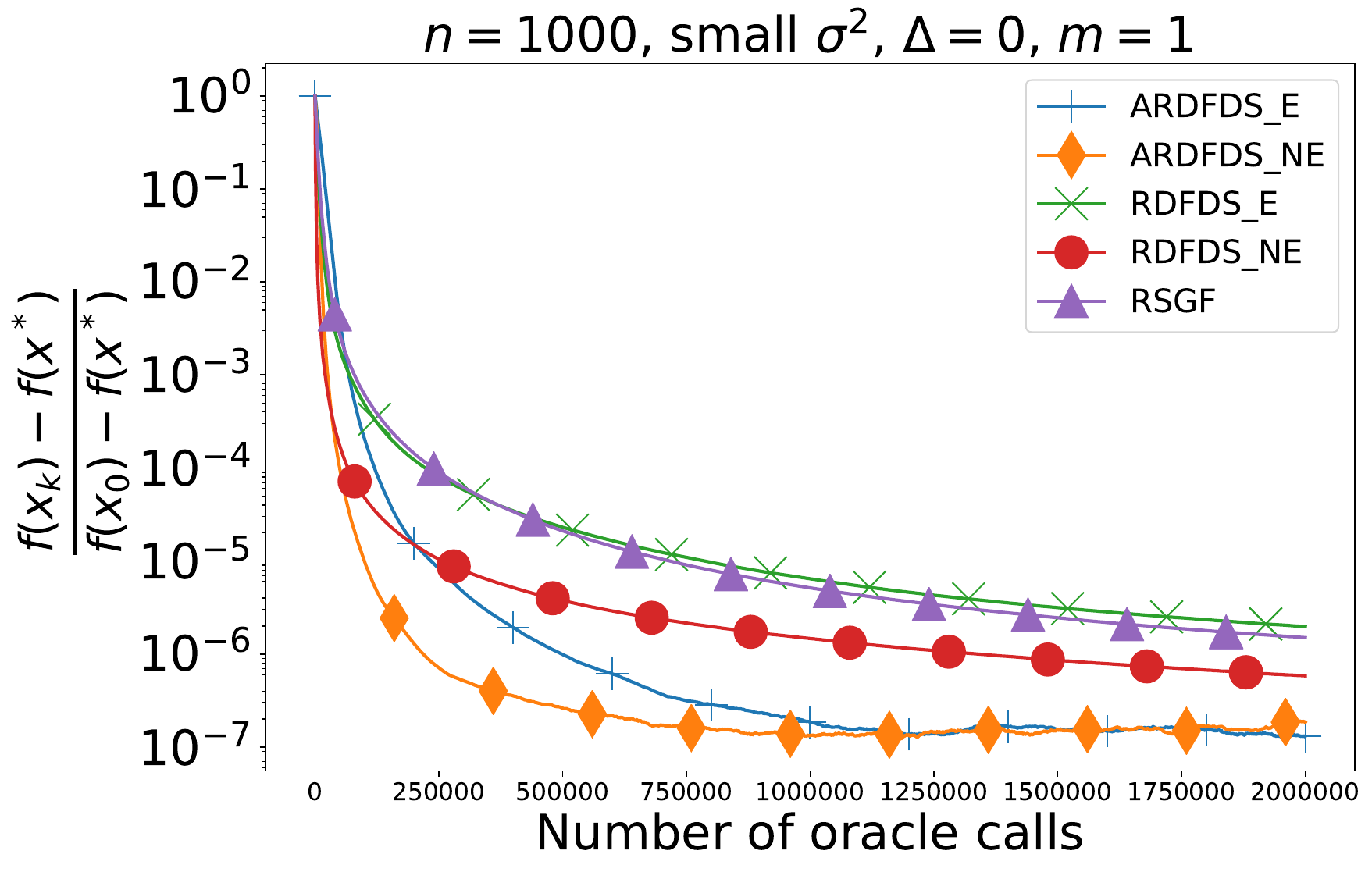}
\includegraphics[scale=0.14]{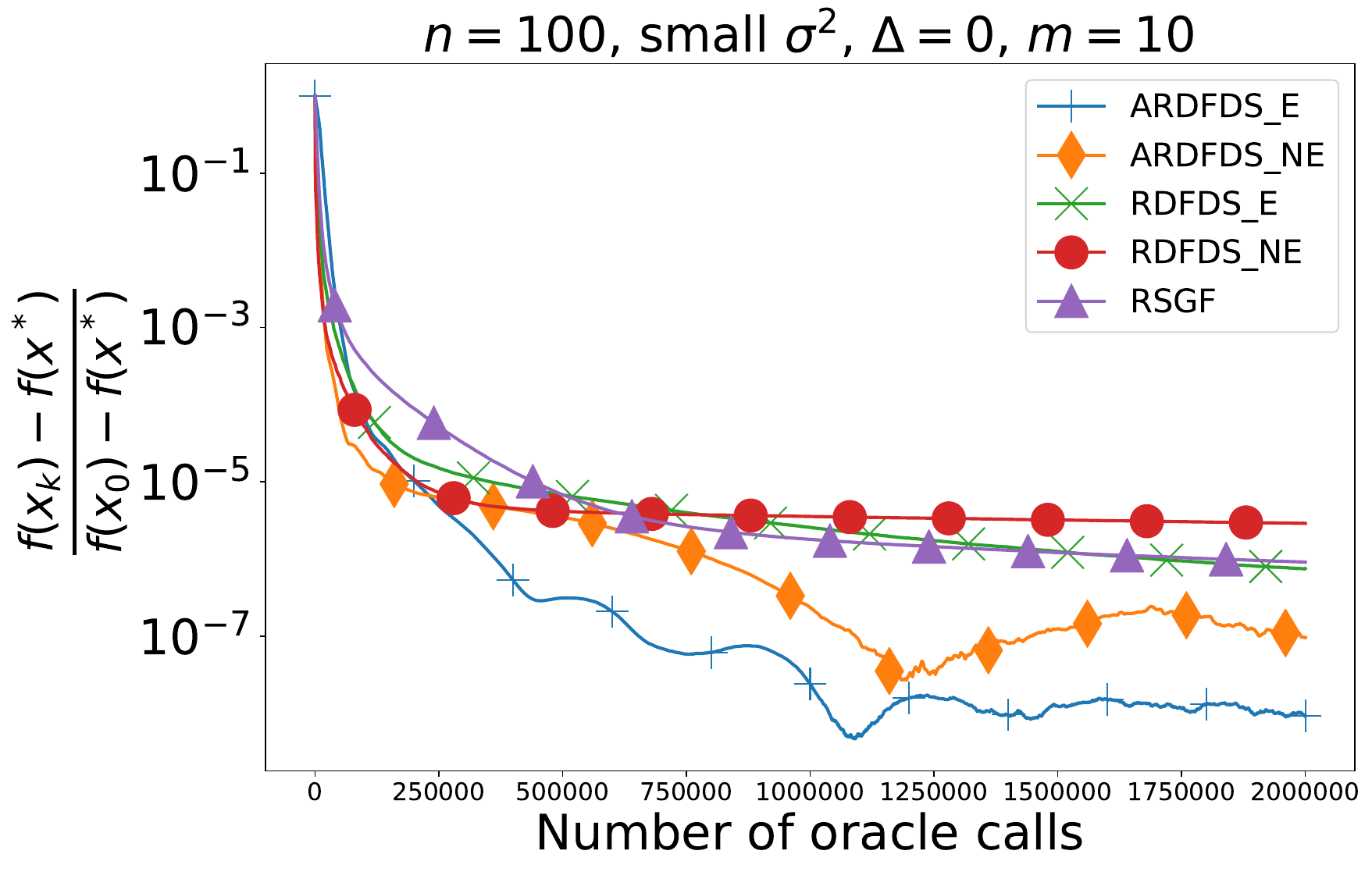}
\includegraphics[scale=0.14]{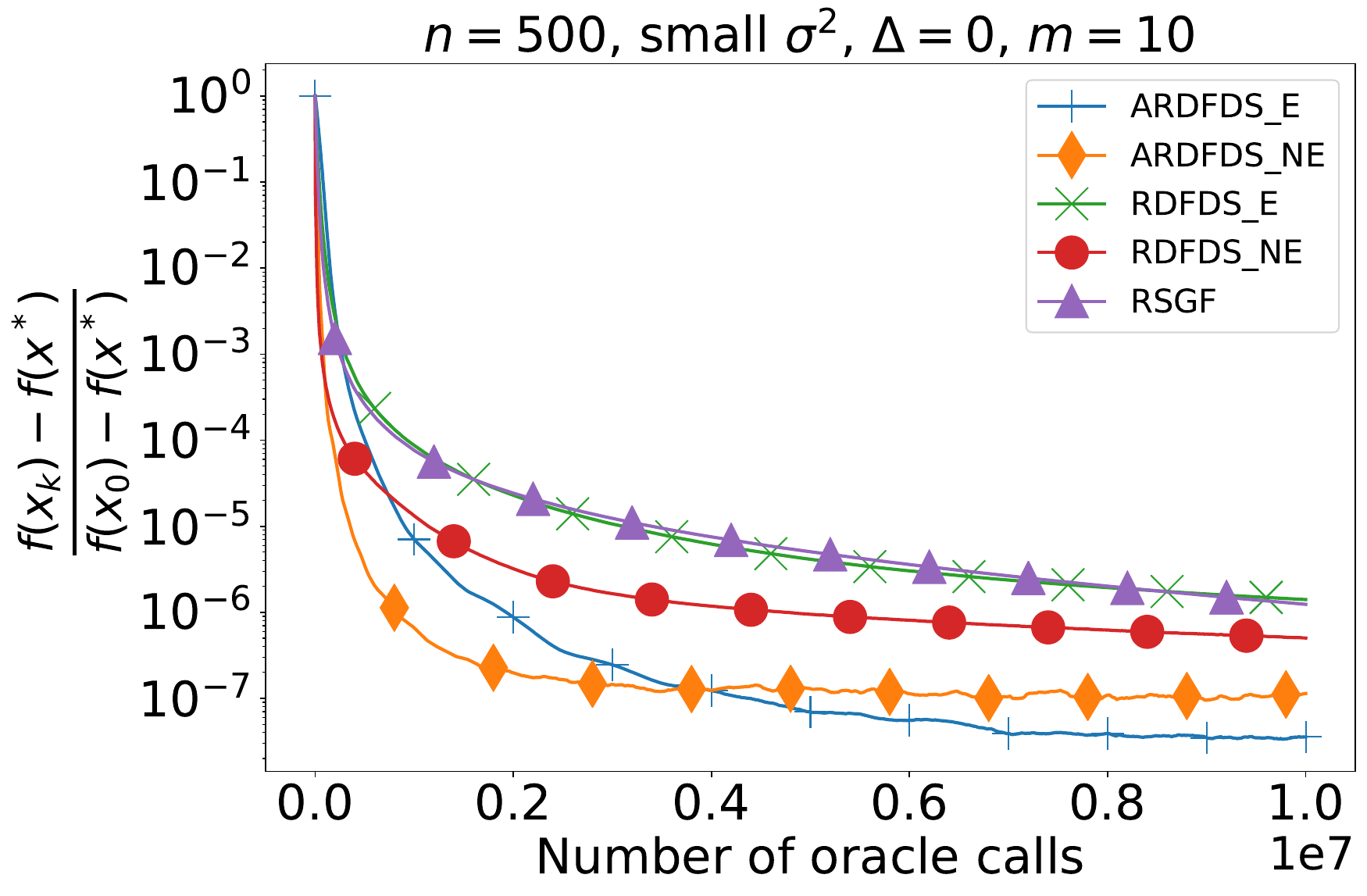}
\includegraphics[scale=0.14]{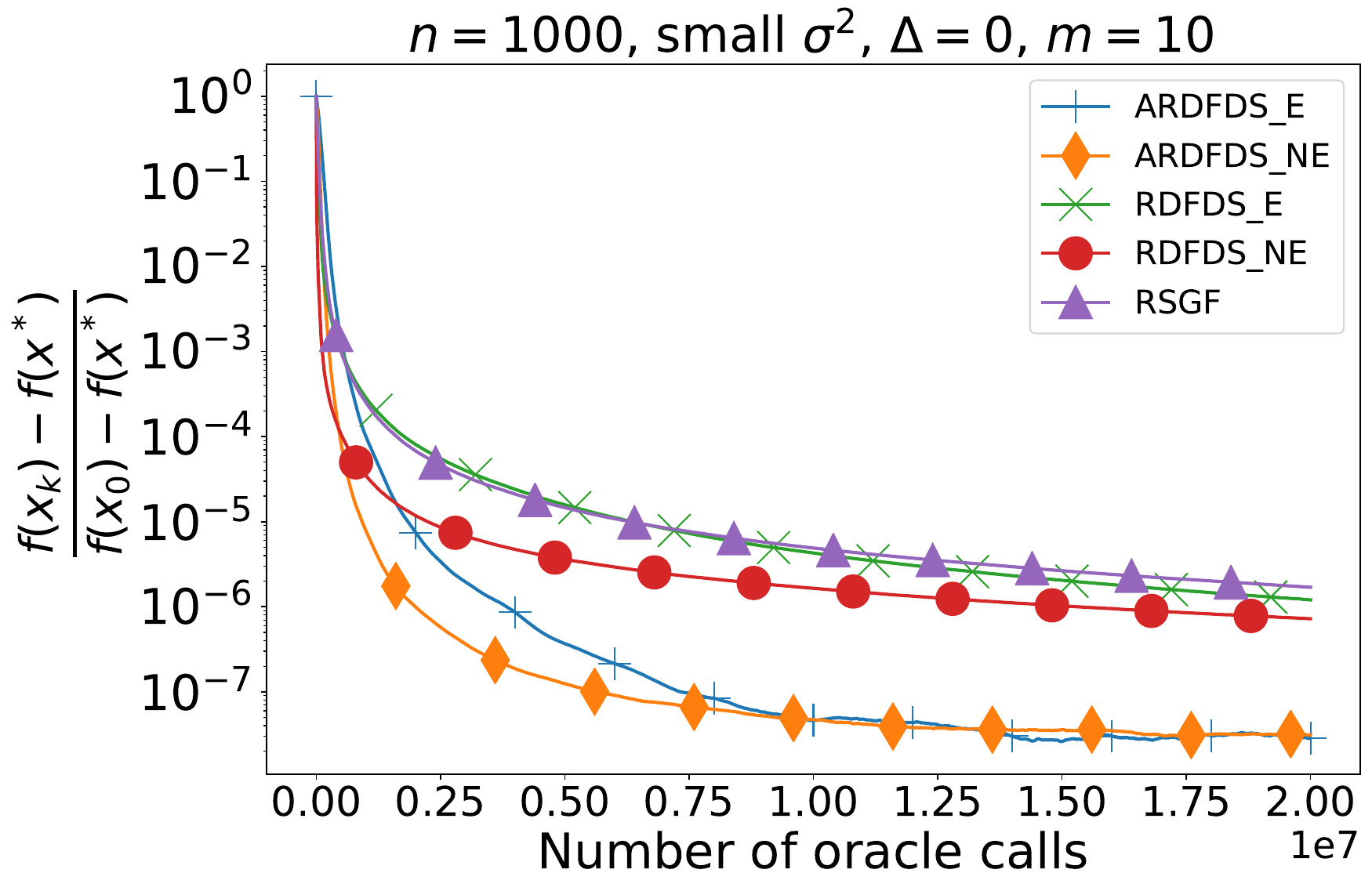}
\includegraphics[scale=0.14]{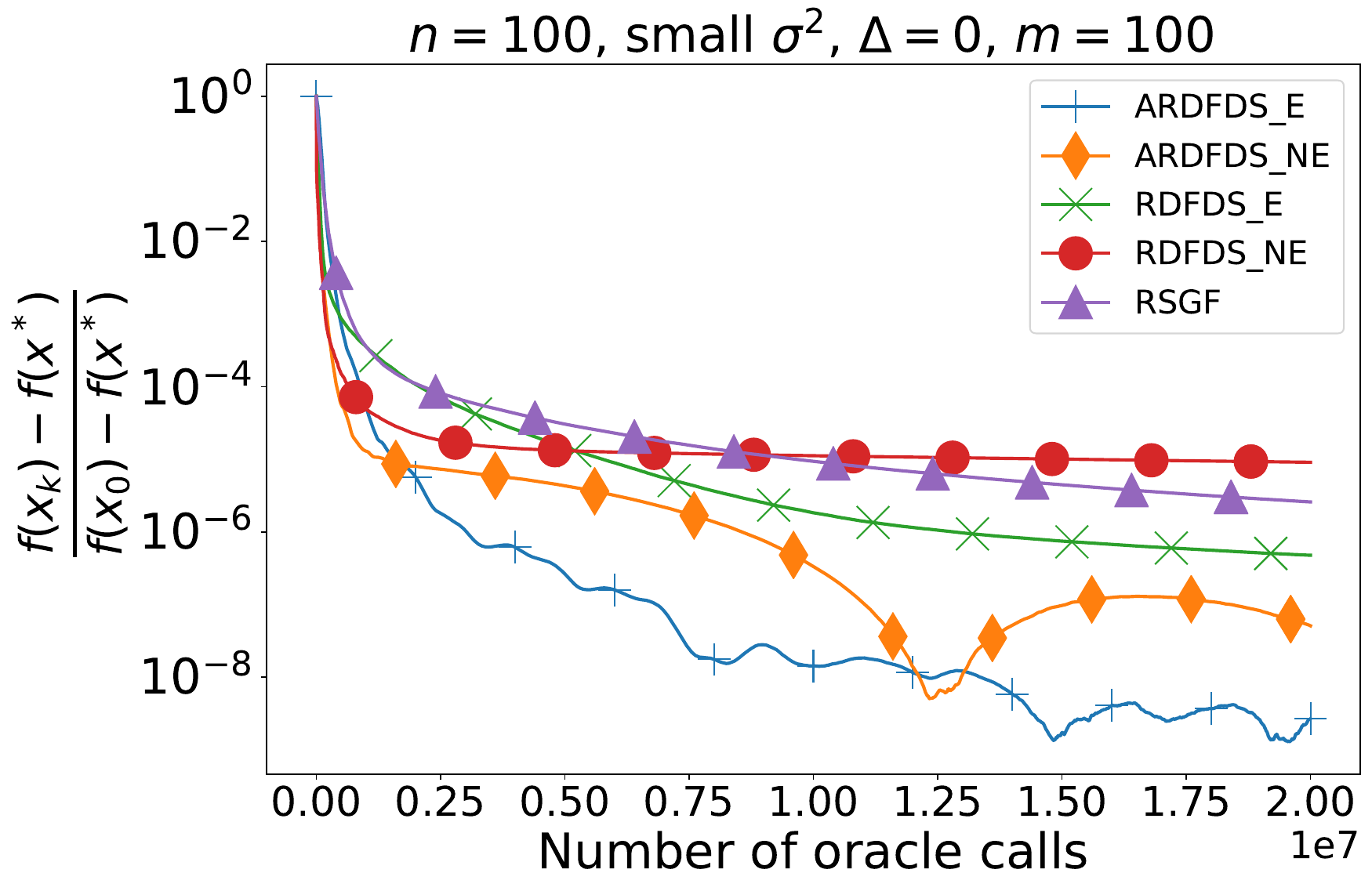}
\includegraphics[scale=0.14]{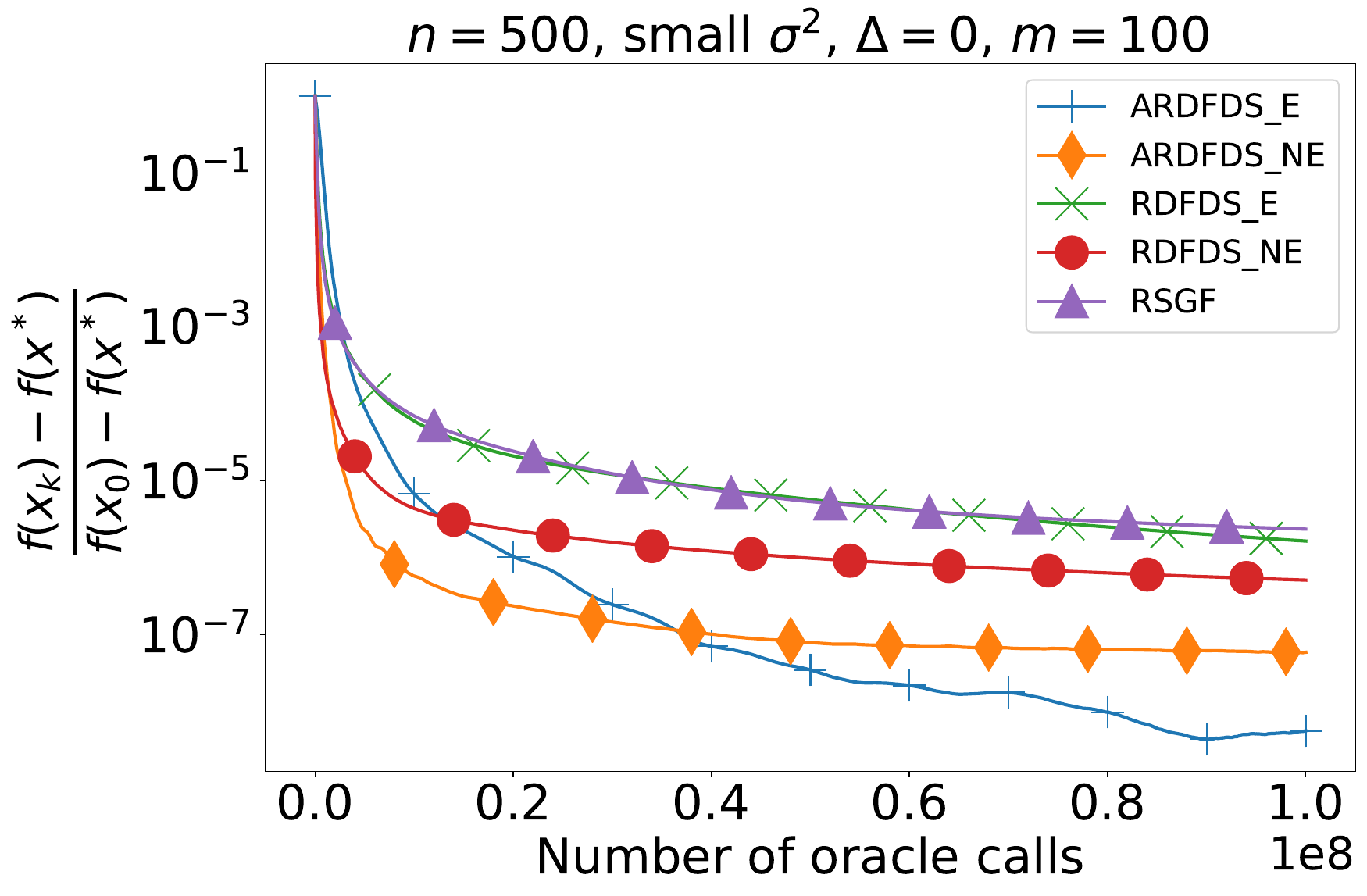}
\includegraphics[scale=0.14]{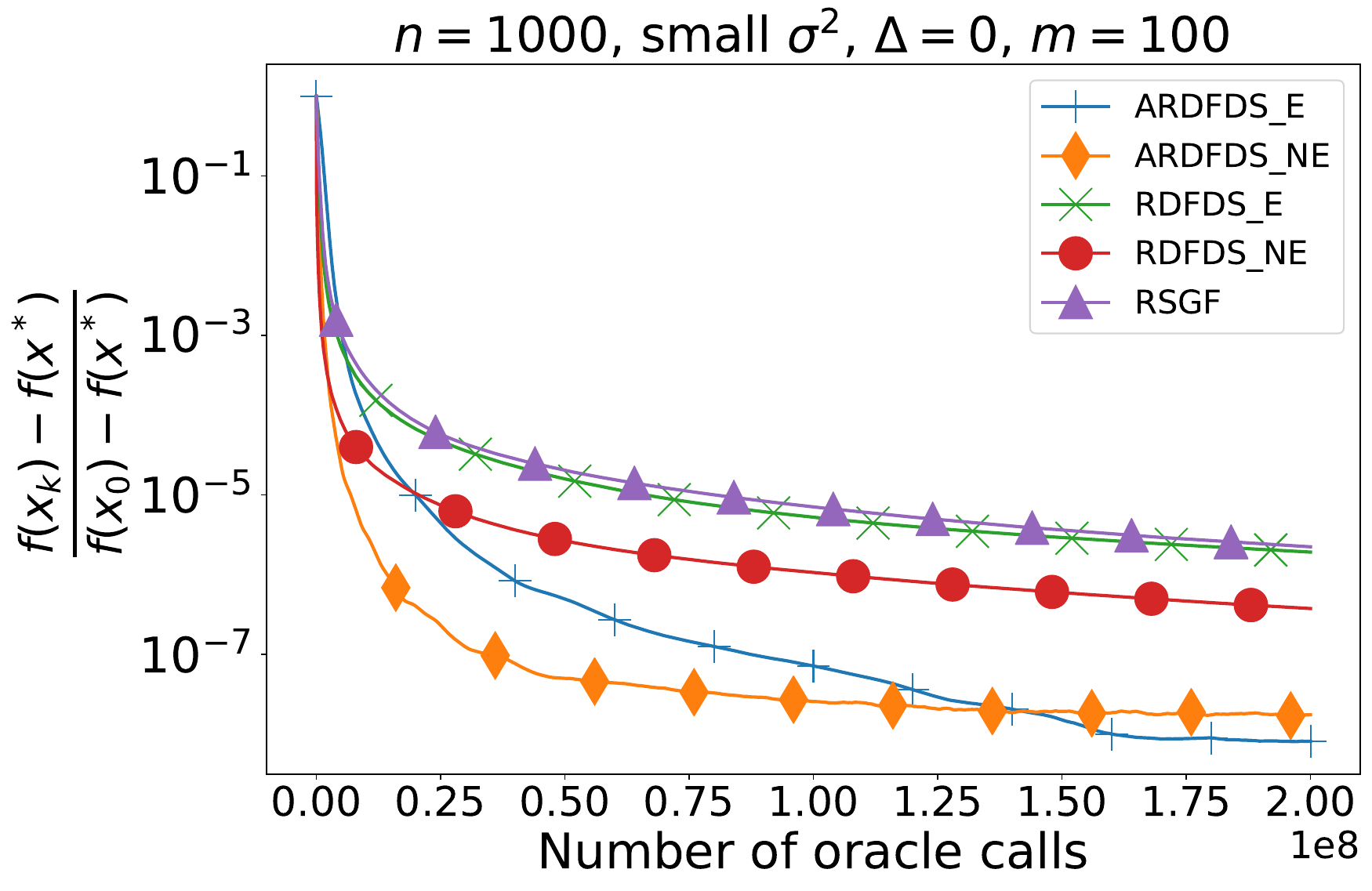}
\caption{\newstuff{Numerical results for minimizing Nesterov's function with noisy stochastic oracle having $\sigma^2 = \sigma_{\text{small}}^2$ for different sizes of mini-batch $m$ and dimensions of the problem $n$.}}
\label{fig:nesterov_small_sigma}
\end{figure}
\begin{figure}
\centering 
\includegraphics[scale=0.14]{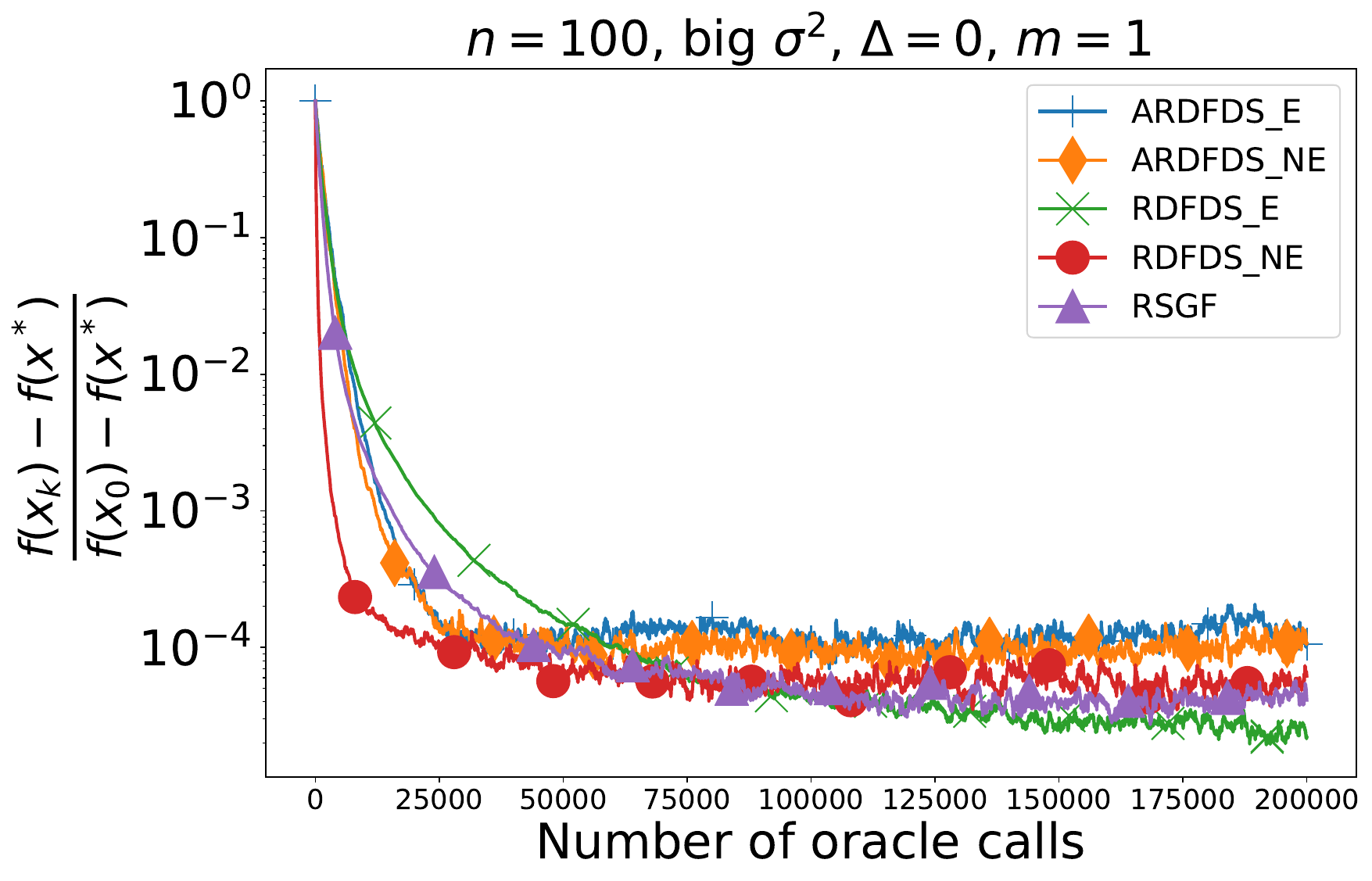}
\includegraphics[scale=0.14]{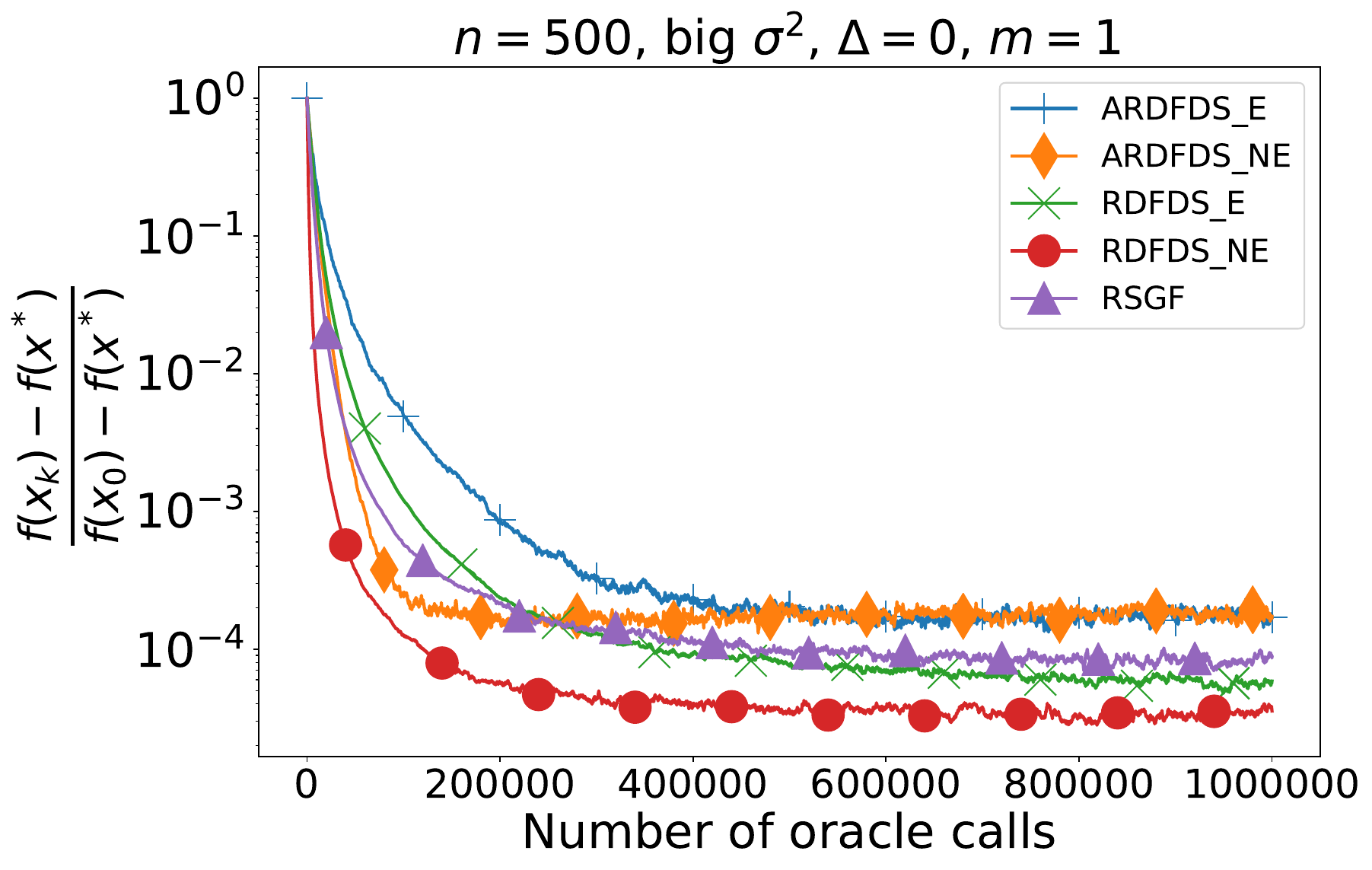}
\includegraphics[scale=0.14]{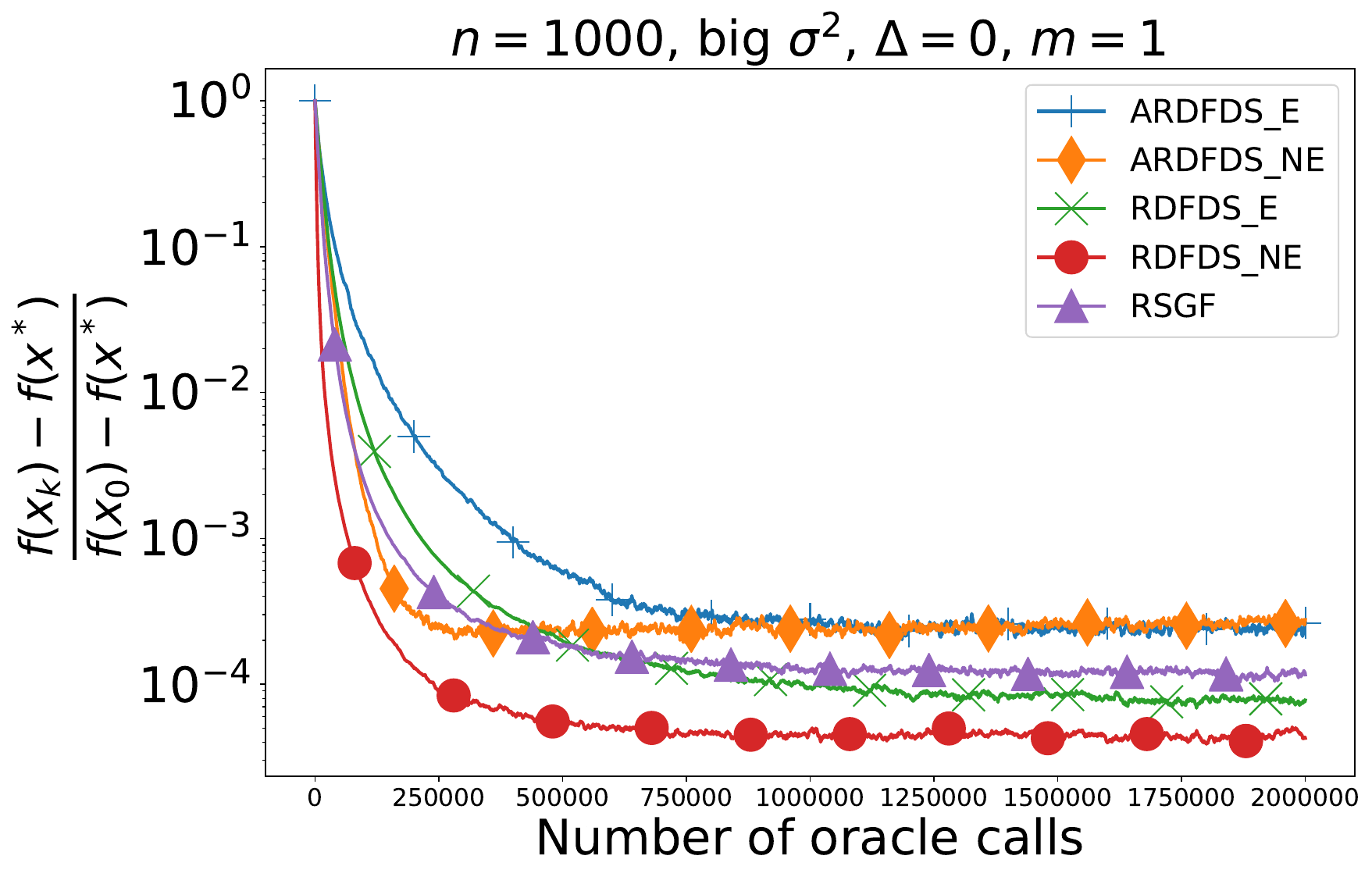}
\includegraphics[scale=0.14]{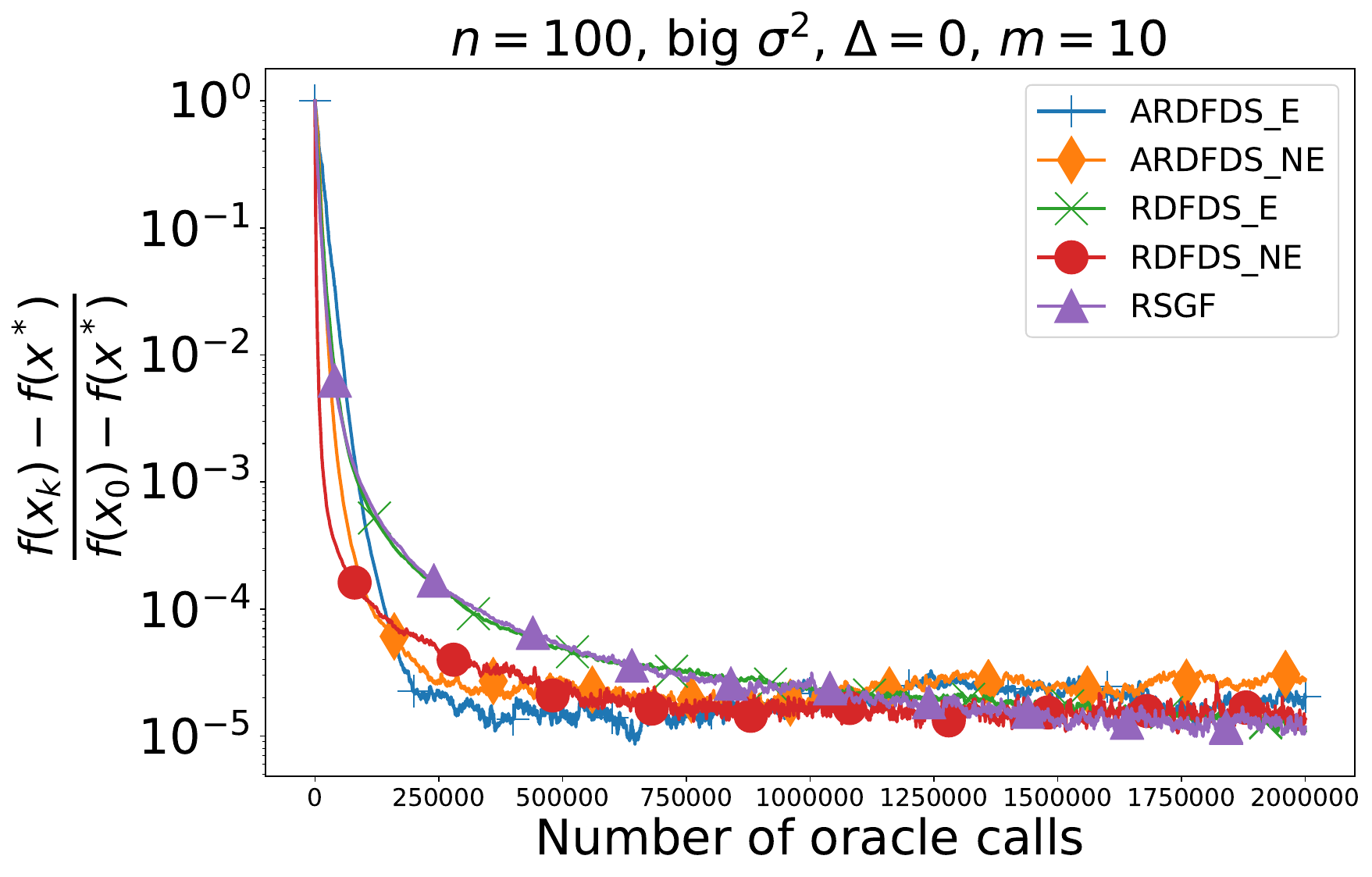}
\includegraphics[scale=0.14]{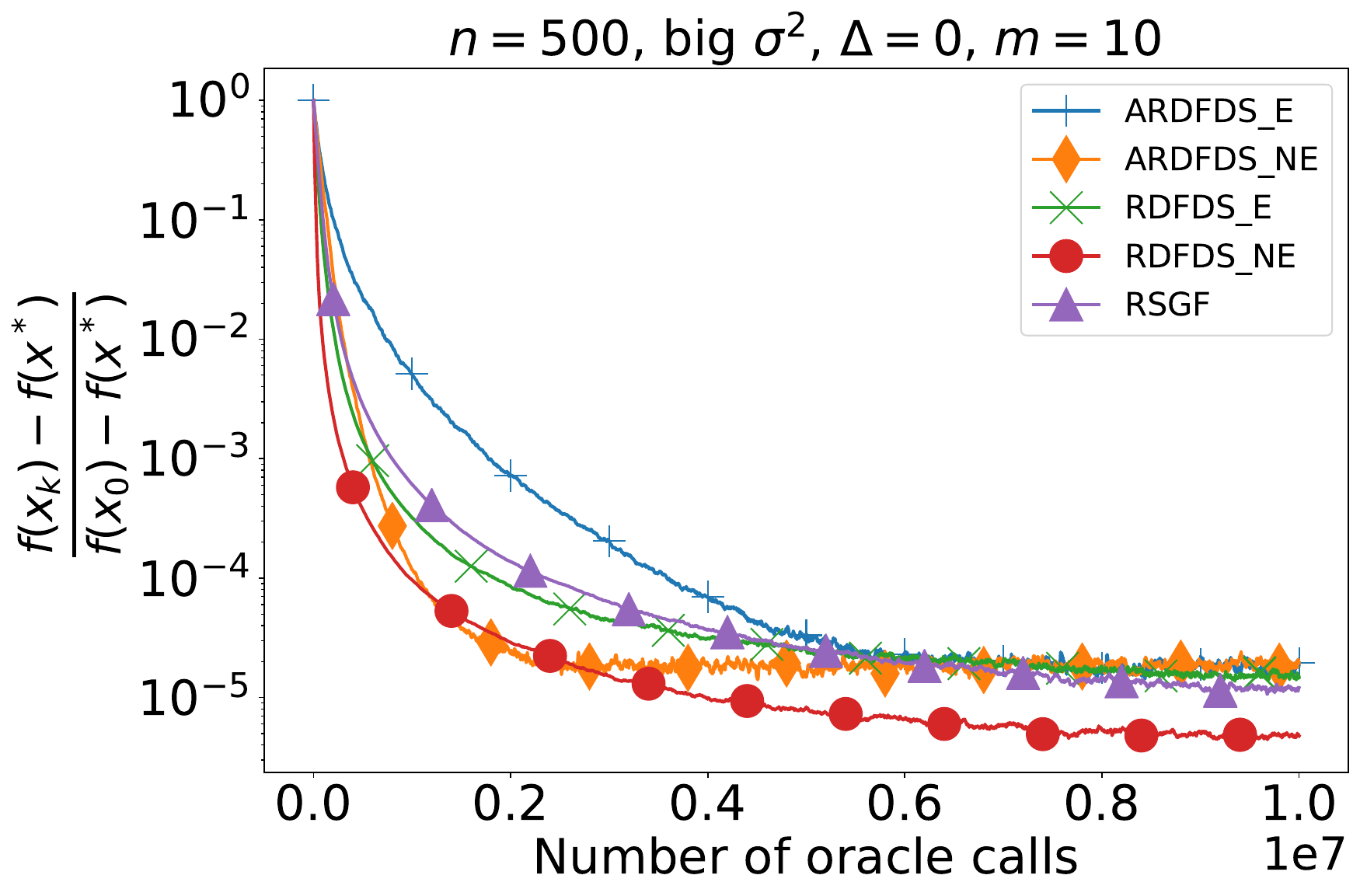}
\includegraphics[scale=0.14]{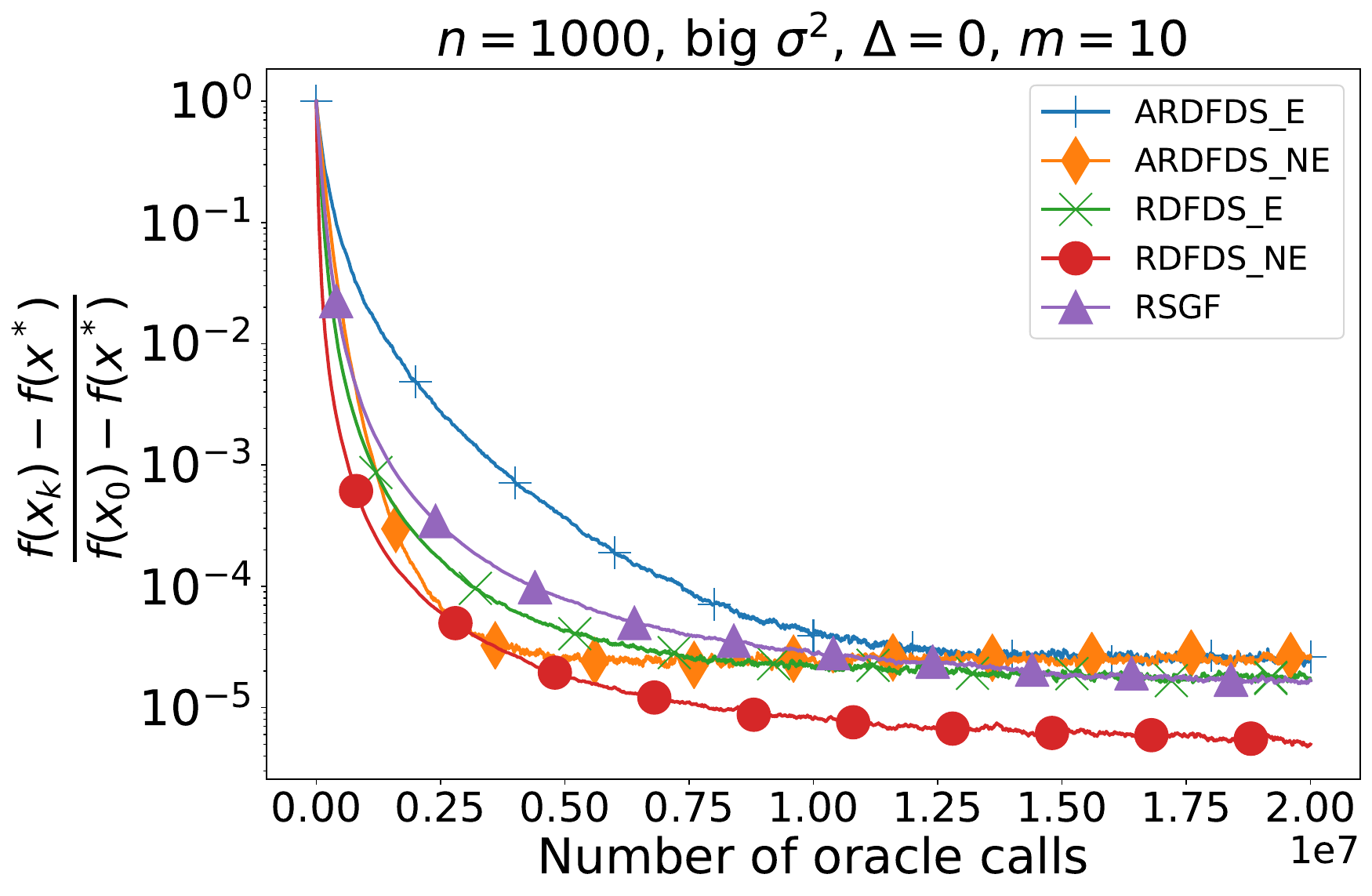}
\includegraphics[scale=0.14]{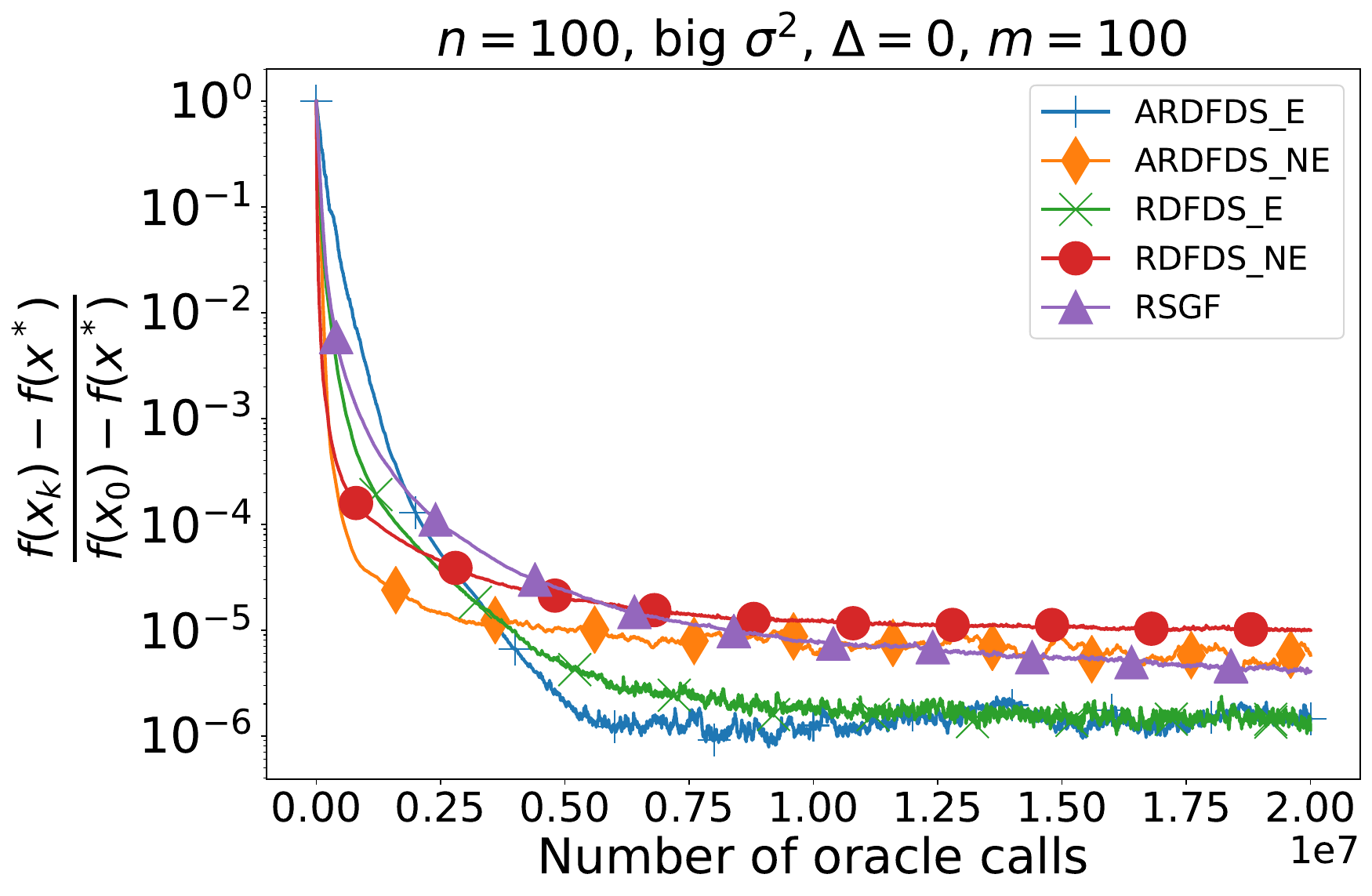}
\includegraphics[scale=0.14]{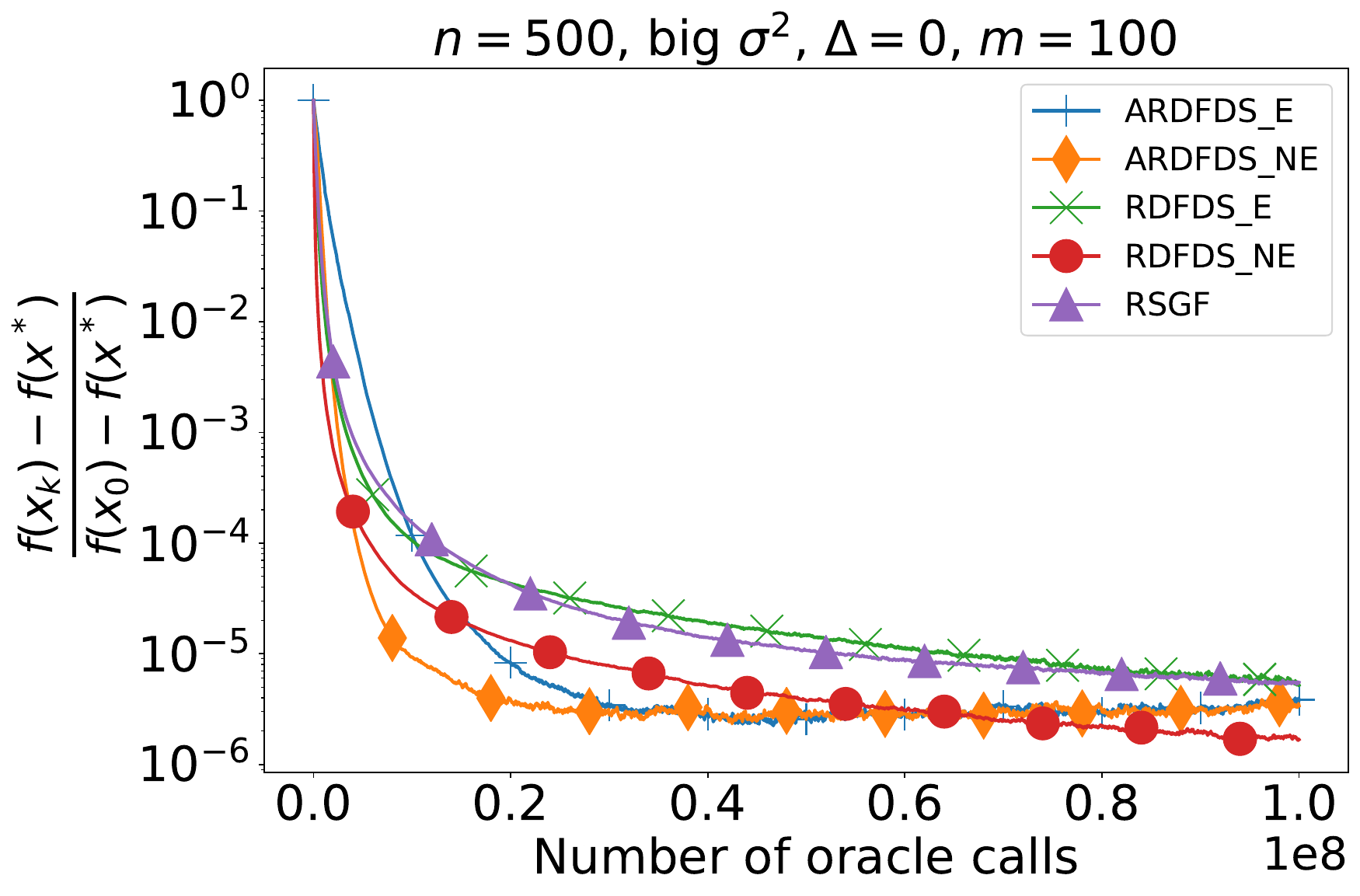}
\includegraphics[scale=0.14]{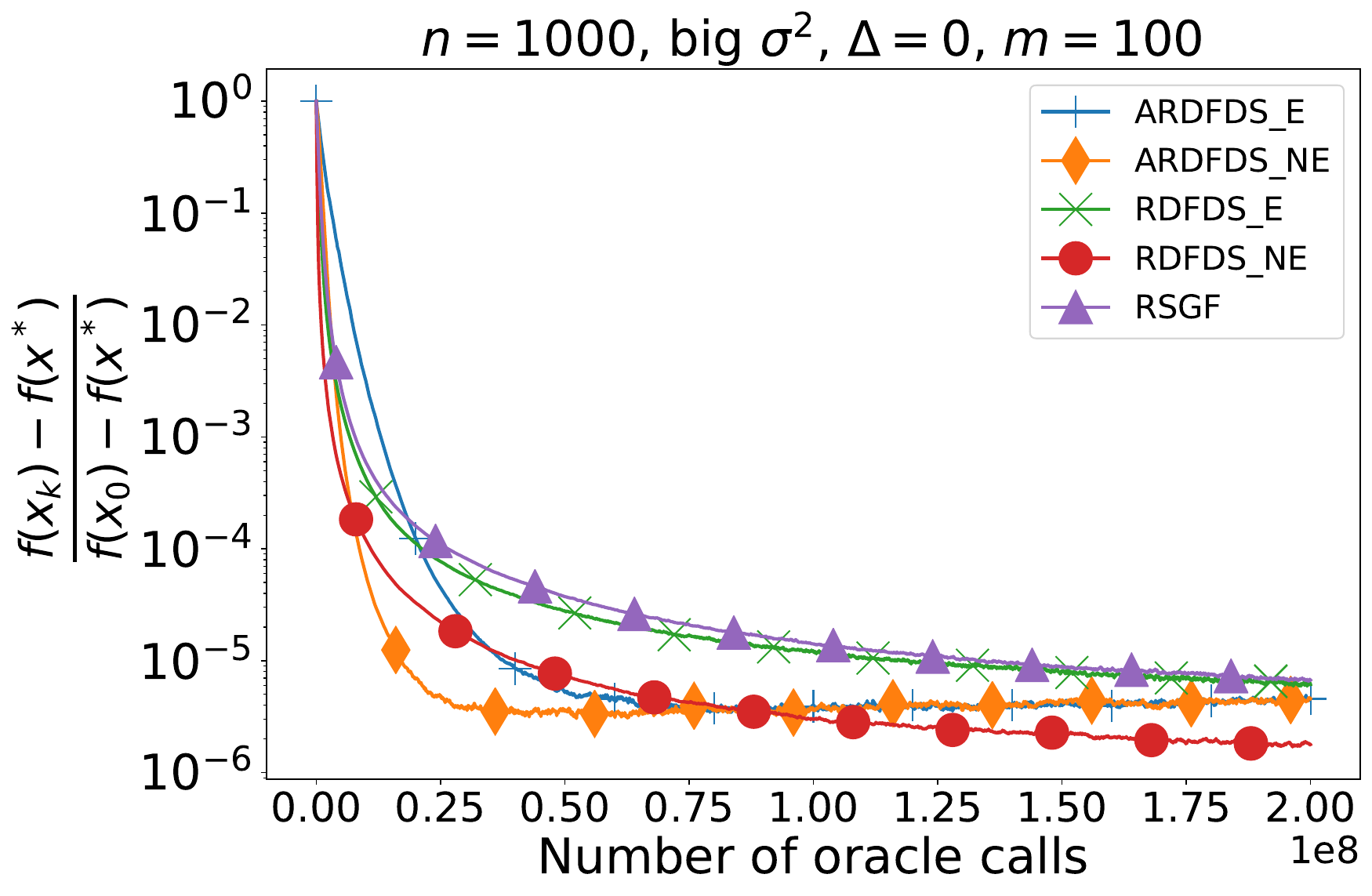}
\includegraphics[scale=0.14]{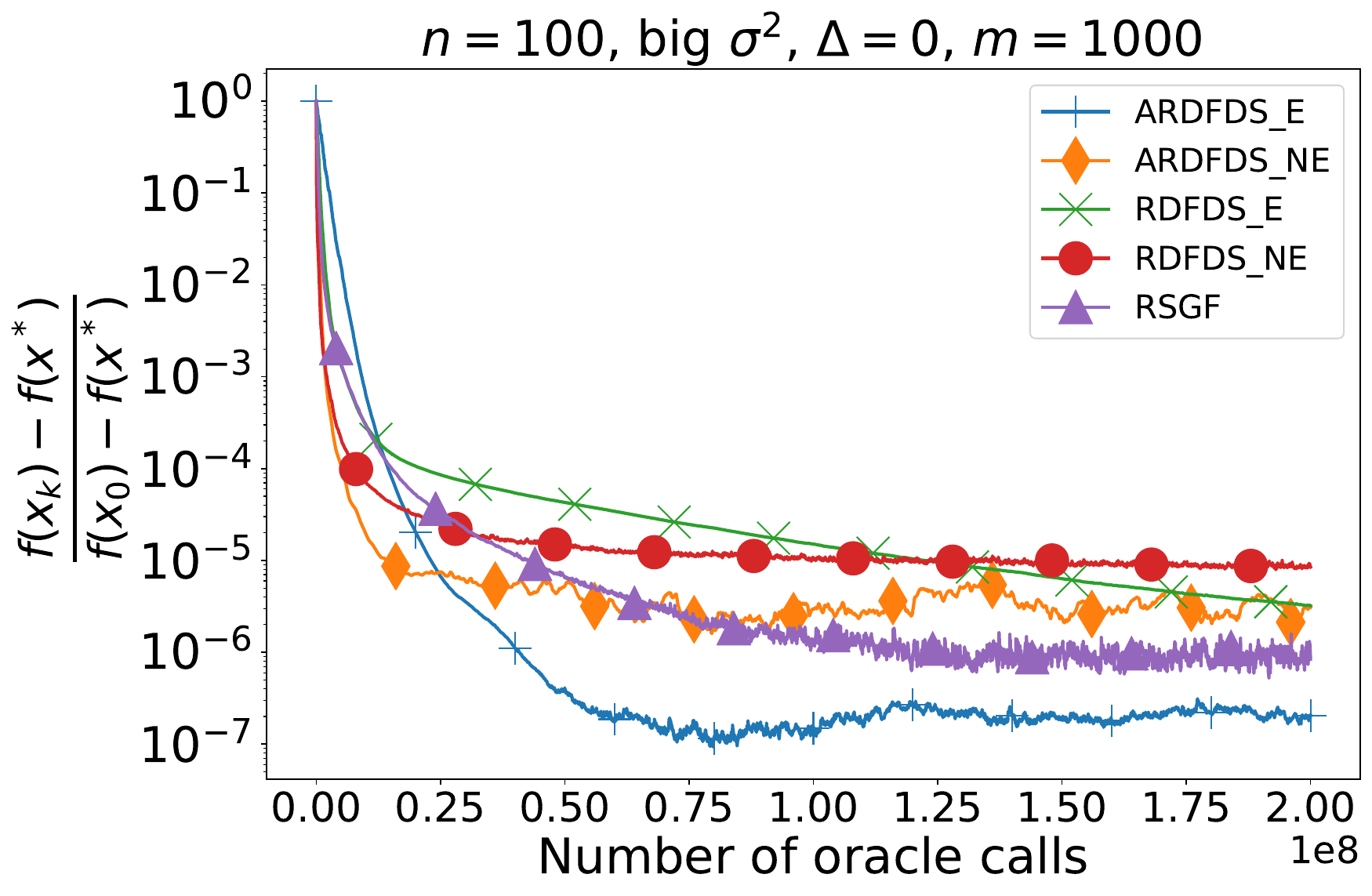}
\includegraphics[scale=0.14]{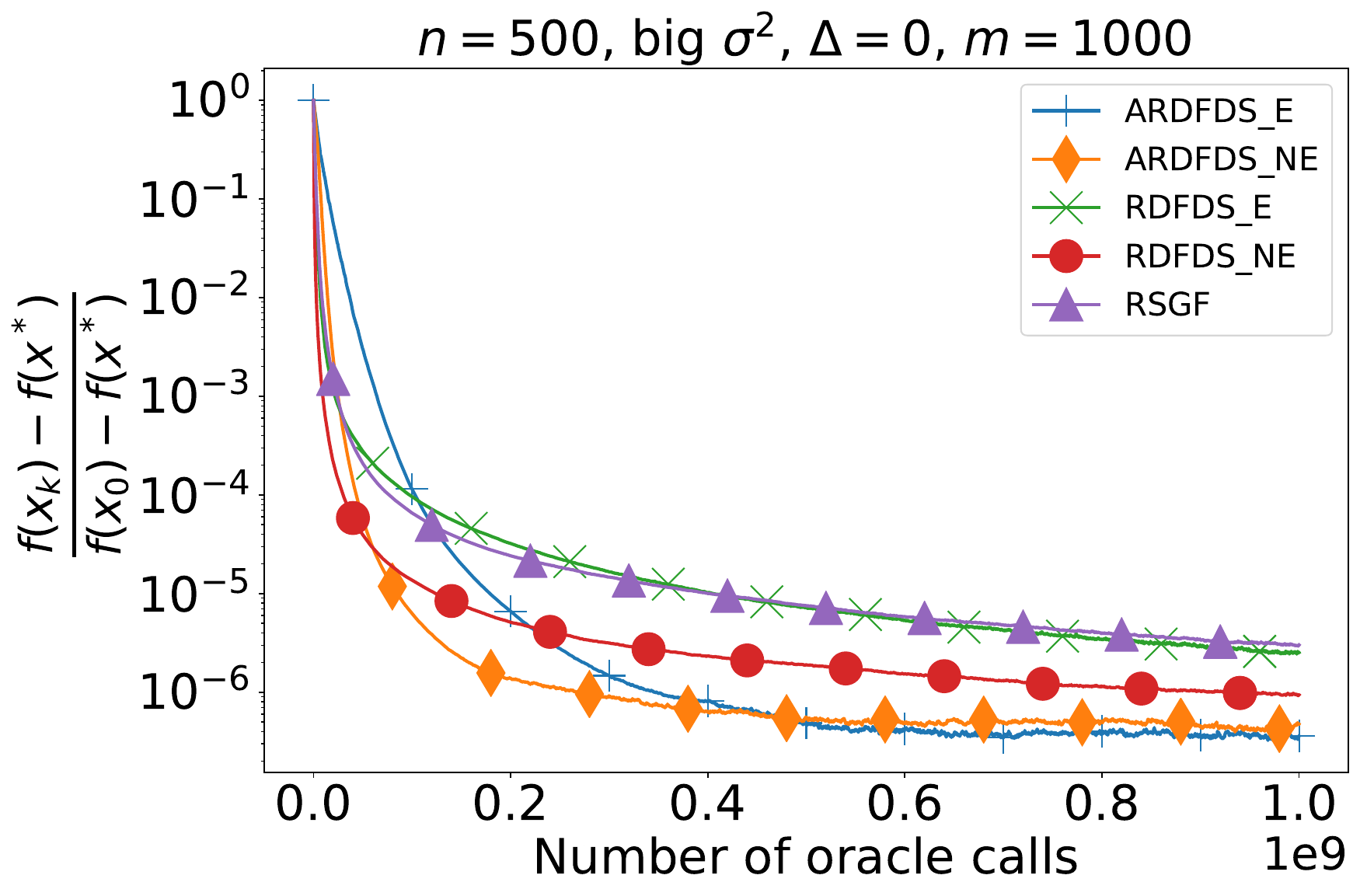}
\includegraphics[scale=0.14]{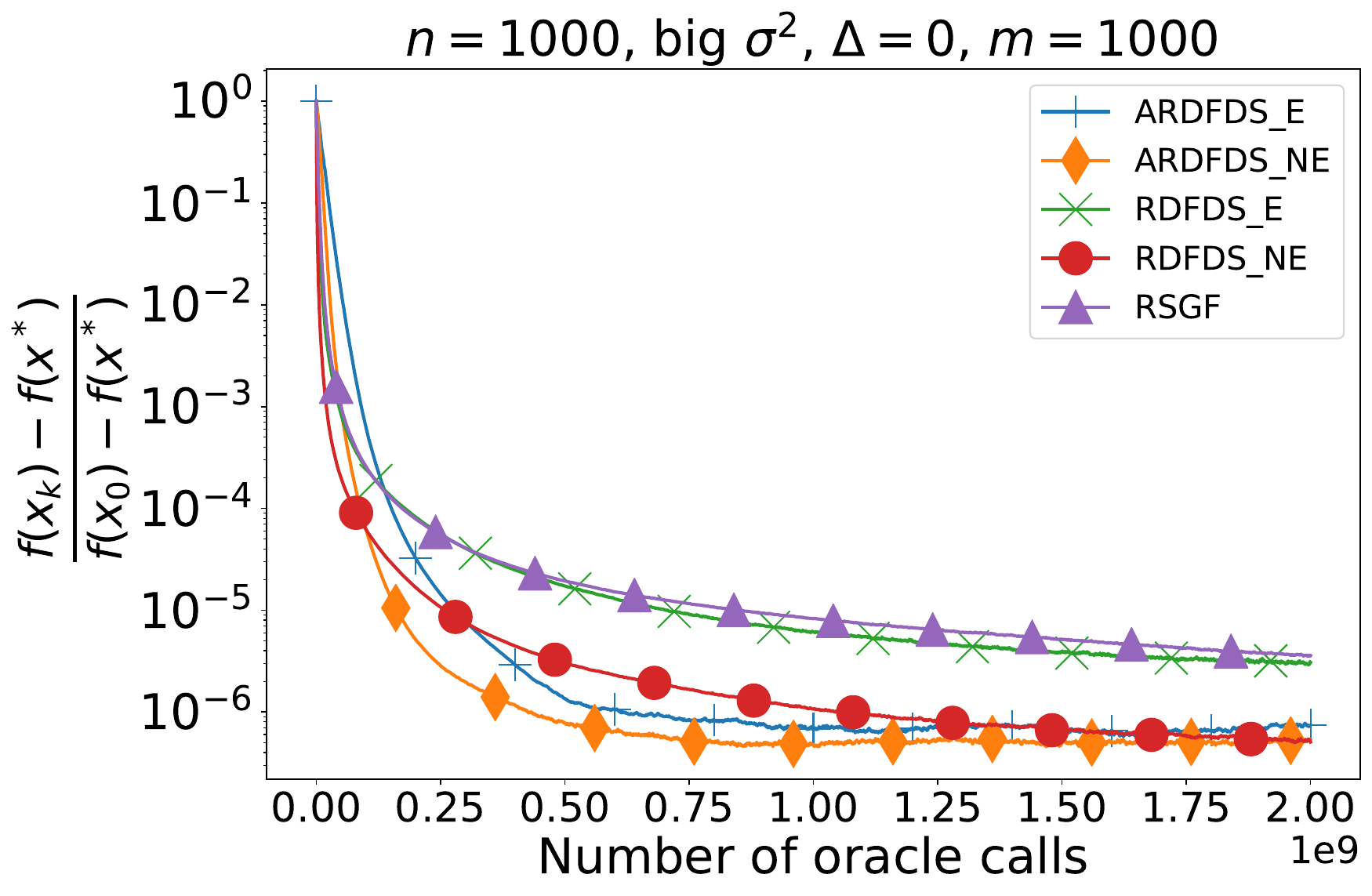}
\caption{\newstuff{Numerical results for minimizing Nesterov's function with noisy stochastic oracle having $\sigma^2 = \sigma_{\text{big}}^2$ for different sizes of mini-batch $m$ and dimensions of the problem $n$.}}
\label{fig:nesterov_big_sigma}
\end{figure}
We see in Figure~\ref{fig:nesterov_small_sigma} that for $\sigma^2 = \sigma_{\text{small}}^2$ it is sufficient to use mini-batches of the size $m=1$ to reach accuracy $\e = 10^{-3}$ and the overall picture is very similar to the one presented in Figure~\ref{fig:nesterov_sparsity_experiment}.}

\newstuff{In contrast, when $\sigma^2 = \sigma_{\text{big}}^2$ (Figure~\ref{fig:nesterov_big_sigma}) and $m=1$ the methods fail to reach the target accuracy. In these tests accelerated methods show higher sensitivity to the noise and, as a consequence, we see that for $n=500,1000$ and $m=10$ RDFDS{\_}NE reaches the accuracy $\e=10^{-3}$ faster than competitors justifying the following insight that we have from our theory: \pdd{when the variance is large, non-accelerated methods require smaller mini-batch size $m$ and are able to find} $\e$-solution faster than their accelerated counterparts.}

\subsubsection{Experiments with different noise level of the oracle}
\newstuff{Here we present the numerical experiments with different values of $\Delta$. To isolate the effect of the non-stochastic noise, we set $\sigma = 0$ for all tests reported in this subsection. We run the methods for problems with $n=100,500, 1000$ and chose the starting point in the same way as in Subsection~\ref{sec:exp_nesterov_sigma}. For each choice of the dimension $n$ we used three values of $\Delta$: $\Delta_{\text{small}} = \min\left\{\tfrac{\varepsilon^{3/2}\sqrt{2}}{\sqrt{ L_2\|x_0-x^*\|_1^2n\ln n}},\, \tfrac{2\varepsilon^2}{nL_2\|x_0-x^*\|_1^2}\right\}$, $\Delta_{\text{medium}} = 10^3\cdot\Delta_{\text{small}}$ and $\Delta_{\text{large}} = 10^6\cdot\Delta_{\text{small}}$ with $\e = 10^{-3}$. As one can see from Table~\ref{tbl:ARDFDS_params} when $\Delta = \Delta_{\text{small}}$ ARDFDS with $p=1$ is guaranteed to find an $\e$-solution. The results are reported in Figure~\ref{fig:nesterov_noise_exp}.
\begin{figure}
\centering 
\includegraphics[scale=0.14]{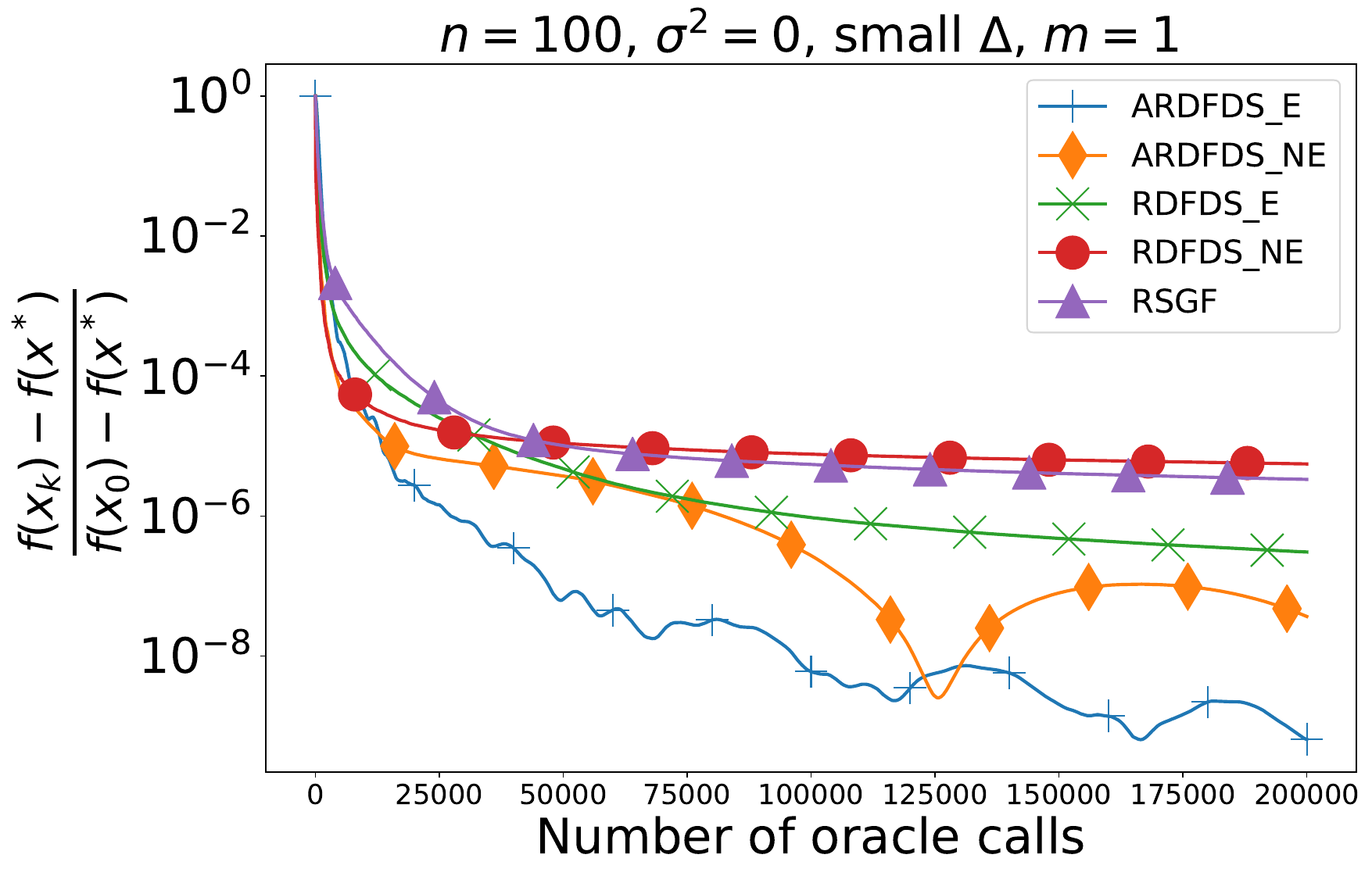}
\includegraphics[scale=0.14]{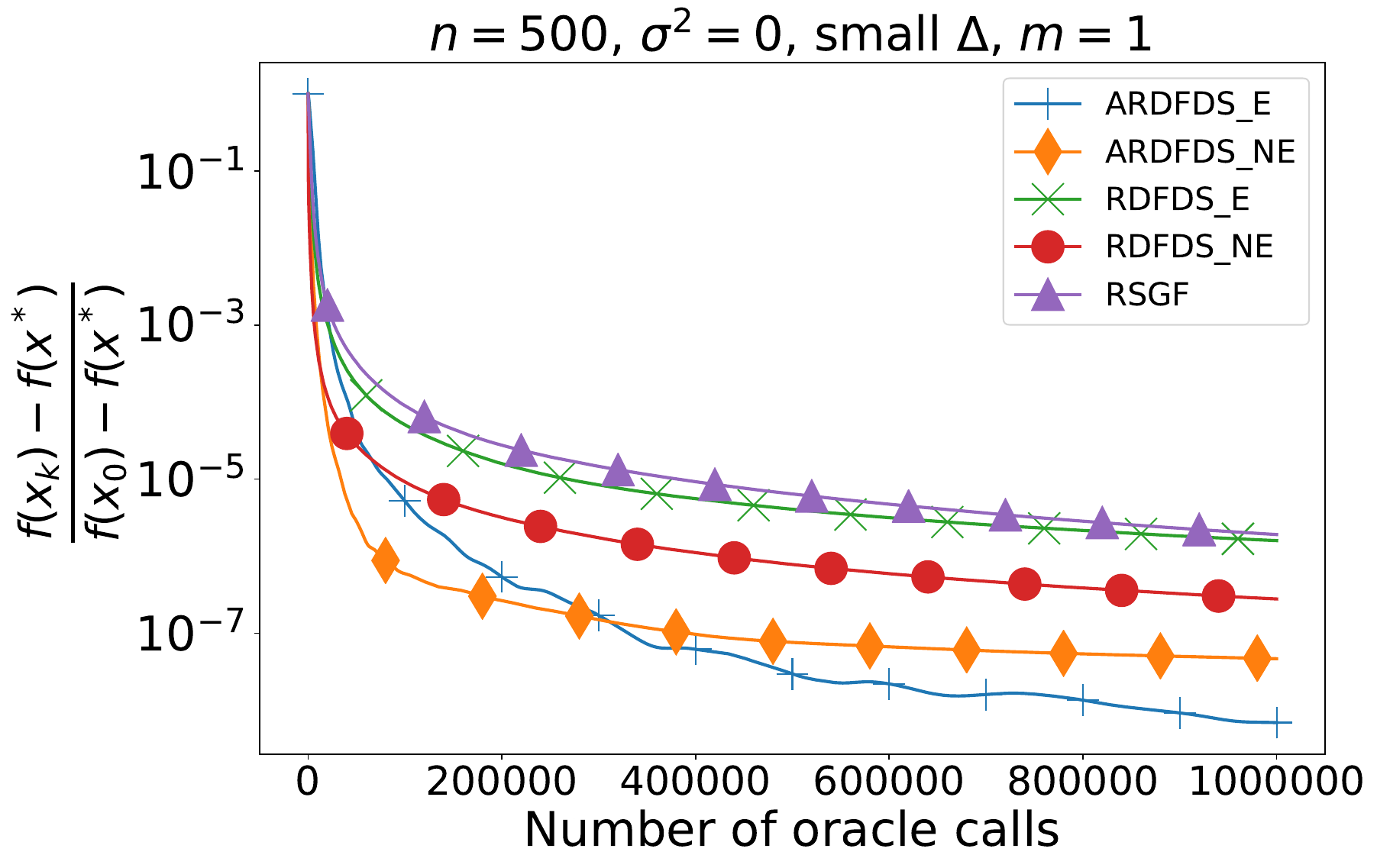}
\includegraphics[scale=0.14]{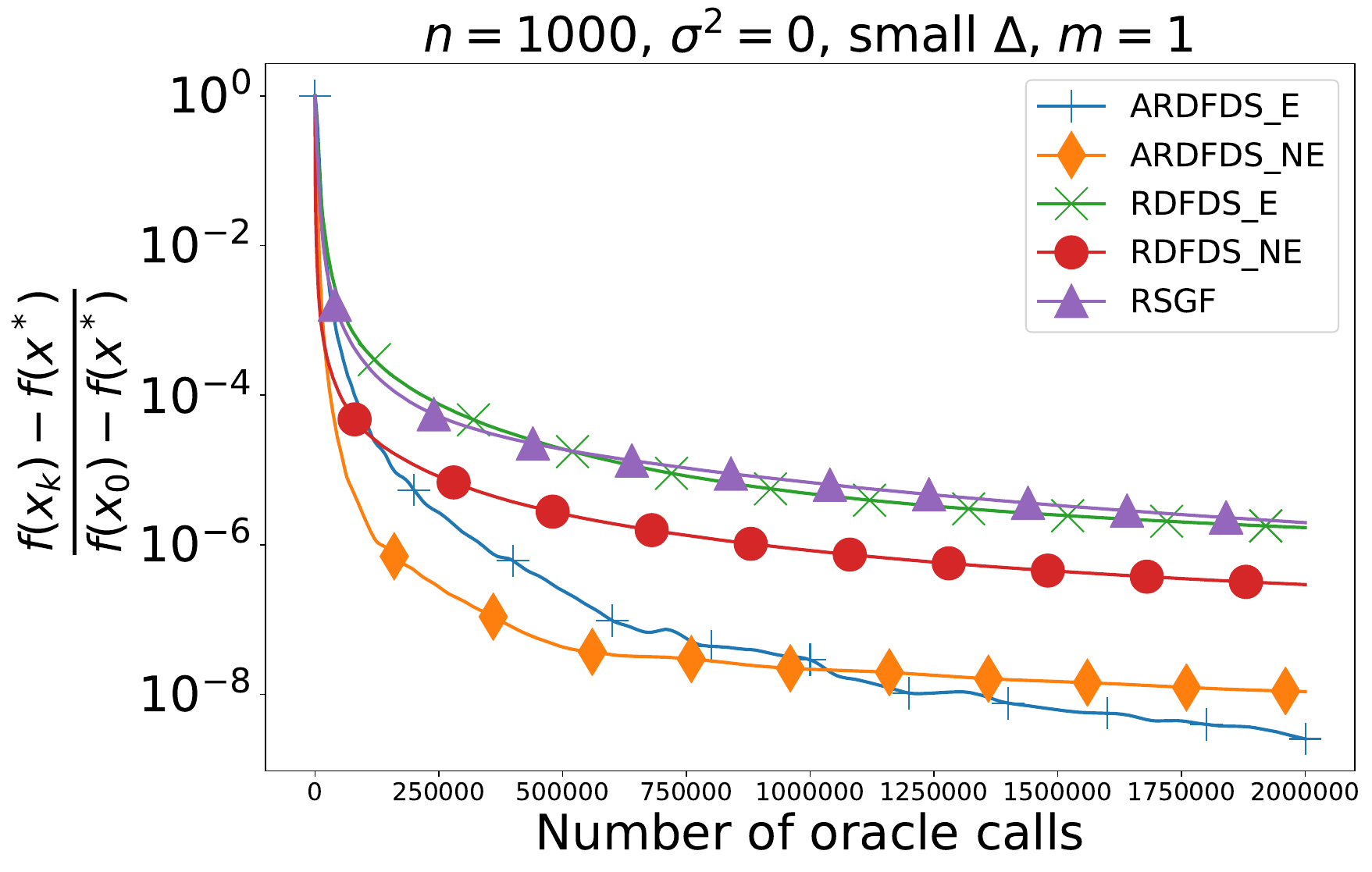}
\includegraphics[scale=0.14]{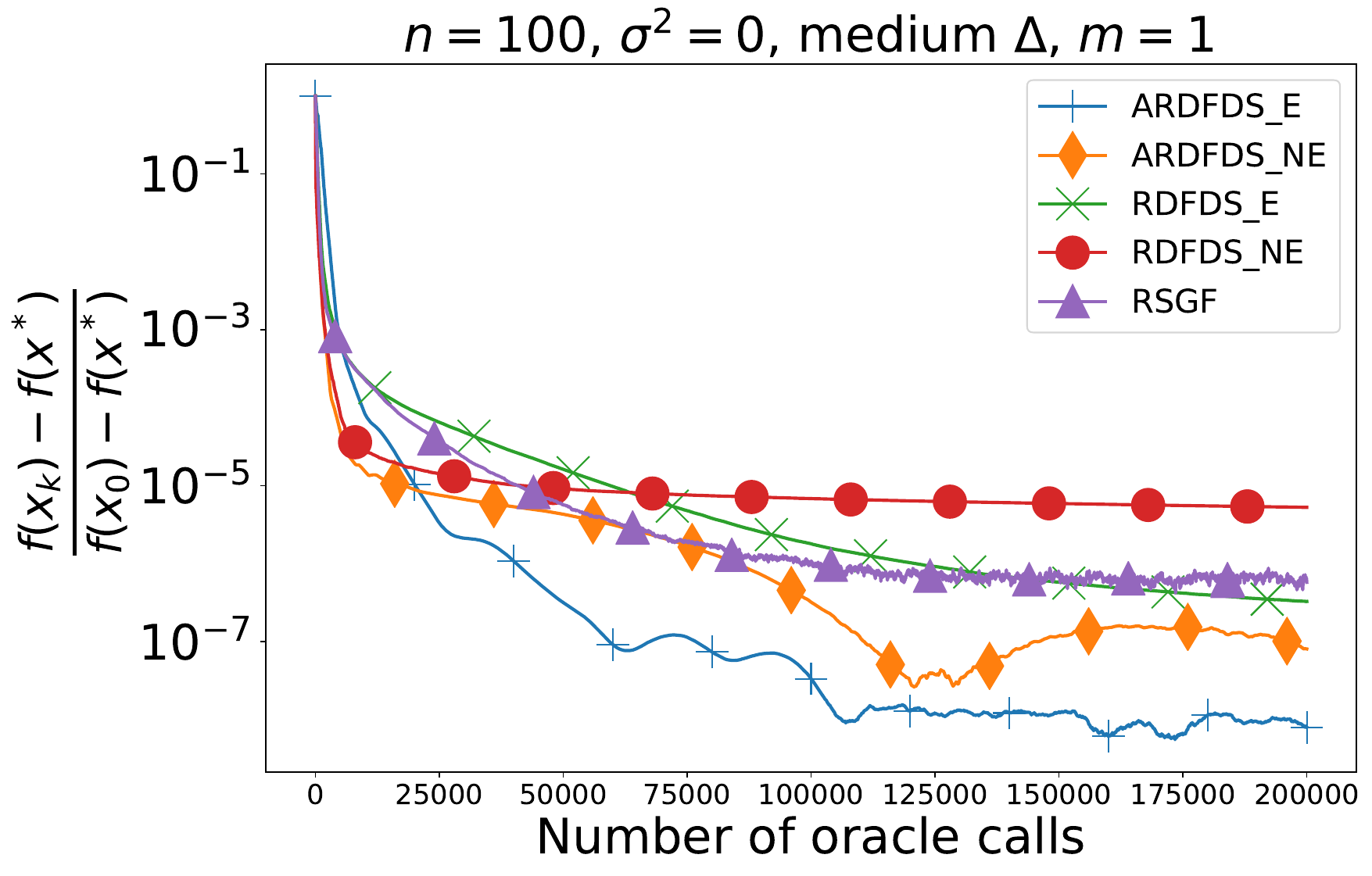}
\includegraphics[scale=0.14]{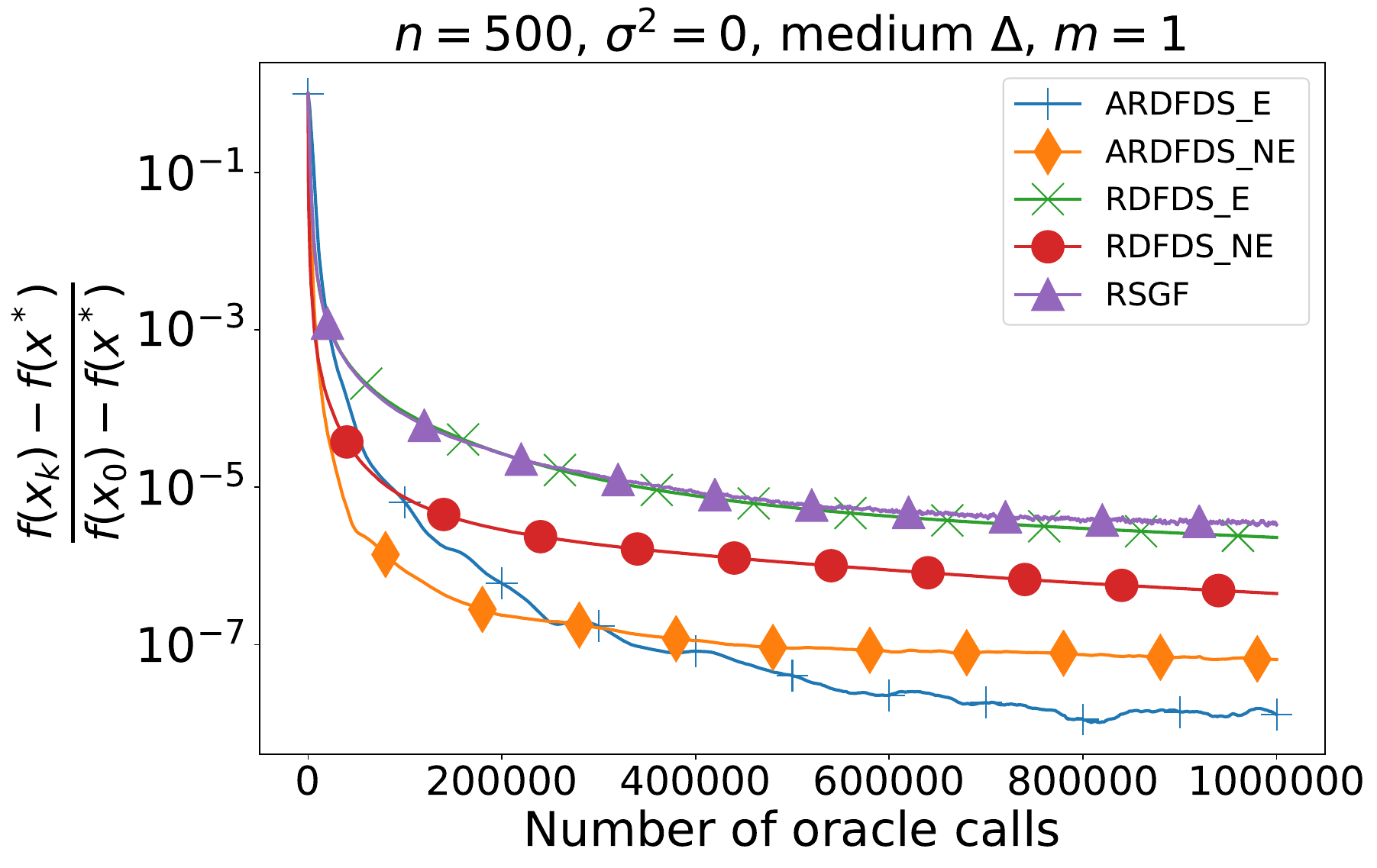}
\includegraphics[scale=0.14]{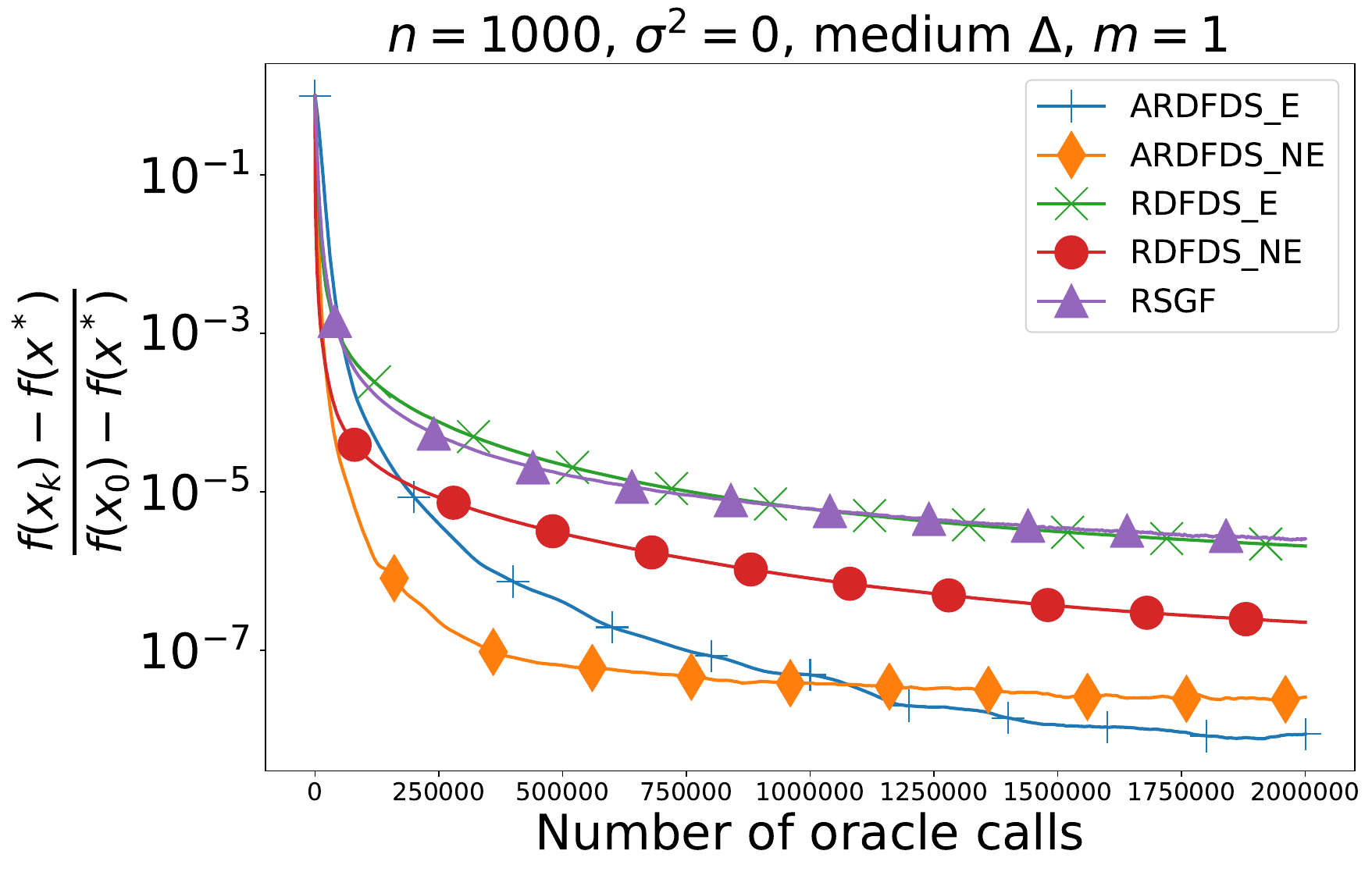}
\includegraphics[scale=0.14]{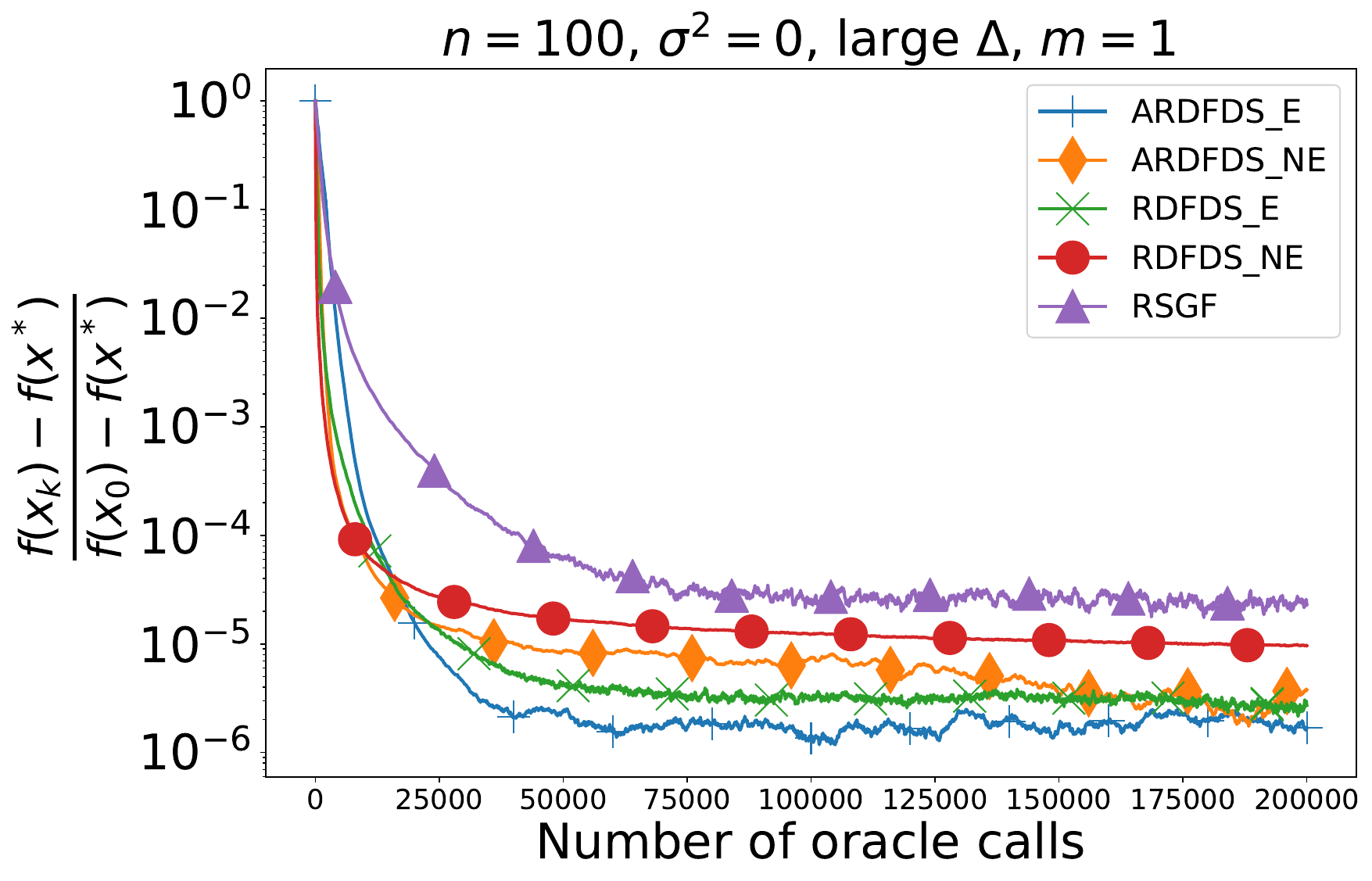}
\includegraphics[scale=0.14]{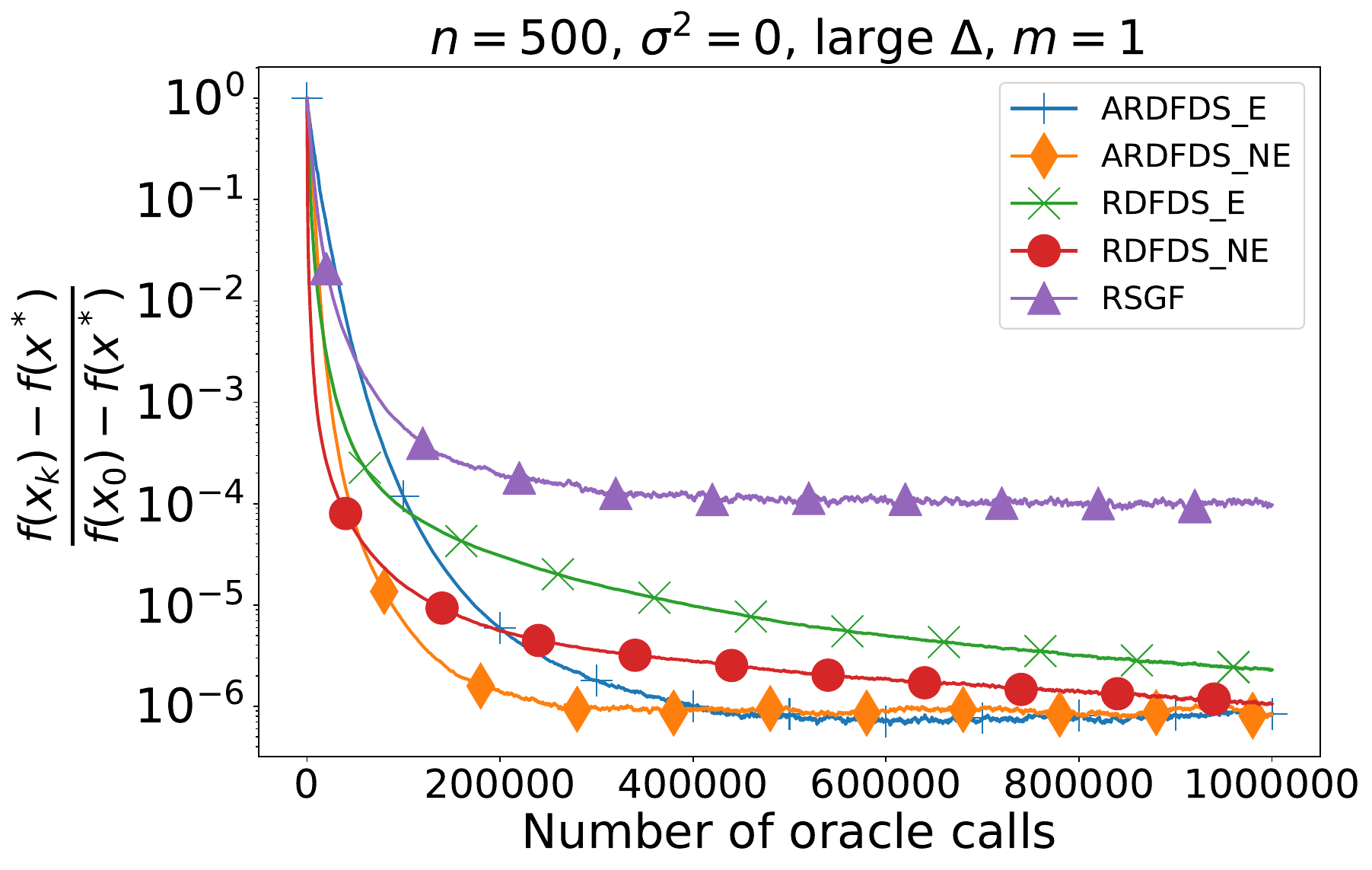}
\includegraphics[scale=0.14]{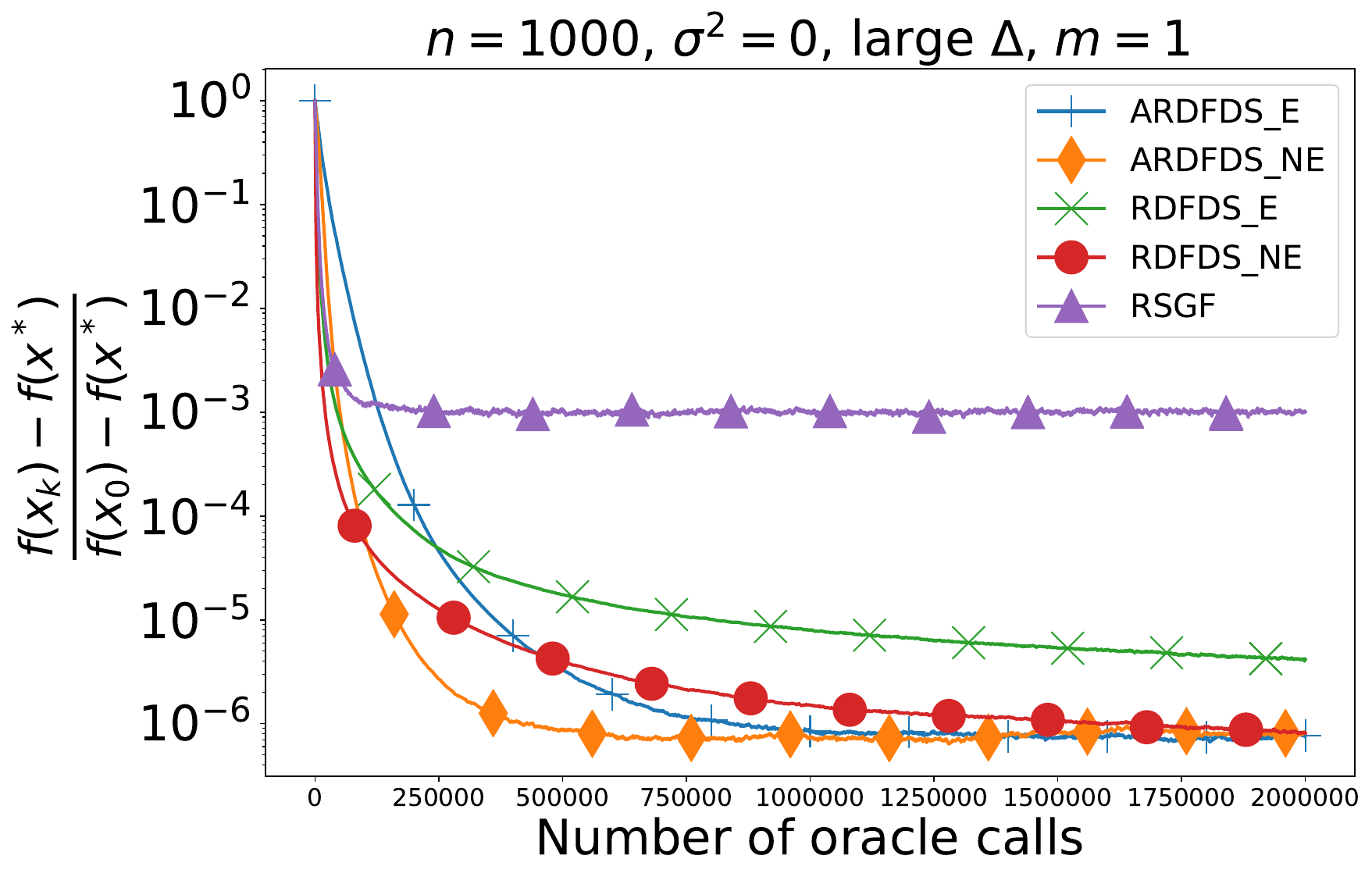}
\caption{\newstuff{Numerical results for minimizing Nesterov's function with noisy stochastic oracle having $\sigma^2 = 0$ for different $\Delta$ and dimensions of the problem $n$.}}
\label{fig:nesterov_noise_exp}
\end{figure}
We see that for larger values of $\Delta$ accelerated methods achieve worse accuracy than for small values of $\Delta$. However, in all experiments our methods succeeded to reach $\e$-solution with $\e=10^{-3}$ meaning that the noise level $\Delta$ in practice can be much larger than it is prescribed by our theory.}

\subsubsection{Experiment with large dimension}
\newstuff{In Figure~\ref{fig:nesterov_big_dimension} we report the experiment results for $n = 5000$, $\sigma^2 = \sigma_{\text{small}}^2 = \tfrac{\e^{\nicefrac{3}{2}}\sqrt{nL_2}}{\|x_0-x^*\|_1}$, \\
\noindent $\Delta = \Delta_{\text{small}} = \min\left\{\tfrac{\varepsilon^{3/2}\sqrt{2}}{\sqrt{ L_2\|x_0-x^*\|_1^2n\ln n}},\, \tfrac{2\varepsilon^2}{nL_2\|x_0-x^*\|_1^2}\right\}$, and $\e = 10^{-3}$.
\begin{figure}
\centering 
\includegraphics[scale=0.3]{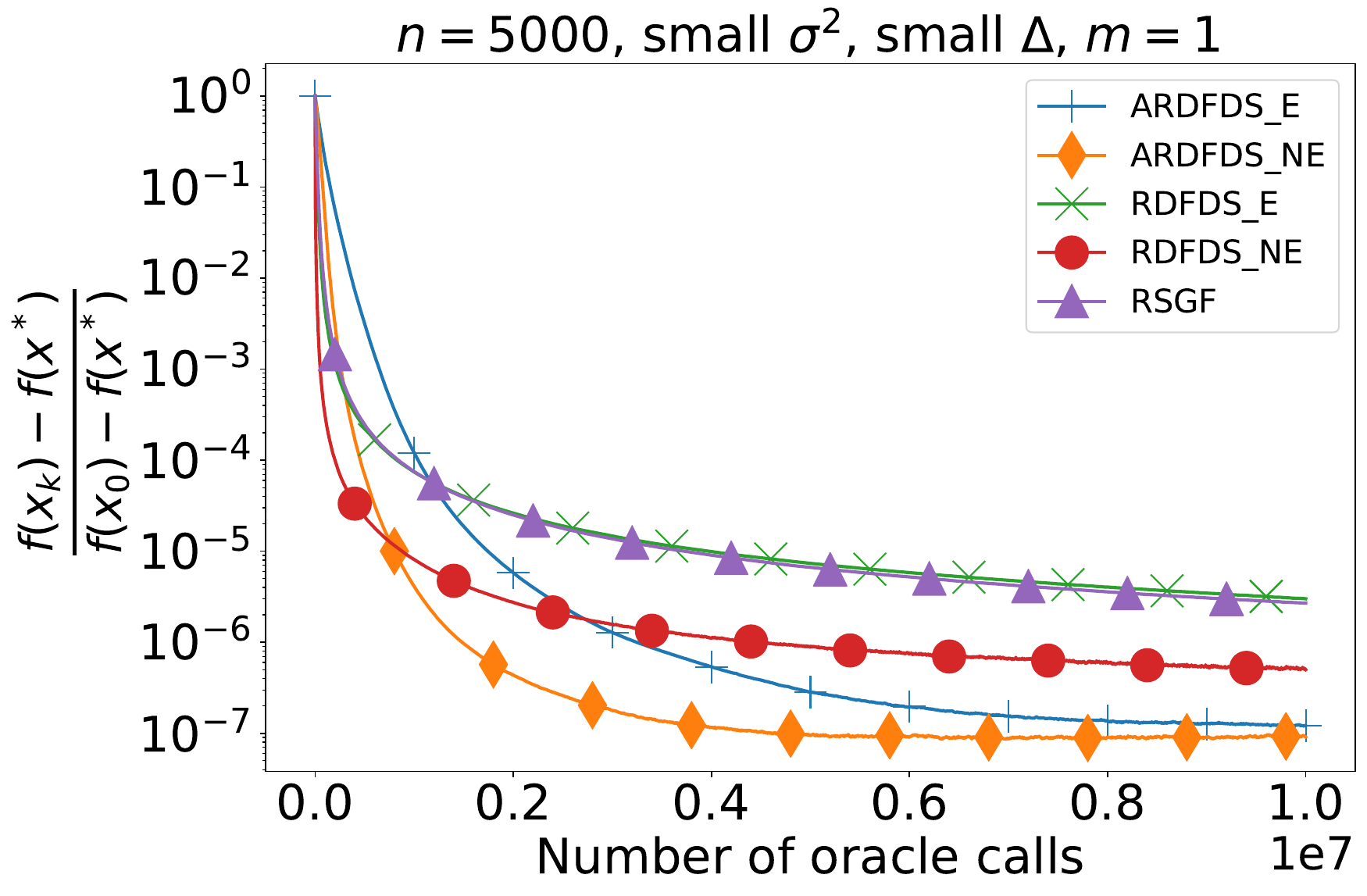}
\caption{\newstuff{Numerical results for minimizing Nesterov's function with noisy stochastic oracle having $\sigma^2 = \sigma_{\text{small}}^2$ and $\Delta = \Delta_{\text{small}}$ for the dimension of the problem $n = 5000$.}}
\label{fig:nesterov_big_dimension}
\end{figure}
The obtained results are in a good agreement with our theory and experiment results for smaller dimensions.}

\subsection{Experiments with logistic regression}
\newstuff{In this subsection we report the numerical results for our methods applied to the logistic regression problem:
\begin{equation}
    \min\limits_{x\in\R^n}\left\{f(x) = \frac{1}{M}\sum\limits_{i=1}^Mf_i(x)\right\},\quad f_i(x) = \log\left(1 + \exp\left(-y_i\cdot (Ax)_i\right)\right). \label{eq:logreg_loss}
\end{equation}
Here $f_i(x)$ is the loss on the $i$-th data point, $A\in\R^{M\times n}$ is a matrix of instances, $y \in\{-1,1\}^M$ is a vector of labels and $x\in\R^n$ is a vector of parameters (or weights). It can be easily shown that $f(x)$ is convex and $L_2$-smooth w.r.t. the Euclidean norm with $L_2 = \nicefrac{\sqrt{\lambda_{\max}(A^\top A)}}{4M}$ where $\lambda_{\max}(A^\top A)$ denotes the maximal eigenvalue of $A^\top A$. Moreover, problem \eqref{eq:logreg_loss} is a special case of \eqref{eq:PrSt} with $\xi$ being a random variable with the uniform distribution on $\{1,\ldots,M\}$.}

\newstuff{For our experiments we use the data from LIBSVM library \cite{chang2011libsvm}, see also Table~\ref{tab:logreg_datasets} summarizing the information about the datasets we used.
\begin{table}[h]
    \centering
    \begin{tabular}{|c|c|c|c|c|c|}
        \hline
         & {\tt heart} & {\tt diabetes} & {\tt a9a} & {\tt phishing} & {\tt w8a} \\
        \hline
        Size \pdd{$M$} & $270$ & $768$ & $32561$ & $11055$ & $49749$ \\
        \hline
        Dimension \pdd{$n$} & $13$ & $8$ & $123$ & $68$ & $300$ \\
        \hline
    \end{tabular}
    \caption{Summary of used datasets.}
    \label{tab:logreg_datasets}
\end{table}
In all test we chose $t=10^{-8}$ and the starting point $x_0$ such that it differs from $x^*$ only in the first component and $f(x_0) - f(x^*) \sim 10$. We use standard solvers from scipy library to obtain a very good approximation of a solution $x^*$ and use it to measure the quality of the approximations by other algorithms. The results for the batch (and hence deterministic) methods with $m=M$ and mini-batch stochastic methods are presented in Figures~\ref{fig:logreg_full_batch}~and~\ref{fig:logreg_mini_batch} respectively.
\begin{figure}
\centering 
\includegraphics[scale=0.14]{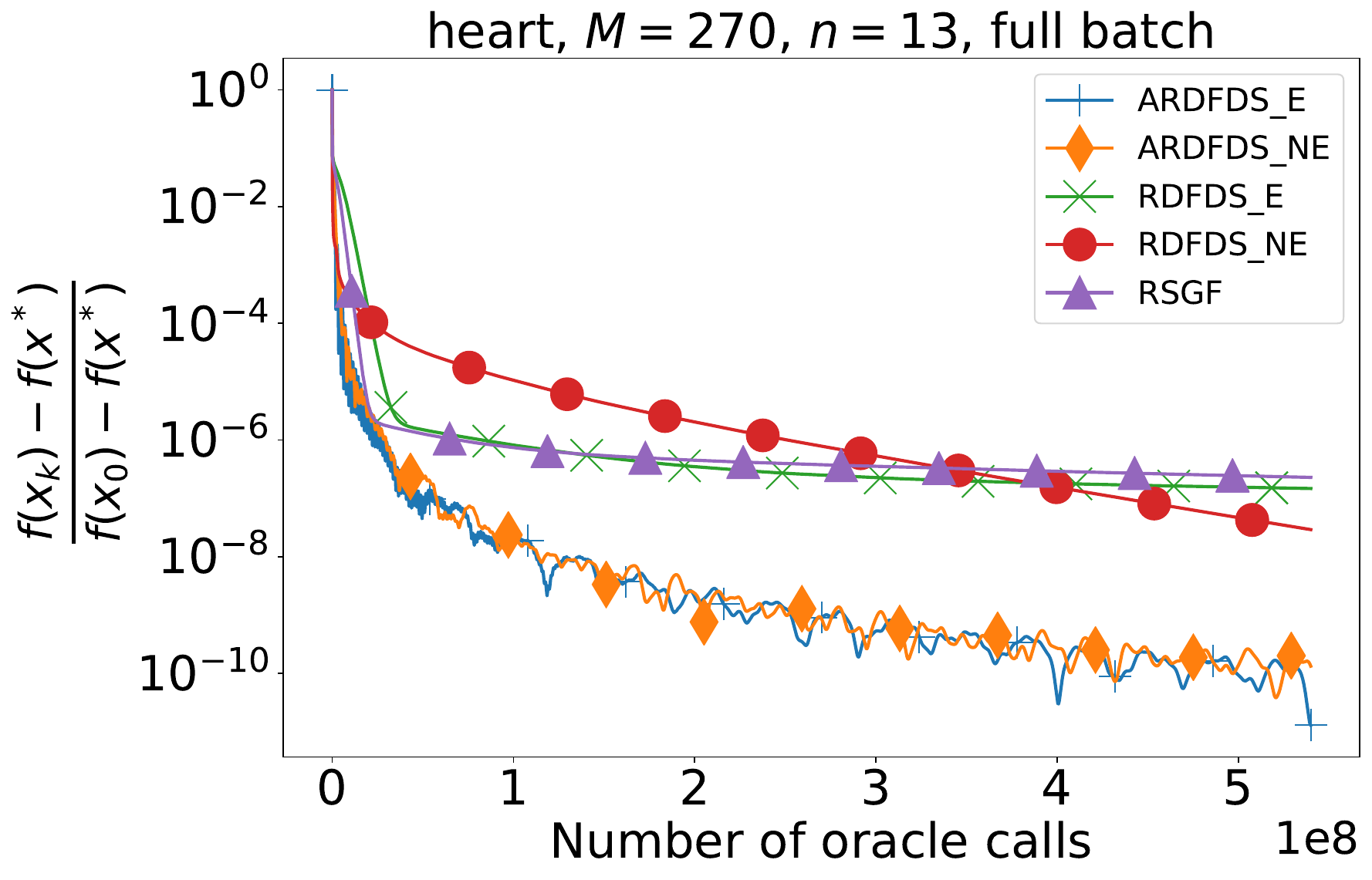}
\includegraphics[scale=0.14]{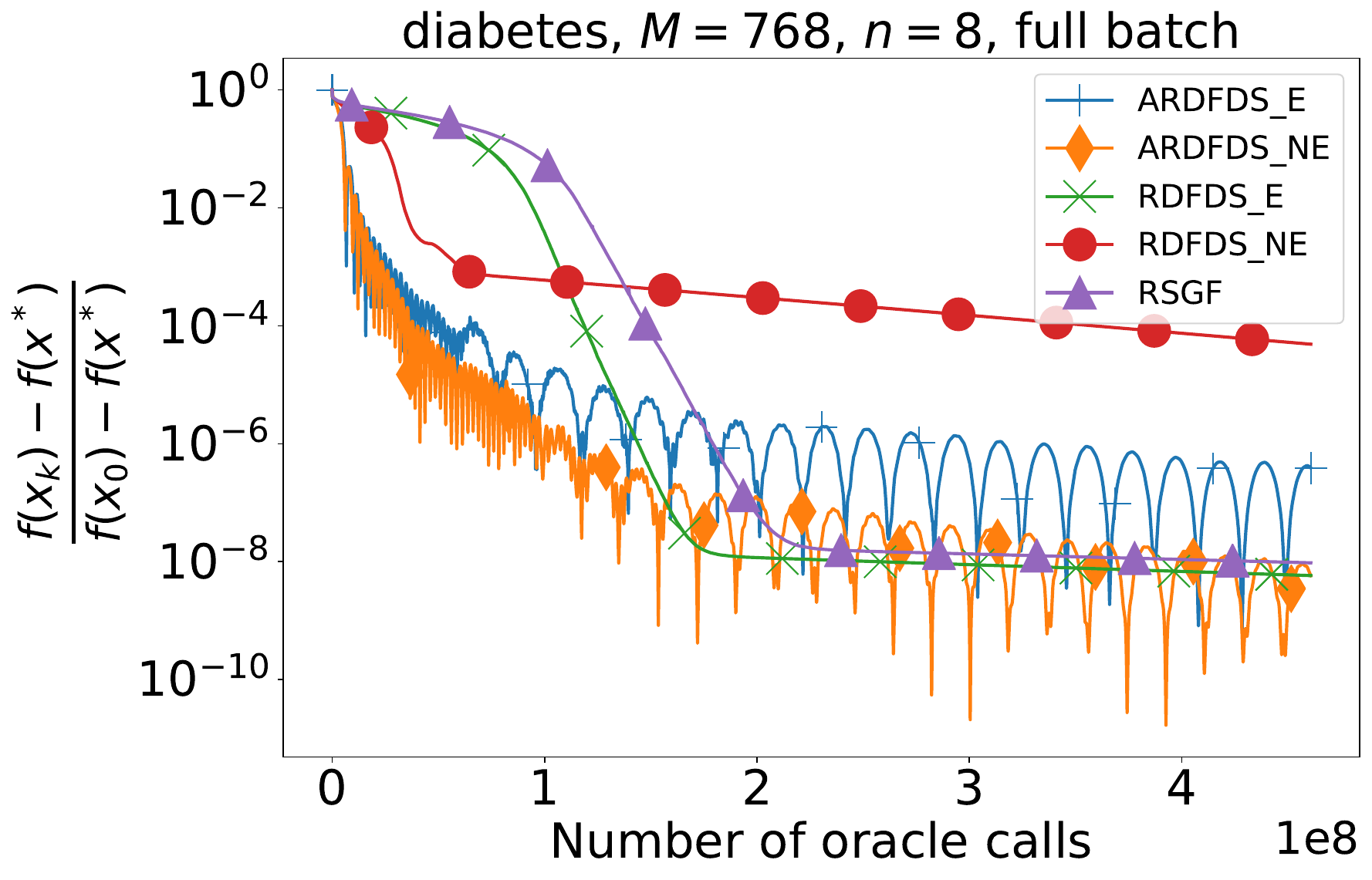}
\includegraphics[scale=0.14]{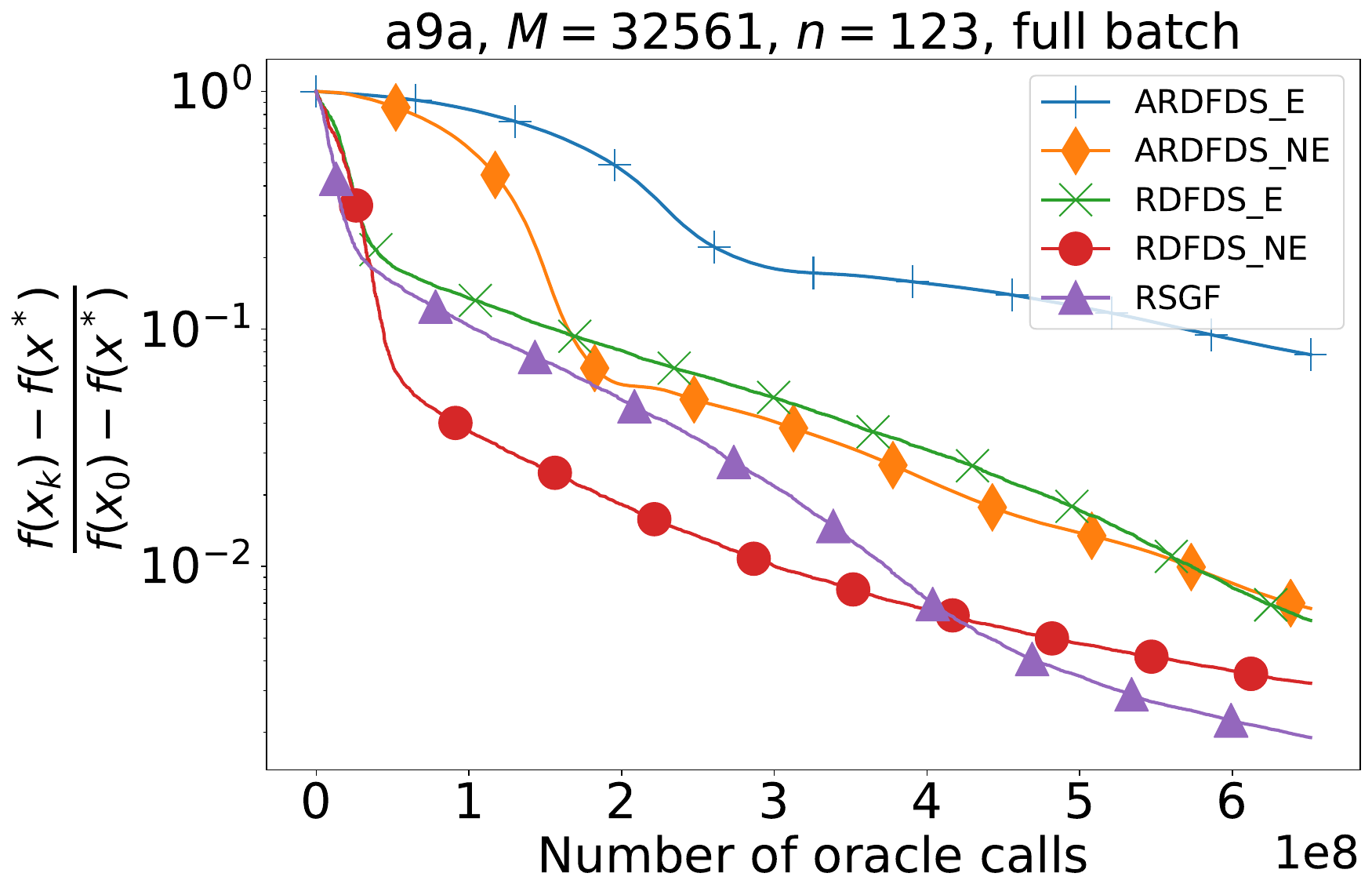}
\includegraphics[scale=0.14]{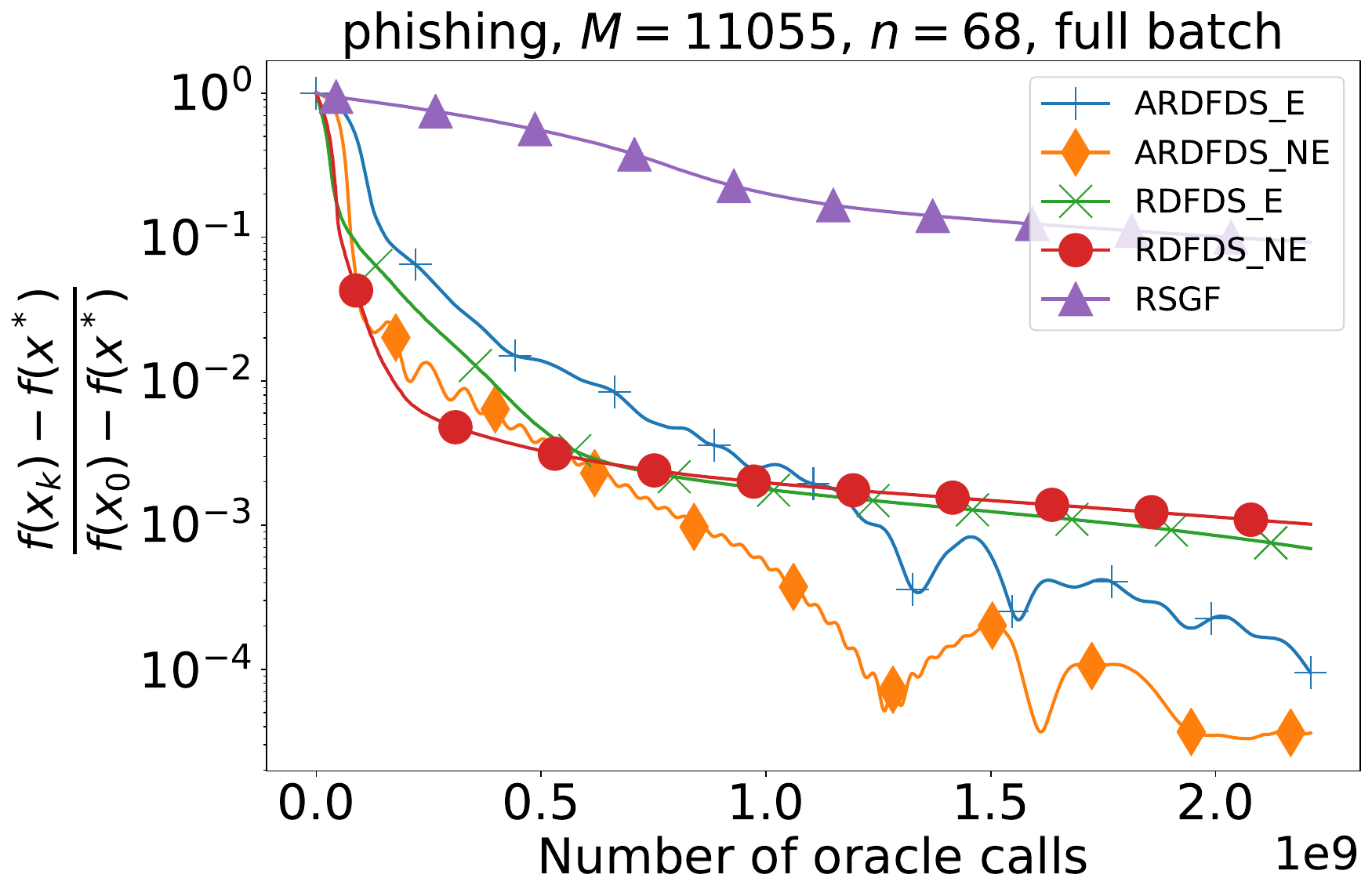}
\includegraphics[scale=0.14]{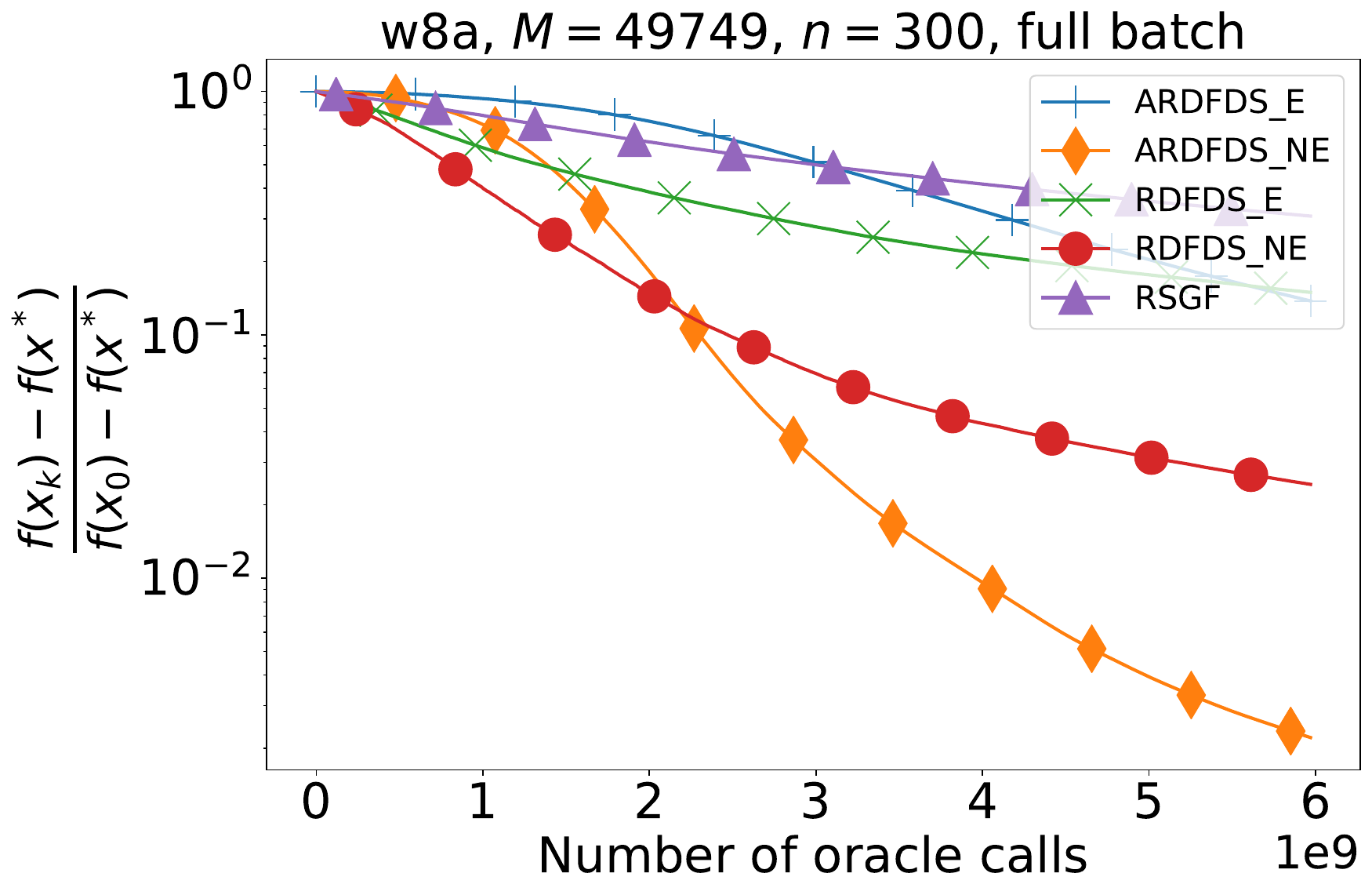}
\caption{\newstuff{Numerical results for solving logistic regression problem \eqref{eq:logreg_loss} for different datasets using batch methods with $m=M$.}}
\label{fig:logreg_full_batch}
\end{figure}
\begin{figure}
\centering 
\includegraphics[scale=0.14]{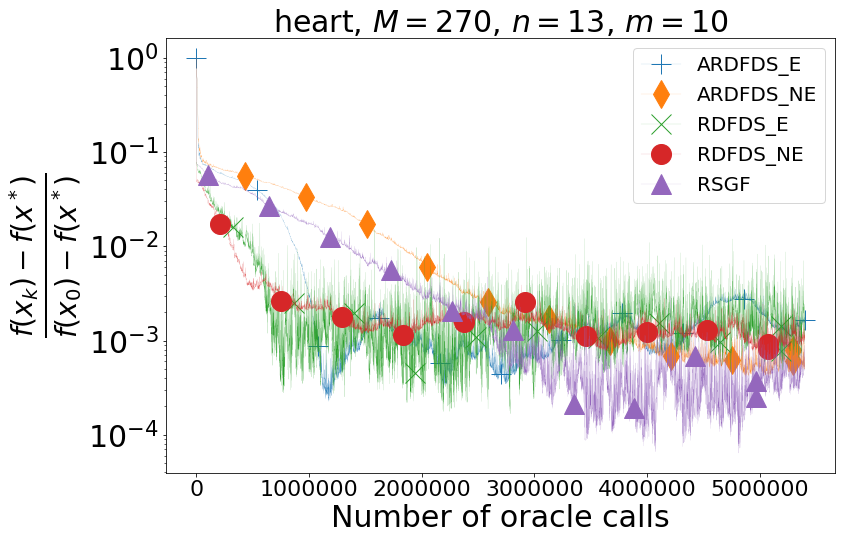}
\includegraphics[scale=0.14]{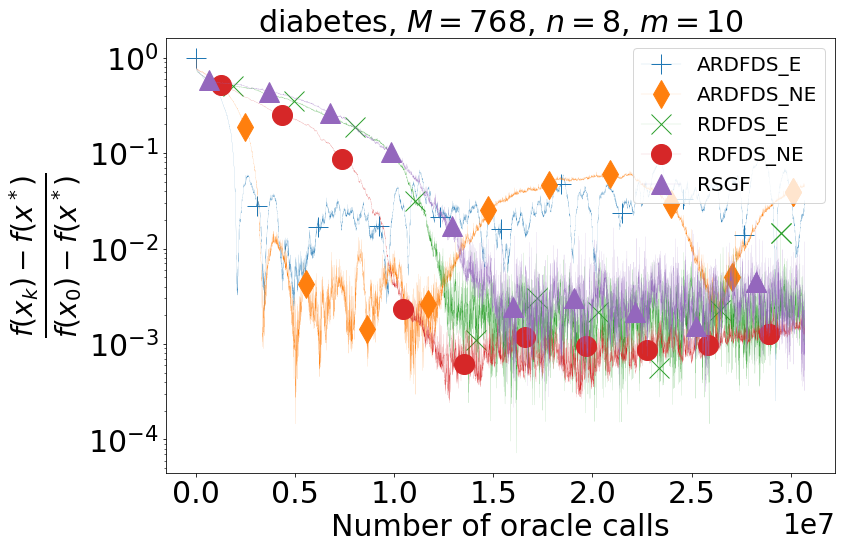}
\includegraphics[scale=0.14]{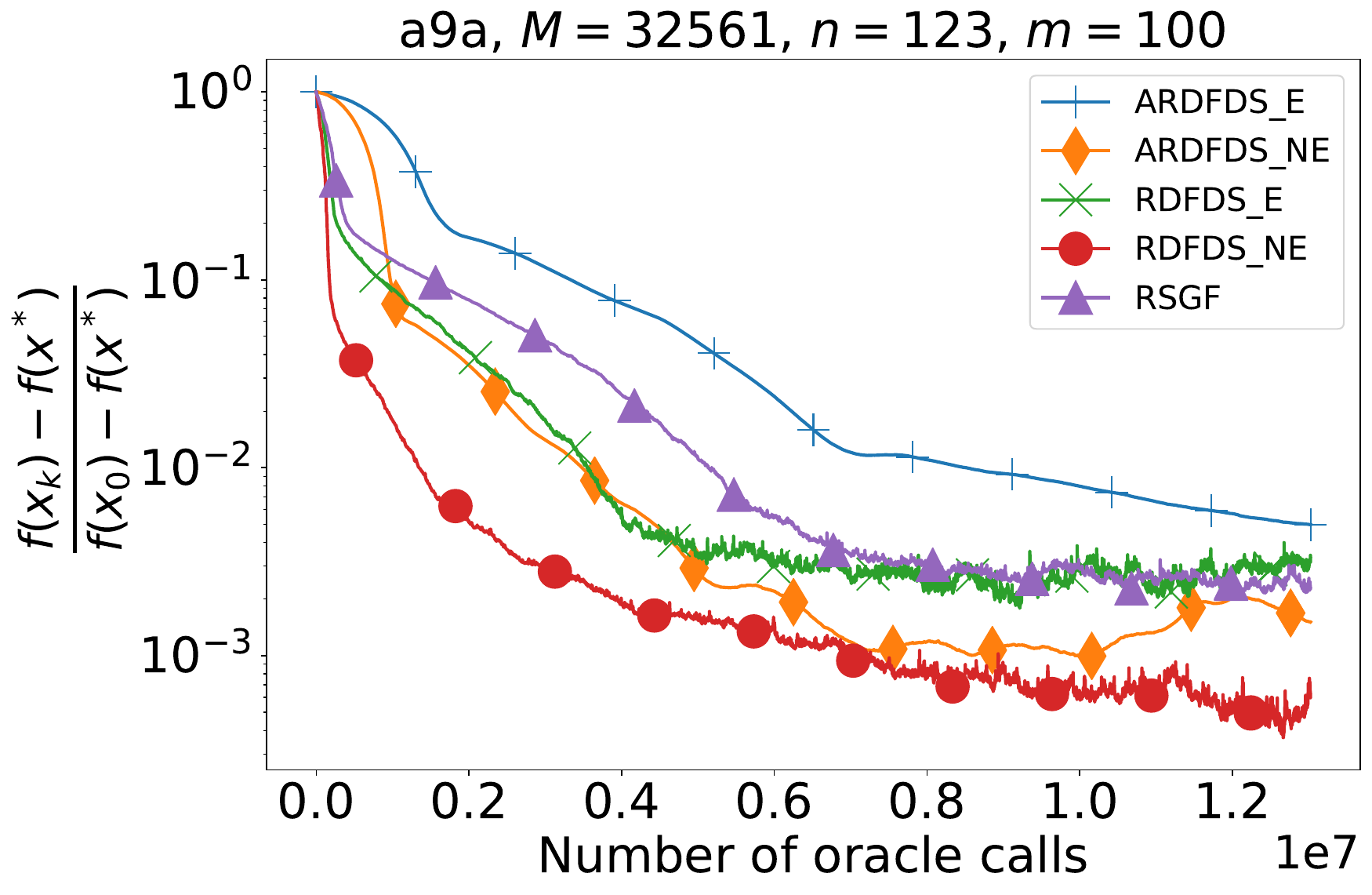}
\includegraphics[scale=0.14]{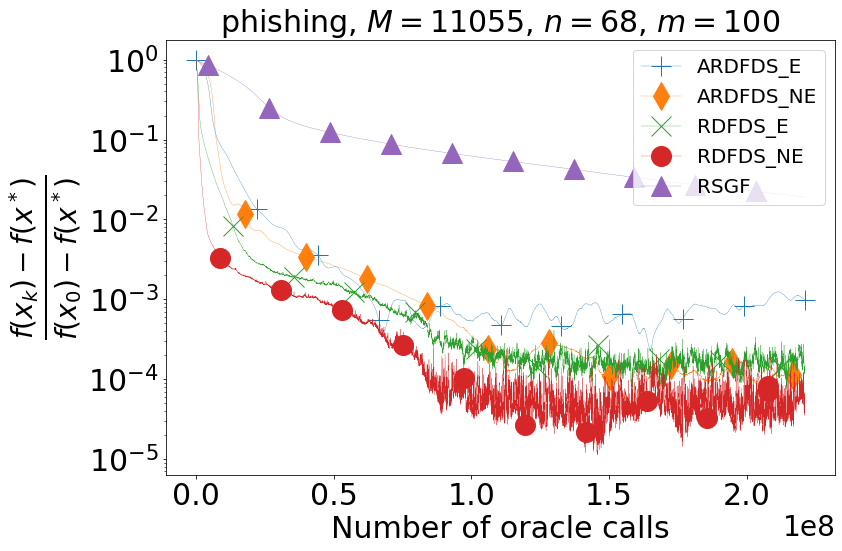}
\includegraphics[scale=0.14]{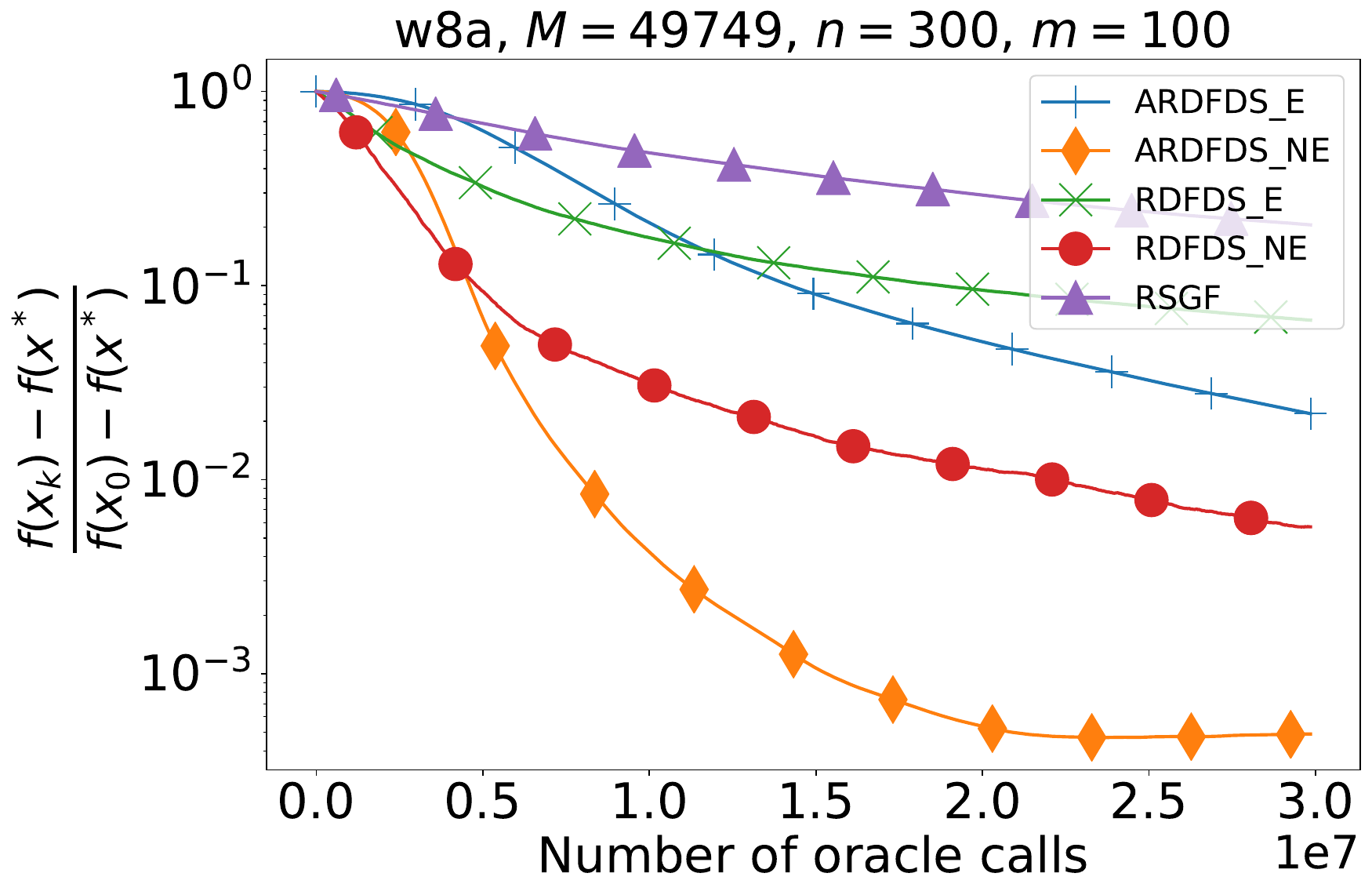}
\caption{\newstuff{Numerical results for solving logistic regression problem \eqref{eq:logreg_loss} for different datasets using mini-batch stochastic methods.}}
\label{fig:logreg_mini_batch}
\end{figure}
In all cases methods with the $1$-norm proximal setup show the best or comparable with the best results.}

\section{Conclusion}
In this paper, we propose two new algorithms for stochastic smooth derivative-free \pd{convex} optimization with two-point feedback and \pdd{inexact function values oracle}. Our first algorithm is an accelerated one and the second one is a non-accelerated one. \pdd{Notably}, despite the traditional choice of \pd{$2$}-norm proximal setup for unconstrained optimization problems, 
\pd{our analysis has yielded better complexity bounds for the method with \pd{$1$}-norm proximal setup than the ones with \pd{$2$}-norm proximal setup. \pdd{This is also confirmed by numerical experiments. }}

\bibliographystyle{siamplain}
\bibliography{references}

\appendix

{\color{red}
}

\section{Proof of Lemma \ref{Lm:MainTechLM}}
\pd{In this appendix} we prove that, for $e \in RS_2\left( 1 \right)$,  \pd{$q \geqslant 2$, and $n\geqslant 8$,}
\begin{align}
        \EE[\|e\|_q^2] & \leqslant \min\{q-1,\,16\ln n - 8\}n^{\tfrac{2}{q}-1}, \label{lemm1:expect_q_norm} \\
        \EE[\langle s,\, e\rangle^2\|e\|_q^2] & \leqslant 6\|s\|_2^2\min\{q-1,16\ln n -8\}n^{\tfrac{2}{q}-2}. \label{lemm1:expect_inner_product}
\end{align}
\pd{Throughout this appendix, to simplify the notation, we denote by $\EE$ the expectation w.r.t. random vector $e \in RS_2\left( 1 \right)$.}

We start \pdd{by} proving the following inequality, which could be \pd{not tight} for \pd{large} $q$:
\begin{equation}\label{lemm1:rough_estimation_of_expectation_q_norm}
    \EE[\|e\|_q^2] \leqslant (q-1)n^{\tfrac{2}{q}-1},\quad 2\leqslant q < \infty.
\end{equation}
We have
\begin{equation}\label{lemm1:jensen}
    \begin{array}{c}
        \EE[\|e\|_q^2] = \EE\left[\left(\sum\limits_{k=1}^{n}|e_k|^q\right)^{\tfrac{2}{q}} \right] \overset{\circledOne}{\leqslant} \left(\EE\left[\sum\limits_{k=1}^{n}|e_k|^q\right]\right)^{\tfrac{2}{q}} \overset{\circledTwo}{=} \left(n\EE[|e_2|^q]\right)^{\tfrac{2}{q}},
    \end{array}
\end{equation}
where $\circledOne$ is due to probabilistic version of Jensen's inequality (function $\varphi(x) = x^{\tfrac{2}{q}}$ is concave, because $q\geqslant2$) and $\circledTwo$ is because expectation is linear and components of the vector $e$ are identically distributed. \pd{We also denote by $e_k$ the $k$-th component of $e$. In particular, $e_2$ is the second component.}

\pdd{By the} Poincare lemma, $e$ \pdd{has the same distribution as} $\tfrac{\xi}{\sqrt{\xi_1^2 + \dots + \xi_n^2}}$, where $\xi$ \pdd{is the standard Gaussian random vector with zero mean and identity covariance matrix.} Then
\begin{equation*}
    \begin{array}{rl}
        \EE[|e_2|^q] &= \EE\left[\tfrac{|\xi_2|^q}{\left(\xi_1^2+\ldots+\xi_n^2\right)^{\tfrac{q}{2}}}\right]\\
        &= \idotsint\limits_{\R^n}|x_2|^q\left(\sum\limits_{k=1}^{n}x_k^2\right)^{-\tfrac{q}{2}}\cdot\tfrac{1}{(2\pi)^{\tfrac{n}{2}}}\cdot \exp\left(-\tfrac{1}{2}\sum\limits_{k=1}^{n}x_k^2\right)dx_1\ldots dx_n.
    \end{array}
\end{equation*}
\pdd{For the transition to the spherical coordinates}
\begin{equation*}
    \begin{array}{rl}
        x_1 &= r\cos\varphi\sin\theta_1\ldots\sin\theta_{n-2},\quad x_2 = r\sin\varphi\sin\theta_1\ldots\sin\theta_{n-2},\\
        x_3 &= r\cos\theta_1\sin\theta_2\ldots\sin\theta_{n-2},\quad x_4 = r\cos\theta_2\sin\theta_3\ldots\sin\theta_{n-2},\\
        &\ldots\\
        x_n &= r\cos\theta_{n-2},\quad r>0,\,\varphi \in [0,2\pi),\, \theta_i \in [0,\pi],\, \pd{i = 1,...,n-2} 
    \end{array}
\end{equation*}
the Jacobian  \pd{satisfies} 
\begin{equation*}
    \det\left(\tfrac{\partial(x_1,\ldots,x_n)}{\partial(r,\varphi,\theta_1,\theta_2,\ldots,\theta_{n-2})}\right) = r^{n-1}\sin\theta_1(\sin\theta_2)^2\ldots(\sin\theta_{n-2})^{n-2}.
\end{equation*}
\pdd{In the new coordinates we have } 
\begin{equation*}
    \begin{array}{rl}
        \EE[|e_2|^q] &= \idotsint\limits_{\substack{r>0,\,\varphi \in [0,2\pi),\\ \theta_i \in [0,\pi],\, i = \pd{1,...,n-2}}}r^{n-1}|\sin\varphi|^q|\sin\theta_1|^{q+1}|\sin\theta_2|^{q+2}\ldots|\sin\theta_{n-2}|^{q+n-2}\\
        &\cdot\tfrac{\pd{\exp{(-\tfrac{r^2}{2}})}}{(2\pi)^{\tfrac{n}{2}}}dr\ldots d\theta_{n-2}=\tfrac{1}{(2\pi)^{\tfrac{n}{2}}} I_r\cdot I_\varphi \cdot I_{\theta_1}\cdot I_{\theta_2}\cdot\ldots\cdot I_{\theta_{n-2}},
    \end{array}
\end{equation*}
where $I_r = \int\limits_{0}^{+\infty}r^{n-1}\pd{\exp\left({-\tfrac{r^2}{2}}\right)}dr$, \\ 
$I_\varphi = \int\limits_{0}^{2\pi}|\sin\varphi|^qd\varphi = 2\int\limits_{0}^{\pi}|\sin\varphi|^qd\varphi$, $I_{\theta_i} = \int\limits_{0}^{\pi}|\sin\theta_i|^{q+i}d\theta_i$ for $\pd{i=1,...,n-2}$. \pdd{Next we calculate these integrals starting with} $I_r$:
\begin{equation*}
    \begin{array}{c}
        I_r = \int\limits_{0}^{+\infty}r^{n-1}\pd{\exp\left({-\tfrac{r^2}{2}}\right)}dr \overset{\newstuff{r = \sqrt{2t}}}{=} \int\limits_{0}^{+\infty}(2t)^{\tfrac{n}{2}-1}\pd{\exp\left(-t\right)}dt = 2^{\tfrac{n}{2}-1}\Gamma(\tfrac{n}{2}).
    \end{array}
\end{equation*}
To compute  \pdd{the} other integrals, we consider the following integral for $\alpha > 0$:
\begin{equation*}
    \begin{array}{rl}
        \int\limits_{0}^{\pi} |\sin\varphi|^\alpha d\varphi &= 2\int\limits_{0}^{\tfrac{\pi}{2}}|\sin\varphi|^\alpha d\varphi = 2\int\limits_{0}^{\tfrac{\pi}{2}}(\sin^2\varphi)^{\tfrac{\alpha}{2}}d\varphi \\ 
        &\overset{\newstuff{t = \sin^2\varphi}}{=} \int\limits_{0}^{1}t^{\tfrac{\alpha-1}{2}}(1-t)^{-\tfrac{1}{2}}dt = B(\tfrac{\alpha+1}{2},\,\tfrac{1}{2}) = \tfrac{\Gamma(\tfrac{\alpha+1}{2})\Gamma(\tfrac{1}{2})}{\Gamma(\tfrac{\alpha+2}{2})} = \sqrt{\pi} \tfrac{\Gamma(\tfrac{\alpha+1}{2})}{\Gamma(\tfrac{\alpha+2}{2})}.
    \end{array}
\end{equation*}
\pdd{This gives}
\begin{equation}\label{lemm1:expectation_component}
    \begin{array}{rl}
        \EE[|e_2|^q] = \tfrac{1}{(2\pi)^{\tfrac{n}{2}}} I_r\cdot I_\varphi \cdot I_{\theta_1}\cdot I_{\theta_2}\cdot\ldots\cdot I_{\theta_{n-2}}\\
        \hspace{-5em} = \tfrac{1}{(2\pi)^{\tfrac{n}{2}}}\cdot 2^{\tfrac{n}{2}-1}\Gamma(\tfrac{n}{2})\cdot 2\sqrt{\pi}\tfrac{\Gamma(\tfrac{q+1}{2})}{\Gamma(\tfrac{q+2}{2})}\cdot\sqrt{\pi}\tfrac{\Gamma(\tfrac{q+2}{2})}{\Gamma(\tfrac{q+3}{2})}\cdot\ldots\cdot\sqrt{\pi}\tfrac{\Gamma(\tfrac{q+n-1}{2})}{\Gamma(\tfrac{q+n}{2})}= \tfrac{1}{\sqrt{\pi}}\cdot\tfrac{\Gamma(\tfrac{n}{2})\Gamma(\tfrac{q+1}{2})}{\Gamma(\tfrac{q+n}{2})}.
    \end{array}
\end{equation}
\pdd{The next step is} to show that, for all $q\geqslant 2$, 
\begin{equation}\label{lemm1:key_estimation}
    \tfrac{1}{\sqrt{\pi}}\cdot\tfrac{\Gamma(\tfrac{n}{2})\Gamma(\tfrac{q+1}{2})}{\Gamma(\tfrac{q+n}{2})} \leqslant \left(\tfrac{q-1}{n}\right)^{\tfrac{q}{2}}.
\end{equation}
\pdd{First we} show that \eqref{lemm1:key_estimation} holds for $q=2$ and arbitrary $n$:
\begin{equation*}
    \tfrac{1}{\sqrt{\pi}}\cdot\tfrac{\Gamma(\tfrac{n}{2})\Gamma(\tfrac{2+1}{2})}{\Gamma(\tfrac{2+n}{2})} - \tfrac{1}{n} = \tfrac{1}{\sqrt{\pi}} \cdot \tfrac{\Gamma(\tfrac{n}{2})\cdot\tfrac{1}{2}\Gamma(\tfrac{1}{2})}{\tfrac{n}{2}\Gamma(\tfrac{n}{2})} - \tfrac{1}{n} = \tfrac{1}{n}-\tfrac{1}{n} = 0 \leqslant 0.
\end{equation*}
\pdd{Next, we c}onsider the function $f_n(q) = \tfrac{1}{\sqrt{\pi}}\cdot\tfrac{\Gamma(\tfrac{n}{2})\Gamma(\tfrac{q+1}{2})}{\Gamma(\tfrac{q+n}{2})} - \left(\tfrac{q-1}{n}\right)^{\tfrac{q}{2}}$, where $q \geqslant 2$, \pdd{and digamma function} $\psi(x) = \tfrac{d(\ln(\Gamma(x)))}{dx}$ with \pd{scalar argument $x > 0$. For the gamma} function it holds \pdd{that} $\Gamma(x+1) = x\Gamma(x),\, x>0.$ Taking natural logarithm \pdd{in both sides} and derivative w.r.t. $x$, we get $\tfrac{d(\ln(\Gamma(x+1)))}{dx} = \tfrac{d(\ln(\Gamma(x)))}{dx} + \tfrac{1}{x}$,
\pdd{meaning that}
$\psi(x+1) = \psi(x) + \tfrac{1}{x}.$
\pdd{To prove that the digamma function monotonically increases for $x>0$, we} show that
\begin{equation}\label{lemm1:digamma_decrease}
    \left(\Gamma'(x)\right)^2 < \Gamma(x)\Gamma''(x).
\end{equation}
\pdd{Indeed,}
\begin{equation*}
    \begin{array}{rl}
        \left(\Gamma'(x)\right)^2 &= \left(\int\limits_0^{+\infty}\pd{\exp(-t)}\ln t\cdot t^{x-1}dt\right)^2 \\
        &\overset{\circledOne}{<} \int\limits_0^{+\infty}\left(\pd{\exp(-\tfrac{t}{2})}t^{\tfrac{x-1}{2}}\right)^2dt\cdot\int\limits_0^{+\infty}\left(\pd{\exp(-\tfrac{t}{2})}t^{\tfrac{x-1}{2}}\ln t\right)^2dt\\
        &= \int\limits_0^{+\infty}\pd{\exp(-t)}t^{x-1}dt\cdot\int\limits_{0}^{+\infty}\pd{\exp(t)}t^{x-1}\ln^2 tdt = \Gamma(x)\Gamma''(x),
    \end{array}
\end{equation*}
where $\circledOne$ follows from \pdd{the} Cauchy-Schwartz inequality  \pdd{and we have strict inequality since} functions $\pd{\exp(-\tfrac{t}{2})}t^{\tfrac{x-1}{2}}$ and $\pd{\exp(-\tfrac{t}{2})}t^{\tfrac{x-1}{2}}\ln t$ are linearly independent. From \eqref{lemm1:digamma_decrease} \pd{it} follows that $\tfrac{d^2(\ln\Gamma(x))}{dx^2} = \left(\tfrac{\Gamma'(x)}{\Gamma(x)}\right)' = \tfrac{\Gamma''(x)}{\Gamma(x)} - \tfrac{\left(\Gamma'(x)\right)^2}{\left(\Gamma(x)\right)^2} \overset{\eqref{lemm1:digamma_decrease}}{>} 0,$ i.e.\ digamma function is increasing.

Now we show that $f_n(q)$ decreases on the interval $[2,+\infty)$. \pd{To that end, we} consider $\ln(f_{\pd{n}}(q))$
\begin{equation*}
    \begin{array}{rl}
        \ln(f_n(q))
        &= \ln\left(\tfrac{\Gamma(\tfrac{n}{2})}{\sqrt{\pi}}\right) + \ln\left(\Gamma\left(\tfrac{q+1}{2}\right)\right) - \ln\left(\Gamma\left(\tfrac{q+n}{2}\right)\right) - \tfrac{q}{2}\left(\ln(q-1)-\ln n\right),\\
        \tfrac{d(\ln(f_n(q)))}{dq} &= \tfrac{1}{2}\psi\left(\tfrac{q+1}{2}\right)-\tfrac{1}{2}\psi\left(\tfrac{q+n}{2}\right)-\tfrac{1}{2}\ln(q-1)-\tfrac{q}{2(q-1)} + \tfrac{1}{2}\ln n
    \end{array}
\end{equation*}
\pdd{and show that} $\tfrac{d(\ln(f_n(q)))}{dq} < 0$ for $q\geqslant 2$. Let $k = \lfloor\tfrac{n}{2}\rfloor$ (the \pdd{largest} integer which is no greater than $\tfrac{n}{2}$). Then $\psi\left(\tfrac{q+n}{2}\right) > \psi\left(k-1+\tfrac{q+1}{2}\right)$ and $\ln n \leqslant \ln(2k+1)$, whence,
\begin{equation*}
    \begin{array}{rl}
        \tfrac{d(\ln(f_n(q)))}{dq}
        &< \tfrac{1}{2}\left(\psi\left(\tfrac{q+1}{2}\right)-\psi\left(k-1+\tfrac{q+1}{2}\right)\right)-\tfrac{1}{2}\ln(q-1)-\tfrac{q}{2(q-1)} + \tfrac{1}{2}\ln(2k+1)\\
        &=\tfrac{1}{2} \left(\psi\left(\tfrac{q+1}{2}\right) - \sum\limits_{i=1}^{k-1}\tfrac{1}{\tfrac{q+1}{2}+ k - i - 1} - \psi\left(\tfrac{q+1}{2}\right)\right) - \tfrac{q}{2(q-1)} + \tfrac{1}{2}\ln\left(\tfrac{2k+1}{q-1}\right)\\
        &\overset{\circledOne}{\leqslant} -\tfrac{1}{2}\sum\limits_{i=1}^{k-1}\tfrac{2}{q-1+2k-2i} - \tfrac{1}{q-1} + \tfrac{1}{2}\ln\left(\tfrac{2k+1}{q-1}\right)\\
        &= -\tfrac{1}{2}\left(\tfrac{2}{q-1} + \tfrac{2}{q+1} + \tfrac{2}{q+3}+ \ldots + \tfrac{2}{q+2k-3}\right) + \tfrac{1}{2}\ln\left(\tfrac{2k+1}{q-1}\right)\\
        &\overset{\circledTwo}{<} -\tfrac{1}{2}\ln\left(\tfrac{q+2k-1}{q-1}\right)+\tfrac{1}{2}\ln\left(\tfrac{2k+1}{q-1}\right) \overset{\circledThree}{\leqslant}-\tfrac{1}{2}\ln\left(\tfrac{2k+1}{q-1}\right)+\tfrac{1}{2}\ln\left(\tfrac{2k+1}{q-1}\right)= 0,
    \end{array}
\end{equation*}
where $\circledOne$ and $\circledThree$ \pdd{are since} $q\geqslant2$, $\circledTwo$ \pdd{follows from an estimate} of \pdd{the} integral of $\tfrac{1}{x}$ by \pdd{the} integral of \newstuff{the constant functions} $g_i(x) = \tfrac{1}{q-1+2i},\, x\in[q-1+2i,q-1+2i+2],\, i=\pd{0,...,2k-1}$: 
 $\tfrac{2}{q-1} + \tfrac{2}{q+1} + \tfrac{2}{q+3} + \ldots + \tfrac{2}{q+2k-3} > \int\limits_{q-1}^{q+2k-1}\tfrac{1}{x}dx = \ln\left(\tfrac{q+2k-1}{q-1}\right).$

\pdd{Thus}, we \pd{have} shown that $\tfrac{d(\ln(f_n(q)))}{dq} < 0$ for $q\geqslant2$ \pd{and an} arbitrary natural number $n$. Therefore, for any fixed number $n$, the function $f_n(q)$ decreases as $q$ increases, which means that $f_n(q)\leqslant f_n(2) = 0$, i.e., \eqref{lemm1:key_estimation} holds. From this, \eqref{lemm1:jensen}, and \eqref{lemm1:expectation_component} we obtain \pdd{\eqref{lemm1:rough_estimation_of_expectation_q_norm}, i.e. that,} for all $2\leqslant q < \infty$,
\begin{equation}\label{lemm1:pre-final}
    \EE[||e||_q^2] \overset{\eqref{lemm1:jensen}}{\leqslant} \left(n\EE[|e_2|^q]\right)^{\tfrac{2}{q}} \overset{\eqref{lemm1:expectation_component},\eqref{lemm1:key_estimation}}{\leqslant} (q-1)n^{\tfrac{2}{q}-1}.
\end{equation}

\pdd{Next, we analyze separately the case of large $q$, in particular, $q=\infty$.}
\pdd{We} consider \pdd{the r.h.s.} of \eqref{lemm1:pre-final} as \pdd{a} function of $q$ and find its minimum for $q\geqslant2$. \pdd{Denote} $h_n(q) = \ln(q-1) + \left(\tfrac{2}{q}-1\right)\ln n$, \pdd{which is the} logarithm of the r.h.s. of \eqref{lemm1:pre-final}. \pdd{The} derivative of $h_n(q)$ is $\tfrac{dh_n(q)}{dq} = \tfrac{1}{q-1} -\tfrac{2\ln n}{q^2}$, which implies \pdd{that the first-order optimality condition is} $\tfrac{1}{q-1} -\tfrac{2\ln n}{q^2} = 0,$ \pdd{or equivalently} $ q^2-2q\ln n + 2\ln n = 0$.
If $n \geqslant 8$, then \pdd{the function $h_n(q)$} attains its minimum on the set $[2,+\infty)$ at $q_0 = \ln n\left(1+\sqrt{1-\tfrac{2}{\ln n}}\right)$ (for the case $n\leqslant 7$ \pdd{the optimal point is} $q_0 = 2$ and without loss of generality, we assume $n\geqslant8$). Therefore, for all $q>q_0$, \pdd{including $q=\infty$, we have}
\begin{equation}\label{lemm1:pre_final_big_q}
    \begin{array}{rl}
        \EE[||e||_q^2] &\overset{\circledOne}{<} \EE[||e||_{q_0}^2] \overset{\eqref{lemm1:pre-final}}{\leqslant}(q_0-1)n^{\tfrac{2}{q_0}-1}\overset{\circledTwo}{\leqslant}(2\ln n -1)n^{\tfrac{2}{\ln n}-1}\\
        &= (2\ln n -1)\tfrac{\pd{\exp(2)}}{n}\leqslant(16\ln n -8)\tfrac{1}{n} \leqslant(16\ln n -8)n^{\tfrac{2}{q}-1},
    \end{array}
\end{equation}
where $\circledOne$ is since $\|e\|_q < \|e\|_{q_0}$ for $q > q_0$, $\circledTwo$ follows from $q_0 \leqslant 2\ln n,\, q_0 \geqslant \ln n$. \pdd{Combining estimates} \eqref{lemm1:pre-final} and \eqref{lemm1:pre_final_big_q}, we obtain \eqref{lemm1:expect_q_norm}. 

\pdd{It remains to} prove \eqref{lemm1:expect_inner_product}. First, we estimate $\sqrt{\EE[\|e\|_q^4]}$. \pdd{By the} probabilistic Jensen's inequality, for $q\geqslant2$,
\begin{equation*}
    \begin{array}{rl}
        \EE[||e||_q^4] &= \EE\left[\left(\left(\sum\limits_{k=1}^{n}|e_k|^q\right)^2\right)^{\tfrac{2}{q}}\right] \leqslant \left(\EE\left[\left(\sum\limits_{k=1}^{n}|e_k|^q\right)^2\right]\right)^{\tfrac{2}{q}}\\
        &\overset{\circledOne}{\leqslant} \left(\EE\left[\left(n\sum\limits_{k=1}^{n}|e_k|^{2q}\right)\right]\right)^{\tfrac{2}{q}} \overset{\circledTwo}{=} \left(n^2\EE[|e_2|^{2q}]\right)^{\tfrac{2}{q}}\\
        &\overset{\eqref{lemm1:expectation_component},\eqref{lemm1:key_estimation}}{\leqslant} n^{\tfrac{4}{q}}\left(\left(\tfrac{2q-1}{n}\right)^{\tfrac{2q}{2}}\right)^{\tfrac{2}{q}} = (2q-1)^{2}n^{\tfrac{4}{q}-2},
    \end{array}
\end{equation*}
where $\circledOne$ is since $\left(\sum\limits_{k=1}^{n} x_k\right)^2 \leqslant n\sum\limits_{k=1}^{n} x_k^2$ for $x_1,x_2,\ldots,x_n\in\R$ and $\circledTwo$ follows \pdd{from the linearity of expectation and the components of the random vector $e$ being} identically distributed. From this we obtain
\begin{equation}\label{lemm1:pre_final_q_norm_power4}
    \sqrt{\EE[||e||_q^4]} \leqslant (2q-1)n^{\tfrac{2}{q}-1}.
\end{equation}
Next, we consider the r.h.s. of  \eqref{lemm1:pre_final_q_norm_power4} as a function of $q$ and find its minimum for $q\geqslant2$. \pdd{The logarithm of the r.h.s. of \eqref{lemm1:pre_final_q_norm_power4} is $h_n(q) = \ln(2q-1) + \left(\tfrac{2}{q}-1\right)\ln n$ with the derivative $\tfrac{dh_n(q)}{dq} = \tfrac{2}{2q-1} -\tfrac{2\ln n}{q^2}$, which implies the first-order optimality condition $\tfrac{2}{2q-1} -\tfrac{2\ln n}{q^2} = 0$, or equivalently $q^2-2q\ln n + \ln n = 0$.} If $n \geqslant 3$, the point where the function \pdd{$h_n(q)$ attains} its minimum on the set $[2,+\infty)$ is $q_0 = \ln n\left(1+\sqrt{1-\tfrac{1}{\ln n}}\right)$ (for the case $n\leqslant 2$ \pdd{the optimal point is} $q_0 = 2$ and without loss of generality we assume that $n\geqslant3$). Therefore for all $q>q_0$\pdd{, including $q=\infty$,}
\begin{equation}\label{lemm1:pre_final_big_q_power4}
    \begin{array}{rl}
        \sqrt{\EE[\|e\|_q^4]} &\overset{\circledOne}{<} \sqrt{\EE[\|e\|_{q_0}^4]} \overset{\eqref{lemm1:pre_final_q_norm_power4}}{\leqslant}(2q_0-1)n^{\tfrac{2}{q_0}-1}\overset{\circledTwo}{\leqslant}(4\ln n -1)n^{\tfrac{2}{\ln n}-1}\\
        &= (4\ln n -1)\tfrac{\pd{\exp(2)}}{n}\leqslant(32\ln n -8)\tfrac{1}{n} \leqslant(32\ln n -8)n^{\tfrac{2}{q}-1},
    \end{array}
\end{equation}
where $\circledOne$ is since $\|e\|_q < \|e\|_{q_0}$ for $q > q_0$, $\circledTwo$ follows from $q_0 \leqslant 2\ln n,\, q_0 \geqslant \ln n$. \pdd{Combining the estimates} \eqref{lemm1:pre_final_q_norm_power4} and \eqref{lemm1:pre_final_big_q_power4}, we get the inequality
\begin{equation}\label{lemm1:expect_q_norm_power4}
    \sqrt{\EE[\|e\|_q^4]}\leqslant \min\{2q-1,32\ln n -8\}n^{\tfrac{2}{q}-1}.
\end{equation}

\pdd{The next step is to estimate} $\EE[\langle s,\,e\rangle^4]$, where $s\in\R^n$ is some fixed vector. Let $S_n(r)$ be \pdd{the} surface area of $n$-dimensional Euclidean sphere with radius $r$ and $d\sigma(e)$ be unnormalized uniform measure on $n$-dimensional Euclidean sphere. \pdd{Then} $S_n(r) = S_n(1)r^{n-1},\, \tfrac{S_{n-1}(1)}{S_n(1)} = \tfrac{n-1}{n\sqrt{\pi}}\tfrac{\Gamma(\tfrac{n+2}{2})}{\Gamma(\tfrac{n+1}{2})}$. Let $\varphi$ be the angle between $s$ and $e$.
Then
\begin{equation}\label{lemm1:inner_product_power4}
    \begin{array}{rl}
        &\EE[\langle s,\, e\rangle^4] = \tfrac{1}{S_n(1)}\int\limits_{S}\langle s,\, e\rangle^4d\sigma(\varphi) = \tfrac{1}{S_n(1)}\int\limits_0^\pi\|s\|_2^4\cos^3\varphi S_{n-1}(\sin\varphi)d\varphi\\
        &= \|s\|_2^4\tfrac{S_{n-1}(1)}{S_n(1)}\int\limits_{0}^\pi\cos^4\varphi\sin^{n-2}\varphi d\varphi = \|s\|_2^4\cdot\tfrac{n-1}{n\sqrt{\pi}}\tfrac{\Gamma(\tfrac{n+2}{2})}{\Gamma(\tfrac{n+1}{2})}\int\limits_{0}^\pi\cos^4\varphi\sin^{n-2}\varphi d\varphi.
    \end{array}
\end{equation}
Further, denoting the Beta function by $B(\cdot,\cdot)$,
\begin{equation*}
    \begin{array}{rl}
       & \int\limits_0^\pi\cos^4\varphi\sin^{n-2}\varphi d\varphi = 2\int\limits_0^{\tfrac{\pi}{2}}\cos^4\varphi\sin^{n-2}\varphi d\varphi \overset{\newstuff{t=\sin^2\varphi}}{=} \int\limits_0^{\tfrac{\pi}{2}}t^{\tfrac{n-3}{2}}(1-t)^{\tfrac{3}{2}}dt  \\
        &=  B(\tfrac{n-1}{2},\tfrac{5}{2})= \tfrac{\Gamma(\tfrac{5}{2})\Gamma(\tfrac{n-1}{2})}{\Gamma(\tfrac{n+4}{2})} = \tfrac{\tfrac{3}{2}\cdot\tfrac{1}{2}\Gamma(\tfrac{1}{2})\Gamma(\tfrac{n-1}{2})}{\tfrac{n+2}{2}\cdot\Gamma(\tfrac{n+2}{2})} = \tfrac{3}{n+2}\cdot\tfrac{\sqrt{\pi}\Gamma(\tfrac{n-1}{2})}{2\Gamma(\tfrac{n+2}{2})}.
    \end{array}
\end{equation*}
From this and \eqref{lemm1:inner_product_power4}, we obtain
\begin{equation}\label{lemm1:inner_product_power4_final}
    \begin{array}{rl}
        \EE[\langle s,\, e\rangle^4] &= \|s\|_2^4\cdot\tfrac{n-1}{n\sqrt{\pi}}\tfrac{\Gamma(\tfrac{n+2}{2})}{\Gamma(\tfrac{n+1}{2})}\cdot\tfrac{3}{n+2}\cdot\tfrac{\sqrt{\pi}\Gamma(\tfrac{n-1}{2})}{2\Gamma(\tfrac{n+2}{2})}\\
        &= \|s\|_2^4\cdot\tfrac{3(n-1)}{2n(n+2)}\cdot\tfrac{\Gamma(\tfrac{n-1}{2})}{\tfrac{n-1}{2}\Gamma(\tfrac{n-1}{2})}
        = \tfrac{3\|s\|_2^4}{n(n+2)} \overset{\circledOne}{\leqslant} \tfrac{3\|s\|_2^4}{n^2}.
    \end{array}
\end{equation}

To prove \eqref{lemm1:expect_inner_product}, it remains to use \eqref{lemm1:expect_q_norm_power4}, \eqref{lemm1:inner_product_power4_final} and the Cauchy-Schwartz inequality $(\EE[XY])^2\leqslant\EE[X^2]\cdot\EE[Y^2]$:
\begin{equation*}
    \begin{array}{c}
        \EE[\langle s,\,e\rangle^2 ||e||_q^2] 
        \leqslant \sqrt{\EE[\langle s,\,e\rangle^4]\cdot\EE[\|e\|_q^4]} \leqslant \sqrt{3}\|s\|_2^2\min\{2q-1,32\ln n -8\}n^{\tfrac{2}{q}-2}.
    \end{array}
\end{equation*}

\section{Technical Results on Recurrent Sequences}
\label{S:Tech_result}
\begin{lemma}\label{stoh:technical_lemma}
    Let $a_0,\ldots,a_{N-1}, b, R_1,\ldots, R_{N-1}$ be non-negative numbers and
    \begin{equation}\label{technical_lemma_assumption}
        R_{l} \leqslant \sqrt{2}\cdot\sqrt{\left(\sum\limits_{k=0}^{l-1}a_k + b\sum\limits_{k=1}^{l-1}\alpha_{k+1}R_k \right)}\quad l=1,\ldots,N,
    \end{equation}
    where $\alpha_{k+1} = \tfrac{k+2}{96n^2\rho_nL_2}$ for all $k\in\NN$. Then, for $l=1,\ldots,N$,
    \begin{equation}\label{technical_inequality_for_induction}
        \sum\limits_{k=0}^{l-1}a_k + b\sum\limits_{k=1}^{l-1} \alpha_{k+1}R_k \leqslant \left(\sqrt{\sum\limits_{k=0}^{l-1}a_k} + \sqrt{2}b\cdot\tfrac{l^2}{96n^2\rho_nL_2}\right)^2.
    \end{equation}
\end{lemma}
\begin{proof}
For $l=1$ \pdd{the inequality is trivial}. \pdd{Next we assume that \eqref{technical_inequality_for_induction} holds for some $l < N$ and prove this inequality} for $l+1$. From the induction assumption and \eqref{technical_lemma_assumption} we obtain
\begin{equation}\label{technical_distance_estimation_3}
    \begin{array}{c}
        R_l \leqslant \sqrt{2}\left(\sqrt{\sum\limits_{k=0}^{l-1}a_k} + \sqrt{2}b\cdot\tfrac{l^2}{96n^2\rho_nL_2}\right),
    \end{array}
\end{equation}
whence
\begin{equation*}
    \begin{array}{rl}
        &\sum\limits_{k=0}^{l}a_k + b\sum\limits_{k=1}^{l}\alpha_{k+1}R_k = \sum\limits_{k=0}^{l-1}a_k + b\sum\limits_{k=1}^{l-1}\alpha_{k+1}R_k + a_l + b\alpha_{l+1}R_{l}\\
        &\overset{\circledOne}{\leqslant} \left(\sqrt{\sum\limits_{k=0}^{l-1}a_k} + \sqrt{2}b\cdot\tfrac{l^2}{96n^2\rho_nL_2}\right)^2 + a_l+ \sqrt{2}b\alpha_{l+1}\left(\sqrt{\sum\limits_{k=0}^{l-1}a_k} + \sqrt{2}b\cdot\tfrac{l^2}{96n^2\rho_nL_2}\right)\\
        &= \sum\limits_{k=0}^{l}a_k + 2\sqrt{\sum\limits_{k=0}^{l-1}a_k}\cdot \tfrac{l^2 \cdot \sqrt{2}b}{96n^2\rho_nL_2} + \tfrac{l^4 \cdot 2b^2}{(96n^2\rho_nL_2)^2}+ \sqrt{2}b\alpha_{l+1}\left(\sqrt{\sum\limits_{k=0}^{l-1}a_k} + \tfrac{l^2 \cdot \sqrt{2}b}{96n^2\rho_nL_2}\right)\\
        &= \sum\limits_{k=0}^{l}a_k + 2\sqrt{\sum\limits_{k=0}^{l-1}a_k}\cdot\sqrt{2}b\left(\tfrac{l^2}{96n^2\rho_nL_2} + \tfrac{\alpha_{l+1}}{2}\right)+ 2b^2\left(\tfrac{l^4}{(96n^2\rho_nL_2)^2} + \cdot\tfrac{\alpha_{l+1}l^2}{96n^2\rho_nL_2}\right)\\
        &\overset{\circledTwo}{\leqslant} \sum\limits_{k=0}^{l}a_k + 2\sqrt{\sum\limits_{k=0}^{l}a_k}\tfrac{(l+1)^2 \cdot \sqrt{2}b}{96n^2\rho_nL_2} + \tfrac{(l+1)^4 \cdot 2b^2}{(96n^2\rho_nL_2)^2}= \left(\sqrt{\sum\limits_{k=0}^{l}a_k} + \sqrt{2}b\cdot\tfrac{(l+1)^2}{96n^2\rho_nL_2}\right)^2,
    \end{array}    
\end{equation*}
where $\circledOne$ holds by the induction assumption and \eqref{technical_distance_estimation_3}, $\circledTwo$ is since $\sum\limits_{k=0}^{l-1}a_k\leqslant \sum\limits_{k=0}^{l}a_k$ and
\begin{equation*}
    \begin{array}{rl}
        \tfrac{l^2}{96n^2\rho_nL_2} + \tfrac{\alpha_{l+1}}{2} = \tfrac{2l^2 + l + 2}{192n^2\rho_nL_2} &\leqslant \tfrac{(l+1)^2}{96n^2\rho_nL_2},\\
        \tfrac{l^4}{(96n^2\rho_nL_2)^2} + \alpha_{l+1}\cdot\tfrac{l^2}{96n^2\rho_nL_2} &\leqslant \tfrac{l^4 + (l+2)l^2}{(96n^2\rho_nL_2)^2} \leqslant \tfrac{(l+1)^4}{(96n^2\rho_nL_2)^2}.
    \end{array}
\end{equation*}
\end{proof}

\begin{lemma}\label{stoh:technical_lemma_non_acc}
    Let $\pd{\alpha},a_0,\ldots,a_{N-1}, b, R_1,\ldots, R_{N-1}$ be non-negative numbers and
    \begin{equation}\label{technical_lemma_assumption_non_acc}
        R_{l} \leqslant \sqrt{2}\cdot\sqrt{\left(\sum\limits_{k=0}^{l-1}a_k + b\alpha\sum\limits_{k=1}^{l-1}R_k \right)}\quad l=1,\ldots,N.
    \end{equation}
    Then, for $l=1,\ldots,N$,
    \begin{equation}\label{technical_inequality_for_induction_non_acc}
        \sum\limits_{k=0}^{l-1}a_k + b\alpha\sum\limits_{k=1}^{l-1} R_k \leqslant \left(\sqrt{\sum\limits_{k=0}^{l-1}a_k} + \sqrt{2}b\alpha l\right)^2.
    \end{equation}
\end{lemma}
\begin{proof}
For $l=1$  \pdd{the inequality is trivial}. Next we assume that \eqref{technical_inequality_for_induction_non_acc} holds for some $l < N$ and prove it for $l+1$. By the induction assumption and \eqref{technical_lemma_assumption_non_acc} we obtain
\begin{equation}\label{technical_distance_estimation_3_non_acc}
    \begin{array}{c}
        R_l \leqslant \sqrt{2}\left(\sqrt{\sum\limits_{k=0}^{l-1}a_k} + \sqrt{2}b\alpha l\right),
    \end{array}
\end{equation}
whence
\begin{equation*}
    \begin{array}{rl}
        &\sum\limits_{k=0}^{l}a_k + b\alpha\sum\limits_{k=1}^{l}R_k = \sum\limits_{k=0}^{l-1}a_k + b\alpha\sum\limits_{k=1}^{l-1}R_k + a_l + b\alpha R_{l}\\
        &\overset{\circledOne}{\leqslant} \left(\sqrt{\sum\limits_{k=0}^{l-1}a_k} + \sqrt{2}b\alpha l\right)^2 + a_l
        + \sqrt{2}b\alpha\left(\sqrt{\sum\limits_{k=0}^{l-1}a_k} + \sqrt{2}b\alpha l\right)\\
        &= \sum\limits_{k=0}^{l}a_k + 2\sqrt{\sum\limits_{k=0}^{l-1}a_k}\cdot \sqrt{2}b\alpha l + 2b^2\alpha^2 l^2+ \sqrt{2}b\alpha\left(\sqrt{\sum\limits_{k=0}^{l-1}a_k} + \sqrt{2}b\alpha l\right)\\
        &= \sum\limits_{k=0}^{l}a_k + 2\sqrt{\sum\limits_{k=0}^{l-1}a_k}\cdot\sqrt{2}b\alpha\left(l + \tfrac{1}{2}\right) + 2b^2\alpha^2\left(l^2 + l\right)\\
        &
        \hspace{-0.8em}\overset{\circledTwo}{\leqslant} \sum\limits_{k=0}^{l}a_k + 2\sqrt{\sum\limits_{k=0}^{l}a_k}\cdot \sqrt{2}b\alpha (l+1) + 2(b\alpha (l+1))^2= \left(\sqrt{\sum\limits_{k=0}^{l}a_k} + \sqrt{2}b\alpha (l+1)\right)^2,
    \end{array}    
\end{equation*}
where $\circledOne$ is by the induction assumption and \eqref{technical_distance_estimation_3_non_acc}, $\circledTwo$ is since $\sum\limits_{k=0}^{l-1}a_k\leqslant \sum\limits_{k=0}^{l}a_k$.
\end{proof}
\end{document}